\documentclass[moor]{informs1}              % for a regular run

\usepackage{natbib}
 \NatBibNumeric
 \def\bibfont{\small}%
 \bibpunct[, ]{[}{]}{,}{n}{}{,}%
 
 \usepackage{graphicx}
 \usepackage{url,color}
 \usepackage{subfigure}
 \usepackage{amsfonts,mathrsfs}
 \usepackage{amssymb,amsmath}
 \usepackage{verbatim}
 \usepackage{acronym}
 \usepackage{multirow}
 \usepackage{array}
 \usepackage{paralist}
 \usepackage[utf8]{inputenc}
 \usepackage[english]{babel}
 \usepackage[mathscr]{euscript}
 \usepackage{bbold}
 \usepackage{bbm}
 \usepackage{accents}
 
\newcommand{\ubar}[1]{\underaccent{\bar}{#1}}
\newcommand{\norms}[1]{\Vert#1\Vert}
\newcommand{\iprods}[1]{\langle #1\rangle}
\newcommand{\R}{\mathbb{R}}
\newcommand{\Pc}{\mathcal{P}}
\newcommand{\Lc}{\mathcal{L}}

\usepackage{tkz-tab}
\usepackage{xpatch}
\usepackage{enumerate}   
\usepackage{arydshln}
\usepackage{slashbox}
\usepackage{adjustbox}

\usepackage{algorithm}
\usepackage{algpseudocode}

\newcommand{\Tc}{\mathcal{T}}

\newcommand{\Eproof}{\hfill $\square$}
 
 \DeclareUnicodeCharacter{2212}{-}

 \DeclareMathOperator{\dom}{dom}
  \DeclareMathOperator{\epi}{epi}

\usepackage[colorlinks=true,breaklinks=true,bookmarks=true,urlcolor=blue,
     citecolor=blue,linkcolor=blue,bookmarksopen=false,draft=false]{hyperref}

\def\EMAIL#1{\href{mailto:#1}{#1}}% When hyperref is used, otherwise outcomment 
\TheoremsNumberedThrough     % Preferred (Theorem 1, Lemma 1, Theorem 2)
\EquationsNumberedThrough    % Default: (1), (2), ...
\begin{document}
%%%%%%%%%%%%%%%%

\TITLE{Complexity of Linearized Perturbed Augmented Lagrangian Methods for  Nonsmooth Nonconvex Optimization with Nonlinear Equality Constraints}

\ARTICLEAUTHORS{%
\AUTHOR{Lahcen El Bourkhissi}
\AFF{Automatic Control and Systems
	Engineering Department, Politehnica Bucharest, 060042
	Bucharest, Romania, \EMAIL{lahcenelbourkhissi1997@gmail.com} }
    
\AUTHOR{Ion Necoara}
\AFF{Automatic Control and Systems
	Engineering Department, Politehnica Bucharest, 060042
	Bucharest, Romania  and  Gheorghe Mihoc-Caius Iacob  Institute of Mathematical Statistics and Applied Mathematics of the Romanian Academy, 050711 Bucharest, Romania, \EMAIL{ion.necoara@upb.ro}}

\AUTHOR{Panagiotis Patrinos}
\AFF{Department of  Electrical Engineering, K.U. Leuven, B-3001 Leuven, Belgium, \EMAIL{panos.patrinos@esat.kuleuven.be} }

\AUTHOR{Quoc Tran-Dinh}
\AFF{Department of Statistics and Operations Research, University of North Carolina at Chapel Hill, NC 27599 Chapel Hill, USA \EMAIL{quoctd@email.unc.edu} }
} 

\ABSTRACT{
This paper addresses a class of general nonsmooth and nonconvex composite optimization problems subject to nonlinear equality constraints. 
We assume that a part of the objective function and the functional constraints exhibit local smoothness. 
To tackle this challenging class of problems, we propose a novel linearized perturbed augmented Lagrangian method. 
This method incorporates a perturbation in the augmented Lagrangian function by scaling the dual variable with a sub-unitary parameter. 
Furthermore, we linearize the smooth components of the objective and the constraints within the perturbed Lagrangian function at the current iterate, while preserving the nonsmooth components. 
This approach, inspired by prox-linear (or Gauss-Newton) methods, results in a convex subproblem that is typically easy to solve. 
The solution of this subproblem then serves as the next primal iterate, followed by a perturbed ascent step to update the dual variables. Under a newly introduced constraint qualification condition, we establish the boundedness of the dual iterates. 
We derive convergence guarantees for the primal iterates, proving convergence to an $\epsilon$-first-order optimal solution within $\mathcal{O}(\epsilon^{-3})$ evaluations of the problem's functions and their first derivatives. Moreover, when the problem exhibits for example a  semialgebraic property, we derive improved local convergence results. Finally, we validate the theoretical findings and assess the practical performance of our proposed algorithm through numerical comparisons with existing state-of-the-art methods.
}


\MSCCLASS{90C25, 90C06, 65K05.}

\HISTORY{2024}

\maketitle

%%%%%%%%%%%%%%%%%%%%%%%%%%%%%%%%%%%%%%%%%%%%%
\section{Introduction.}
\label{intro}
Numerous applications across diverse fields, including statistics, control and signal processing, and machine learning, can be formulated as nonsmooth nonconvex optimization problems subject to nonlinear equality constraints (see, e.g., \cite{HonHaj:17,LukSab:19,NedNec:14}). 
This paper addresses this class of challenging problems by employing a perturbed augmented Lagrangian approach. 
Originally introduced in \cite{Hes:69,Pow:69} to minimize functions subject to linear equality constraints, the augmented Lagrangian method, also known as the method of multipliers, has proven to be a powerful tool to develop scalable optimization algorithms. 
Even in the context of nonconvex problems, it offers several theoretical advantages, such as a zero duality gap and an exact penalty representation (see \cite{RocWet:98}). 
Furthermore, the augmented Lagrangian framework serves as the foundation for the Alternating Direction Method of Multipliers (ADMM), a highly efficient method for solving large-scale (and distributed) optimization problems, as demonstrated in, e.g., \cite{BoyPar:11,CohHal:21,ElbNec:23,GloTal:89}.

%%% Related work.
\medskip 
\noindent \textit{\textbf{Related work.}} 
While the augmented Lagrangian method has been extensively applied to convex problems, see \cite{Ber:96, BoyPar:11, GloTal:89, NedNec:14, PatNec:17} and non-convex problems with linear equality constraints, see \cite{HajHon:19, HonHaj:17, JiaLin:19, ZhaLuo:20}, its application to non-convex optimization problems with nonlinear equality constraints has been relatively limited (see, e.g., \cite{CohHal:21, ElbNec:23, HalTeb:23, Lu:22, SahEft:19, XieWri:21} as a few examples).
For instance, Hajinezhad and Hong  \cite{HajHon:19} and Lu \cite{Lu:22} proposed  perturbed variants of the augmented Lagrangian method. 
This approach introduces a geometric perturbation of the dual variables, controlled by a coefficient $1-\tau \in (0, 1)$. 
Notably, these perturbed methods can be viewed as an interpolation between the standard augmented Lagrangian method (when $1-\tau \to 1$) and the quadratic penalty method (when $1-\tau \to 0$).
More specifically,  Hajinezhad and Hong in \cite{HajHon:19}  introduced a perturbed augmented Lagrangian method for solving linearly constrained optimization problems. 
These problems involve a composite objective function consisting of a smooth (potentially nonconvex) component and a convex (nonsmooth but Lipschitz) component, along with simple constraints. Their method employed a primal metric proximal gradient step followed by a perturbed dual gradient ascent step. 
Under appropriate conditions, including a penalty parameter of order $\mathcal{O}(\epsilon^{-1})$, the authors demonstrated that this method achieves an $\epsilon$-first-order optimal solution with a complexity of $\mathcal{O}(\epsilon^{-4})$ for a free perturbation parameter $1 - \tau$.  Furthermore, for problems involving nonconvex inequality constraints, a smooth objective function, and simple set constraints, Lu in \cite{Lu:22} proposed a perturbed augmented Lagrangian method. 
This method linearizes the nonlinear constraints within the perturbed augmented Lagrangian function during the primal update, followed by a perturbed dual update. The author established a convergence rate of  $\mathcal{O}(\epsilon^{-3})$ for finding an $\epsilon$-first-order solution when the penalty parameter is set to $\mathcal{O}(\epsilon^{-1})$ and the perturbation parameter $1-\tau$ is below a small threshold. 
This parameter choice makes the method resemble more closely a quadratic penalty approach than a traditional augmented Lagrangian approach. However, the definition of an $\epsilon$-first-order solution presented by the author is questionable; while feasibility is ensured at the proposed point, the satisfaction of the optimality condition remains unclear. 
This ambiguity arises particularly from the vagueness and lack of clear justification for inequality (263) in \cite{Lu:22}.

\medskip 

In \cite{HalTeb:23}, an adaptive augmented Lagrangian method was employed to tackle a class of nonconvex optimization problems characterized by nonsmoothness and a nonlinear functional composite structure in the objective function. 
To address this structure, a slack variable was introduced, transforming the original problem into one with a nonsmooth objective and nonlinear equality constraints. This reformulation involves two sets of variables: the original variables and the newly introduced slack variables. 
Notably, when considering these two variable sets collectively, their problem can be viewed as a specific instance of the problem addressed in this paper.
The authors in \cite{HalTeb:23} then applied an ADMM scheme to solve the reformulated problem. 
This approach involves linearizing the smooth component of the augmented Lagrangian function at each iteration, resulting in proximal updates for each step of the algorithm. 
Under the assumption of bounded dual multipliers -- specifically, assuming that there exists a constant $M > 0$ such that $\Vert y^k \Vert \leq M$ for all iterations $k$ -- their method achieves an $\epsilon$-first-order optimal solution with a complexity of $\mathcal{O}(\epsilon^{-2})$.
However, to ensure that the iterates remain $\epsilon$-feasible after a certain number of iterations, it is crucial that the penalty parameter $\rho_k$ satisfies the condition $\rho_k \geq \frac{2M}{\epsilon}$. 
Consequently, based on the proof of Theorem 2 in \cite{HalTeb:23}, the overall complexity of their method is ultimately $\mathcal{O}(\epsilon^{-4})$. Furthermore, Sahin et al. \cite{SahEft:19} addressed the inexact solution of the augmented Lagrangian subproblem and incorporated a decreasing stepsize for updating the dual multipliers to ensure the boundedness of these multipliers. 
Their analysis established a total complexity of $\mathcal{O}(\epsilon^{-4})$ for their inexact augmented Lagrangian method (iALM) to achieve an $\epsilon$-first-order solution of the original problem, provided that the penalty parameter is scaled as $\mathcal{O}(\epsilon^{-1})$.
Additionally, ADMM methods have been specifically designed for separable problems with structured nonlinear equality constraints in \cite{CohHal:21,ElbNec:23}.

%%% Drawback of bounding multipliers.
\medskip 
\noindent\textit{\textbf{Drawback of bounding multipliers.}}
A primary challenge in employing (perturbed) augmented Lagrangian methods lies in simultaneously ensuring feasibility and satisfying the optimality condition at a test point. 
A common approach to address this challenge involves assuming boundedness of the dual iterates and progressively increasing the penalty parameter, as exemplified in \cite{HalTeb:23}. 
However, this boundedness assumption presents a significant limitation, as it is imposed on the algorithm's generated sequence rather than being an inherent property of the problem itself. 
Indeed, Hallak and Teboulle in \cite{HalTeb:23} acknowledged the difficulty of ensuring boundedness of the multiplier sequence in nonconvex settings by stating that: ``the boundedness of the multiplier sequence $\{y^k\}$ in the nonconvex setting is a very difficult matter and not at all obvious because coercivity arguments do not apply directly, and we are not aware of any breakthrough in this area.''

To circumvent the restrictive assumption of bounded multipliers, the authors in \cite{SahEft:19} employed an augmented Lagrangian algorithm with a sufficiently decreasing stepsize for the dual iterates. 
This approach allows them to control the growth of the dual multipliers. 
Furthermore, they imposed a regularity condition to facilitate the control of the feasibility measure. However, this strategy necessitates short steps in the dual updates, which consequently slows down the overall algorithm. 
Moreover, the regularity condition employed in \cite{SahEft:19} is not entirely accurate, as discussed in Section \ref{sec_reg} in the sequel, where we present our own regularity condition. 
As a result, the algorithms proposed in both \cite{HalTeb:23} and \cite{SahEft:19} exhibit a relatively high computational complexity, specifically of order $\mathcal{O}(\epsilon^{-4})$.

In \cite{Lu:22}, the author derived  $\epsilon$-feasibility but encountered challenges in establishing the boundedness of the dual variables. 
Furthermore, the method in \cite{Lu:22} necessitates a very small perturbation parameter, effectively rendering it similar to a quadratic penalty method. 
In contrast, \cite{HajHon:19} successfully demonstrated the boundedness of dual variables under the Robinson constraint qualification condition. 
However, this result is limited to a specific class of optimization problems: those involving only linear equality constraints, simple constraints, and a nonsmooth objective function that must be convex and Lipschitz continuous. 
Even within this restricted setting, the method exhibits a high computational complexity, specifically $\mathcal{O}(\epsilon^{-4})$.

%%%% Contribution.
\medskip 
\noindent \textit{\textbf{Our contribution.}}  
This paper proposes a novel Linearized Perturbed Augmented Lagrangian method (abbreviated by LIPAL) specifically designed to address the challenges encountered in solving general nonsmooth, nonconvex optimization problems subject to nonlinear equality constraints. 
To ensure the boundedness of dual iterates, we introduce a novel regularity condition that generalizes the well-known LICQ condition in the literature. 
This novel regularity condition, in conjunction with the perturbation of the augmented Lagrangian, enables us to establish rigorous convergence guarantees for the iterates generated by our proposed LIPAL method under various assumptions on the problem structures. 
Our primary contributions can be summarized as follows.
\begin{compactenum}[(i)]
\item 
We focus on a general optimization problem characterized by a composite and possibly nonsmooth and nonconvex objective and nonlinear functional equality constraints. 
This model covers a broad range of problems, including general composite formulations as well as those involving nonconvex inequality constraints. 
We propose a novel constraint qualification condition that extends existing ones in the literature, such as LICQ, and possibly covers new models. 
This new condition is instrumental in bounding the dual iterates, effectively addressing the issue highlighted in \cite{HalTeb:23} as discussed above.

\item 
We develop a novel perturbed augmented Lagrangian method. 
At each iteration, we linearize both the smooth component of the objective function and the functional constraints within the perturbed augmented Lagrangian framework. This linearization is performed using a prox-linear-type (also known as Gauss-Newton-type) mechanism. 
Furthermore, we incorporate a suitable regularization term into the method. 
Our algorithm exhibits several desirable properties. 
Notably, it only necessitates evaluations of the problem's function values and their first-order derivatives. 
Moreover, each iteration involves solving a \textit{simple} subproblem, which is guaranteed to be  convex provided that the nonsmooth component of the objective function is convex. 
This characteristic enables the efficient handling of large-scale nonconvex problems. 
The solution obtained from this subproblem is then utilized to update the dual variables through a new perturbed ascent step.

\item 
To establish the boundedness of the dual iterates, we combine the newly introduced constraint qualification condition with our perturbed ascent step on the dual variables. 
To the best of our knowledge, this represents one of the first results of this nature obtained within this general framework. 
We rigorously establish convergence guarantees for the iterates generated by our method. 
Specifically, we prove global convergence to an $\epsilon$-first-order optimal solution within $\mathcal{O}(\epsilon^{-3})$ evaluations of the problem's functions and their first derivatives, improved over existing results by a factor of $\epsilon^{-1}$. 
Furthermore, by leveraging the Kurdyka-{\L}ojasiewicz (KL)  property, we demonstrate the convergence of the entire sequence generated by the LIPAL algorithm and derive improved local convergence rates that exhibit a dependence on the KL parameter.

%\red{$\mathrm{(iv)}$~
%We acknowledge that the complexity of perturbed augmented Lagrangian algorithms for general problems remains substantial, typically requiring at least $\mathcal{O}(\epsilon^{-3})$ iterations to achieve an $\epsilon$-first-order optimal solution.  Furthermore, these theoretical complexity bounds often do not accurately reflect the observed practical performance of these algorithms.  To address this discrepancy, we specialize our LIPAL algorithm for a class of optimization problems exhibiting benign nonconvex properties. % For this restricted class, we establish improved complexity bounds ranging from $\mathcal{O}(\epsilon^{-2})$ to $\mathcal{O}(\epsilon^{-1})$ to attain even an $\epsilon$-second-order optimal solution. % This significant improvement in complexity bounds bridges the gap between theoretical predictions and practical observations for such algorithms. To the best of our knowledge, this represents the first complexity result for an augmented Lagrangian-type algorithm specifically designed for solving this class of benign nonconvex problems with nonlinear equality constraints.?}

\item 
Finally, in addition to proposing a novel algorithm and providing its convergence guarantees, we demonstrate the algorithm's efficiency through numerical experiments. These experiments involve a comparative analysis with existing methods from the literature on large-scale clustering problems, employing the Burer-Monteiro factorization technique for semidefinite programming.
\end{compactenum}

%%% The structure of paper.
\medskip 
\noindent\textbf{\textit{Paper organization.}}
The rest of this paper is organized as follows. In Section \ref{sec2}, we formally introduce the optimization problem that is the focus of our study, along with the key assumptions that underpin our analysis. We also introduce a novel regularity condition. Section \ref{sec3} presents the proposed algorithm in detail. Subsequently, in Section \ref{sec4}, we provide a rigorous analysis of the algorithm's convergence properties, while in Section \ref{sec_reg} we provide a detailed discussion on the regularity condition introduced in this work. %In Section \ref{sec:improved_results}, we delve deeper into the convergence behavior by establishing improved rates under specific, yet reasonable, nonconvexity assumptions. Section \ref{sec_reg} offers a comprehensive discussion of the newly introduced regularity condition, exploring its implications and connections to existing conditions. 
Finally, in Section \ref{sec5}, we present a comparative numerical study, showcasing the performance of our proposed algorithm relative to existing methods using two different optimization models.

%%%%%%%%%%%%%%%%%%%%%%%%%%%%%%%%%%%
%%% Basic notations and definitions
%\subsection{Basic notations and definitions.}
\medskip 
\noindent\textbf{\textit{Basic notations, concepts, and terminologies.}}
We recall some necessary notions and results from variational analysis used in this paper.  
First, we start by defining the limiting normal cone of a set $C\subset\mathbb{R}^n$ at a point $\bar{x}$. Then, we define the basic and singular subdifferentials of an extended real-valued function through (limiting) normals to its epigraph. 

\begin{definition}\label{subdiffc}  
[Normal vectors (see  Definition 6.3 in \cite{RocWet:98})]
Let $C \subset \mathbb{R}^n$ be a given subset of $\mathbb{R}^n$ and $\bar{x} \in C$. 
Then, the regular normal cone  to $C$ at $\bar{x}$ is defined as
\begin{equation*}
   \widehat{N}_C(\bar{x}) := \Bigg\{v\in\mathbb{R}^n:\limsup_{\substack{x \xrightarrow[C]{} \bar{x} \\ x \neq \bar{x}}} \frac{\langle v, x - \bar{x}   \rangle}{\Vert x - \bar{x}\Vert  } \leq 0   \Bigg\}, 
\end{equation*}
where $x \xrightarrow[C]{} \bar{x}$ denotes  the limit  for $x\in C$.
The limiting normal cone to $C$ at $\bar{x}$ is defined as
\begin{equation*}
  N^{\text{lim}}_C(\bar{x}) := \left\{v\in\mathbb{R}^n: \exists x^{k} \xrightarrow[C]{} \bar{x} \text{ and } v_k \xrightarrow[]{} v \text{ with } v_k \in \widehat{N}_C(x^{k})
  \right\},
\end{equation*}
where $x^{k} \xrightarrow[C]{} \bar{x}$ denotes  the limit of a sequence $\{x^k\}_{k\geq0}$  to $\bar{x}$ while remaining within the set $C$.
        
\noindent If the set $C$ is convex, we have $\widehat{N}_C(\bar{x})= N^{\text{lim}}_C(\bar{x}) =N_C(\bar{x}) $, called the  normal cone
to $C$ at $\bar{x}$.
\end{definition}

\medskip 

\begin{definition}\label{subdiff}
[basic and singular subdifferentials (Definition 1.77 in \cite{Mor:06} and Theorem 8.9 in \cite{RocWet:98})] 
Consider a  function $\phi: \mathbb{R}^n \rightarrow \bar{\mathbb{R}} := [-\infty,\infty]$ and a point $\bar{x} \in \mathbb{R}^n$ with $\vert \phi(\bar{x}) \vert < \infty$.
\begin{compactenum}
\item The set
\begin{equation*}
\partial \phi(\bar{x}) := \{v \in \mathbb{R}^n \mid (v, -1) \in N^{\text{lim}}_{\text{epi} \, \phi}(\bar{x}, \phi(\bar{x})\}
\end{equation*}
is the (basic or limiting) subdifferential of $\phi$ at $\bar{x}$, and its elements are basic subgradients of $\phi$ at this point. We define $\partial \phi(\bar{x}) := \emptyset$ if $|\phi(\bar{x})| = \infty$.

\item The set
\begin{equation*}
\partial^\infty \phi(\bar{x}) := \{v \in \mathbb{R}^n \mid (v, 0) \in N^{\text{lim}}_{\text{epi} \, \phi}(\bar{x}, \phi(\bar{x})\}
\end{equation*}
is the singular subdifferential of $\phi$ at $\bar{x}$, and its elements are singular subgradients of $\phi$ at this point. We define $\partial^\infty \phi(\bar{x}) := \emptyset$ if $|\phi(\bar{x})| = \infty$.
\end{compactenum}
\end{definition}

\medskip 

\noindent For a given proper and lower semicontinuous (lsc) function $g : \mathbb{R}^n \rightarrow \bar{\mathbb{R}}$ and a point $\bar{x} \in \dom g$, we denote the $g$-attentive convergence of a sequence $\{x^k\}$ as $x^k \xrightarrow[]{g} \bar{x}$, which means $x^k \rightarrow \bar{x}$ and $g(x^k) \rightarrow g(\bar{x})$ as $k \to \infty$. 
Let $\Phi:\mathbb{R}^d \to \bar{\mathbb{R}}$ be a proper lsc function. 
For $−\infty < \tau_1 < \tau_2 \leq +\infty$, we define 
\begin{equation*}
[\tau_1 < \Phi < \tau_2] = \{x \in\mathbb{R}^d :\tau_1<\Phi(x)<\tau_2\}.
\end{equation*}
Let $\tau\in(0,+\infty]$.
We denote by $\Psi_{\tau}$ the set of all continuous concave functions $\varphi: [0, \tau] \to [0,+\infty)$ such
that $\varphi(0) = 0$ and $\varphi$ is continuously differentiable on $(0, \tau)$, with $\varphi'(s) > 0$ over $(0, \tau)$. 
Next, let us define class of functions satisfying the Kurdyka-{\L}ojasiewicz (KL) property, which is a generalization of  strongly convex functions \cite{AttBol:13}. 

\medskip 

\begin{definition} \label{def2}
Let $\Phi : \mathbb{R}^d \to \bar{\mathbb{R}}$  be a proper lsc function that takes constant value on $\Omega$. 
We say that $\Phi$ satisfies the KL property on $\Omega$ if there exists $ \epsilon>0, \varepsilon>0$, and $\varphi\in\Psi_{\varepsilon}$ such that 
$\forall x^{*} \in \Omega$ and for all $x$ in the intersection $\{x\in\mathbb{R}^d: \text{ dist}(x,\Omega)<\epsilon\}\cap[\Psi(x^{*})<\Psi(x)<\Psi(x^{*})+\varepsilon]$, we have
\begin{equation*}
 \varphi'\big(\Phi(x) − \Phi(x^{*})\big)\mathrm{dist}\big(0, \partial\Phi(x)\big) \geq1.
\end{equation*}
\end{definition}

\noindent This definition covers many classes of functions arising in practical optimization problems, including the following examples:
\begin{compactitem}
\item If $f$ is a proper closed semialgebraic function, then $f$ is a KL function,  see \cite{AttBol:13}. 
\item The function $g(Ax)$, where $g$ is strongly convex on a compact set and  twice differentiable, and $A \in \mathbb{R}^{m\times n}$, is
a KL function. 
\item Convex piecewise linear/quadratic functions such as $\Vert x\Vert _{1}, \Vert x\Vert _{0}, \gamma\sum_{i=1}^{k}|x_{[i]}|$,  where $|x_{[i]}|$ is the $i$-th largest  entry in $x, \;k\leq n$ and $\gamma \in (0, 1]$; the indicator function $\delta_{\Delta}(x)$, where $\Delta= \{x\in\mathbb{R}^n : e^T x = 1, x \geq 0\}$;  least-squares problems with the Smoothly Clipped Absolute Deviation (SCAD) \cite{Fan:97};  and Minimax Concave Penalty (MCP)  regularized functions \cite{Zha:10} are all KL functions.  
\end{compactitem}
For a twice differentiable function $\varphi : \R^n \to \R$, we denote its gradient at a given point $x$ by $\nabla \varphi(x)\in\mathbb{R}^n$ and its hessian at $x$ by $\nabla^2 \varphi(x)\in\mathbb{R}^{n\times n}$. 

%%%%%%%%%%%%%%%%%%%%%%%%%%%%%%%%%%%%%%%%
%%% 2. Problem formulation and  assumptions
%%%%%%%%%%%%%%%%%%%%%%%%%%%%%%%%%%%%%%%%
\section{Problem Formulation and Model Assumptions.}\label{sec2}
In this paper, we study the following general nonsmooth nonconvex composite optimization problem with nonlinear equality constraints:
\begin{equation}
\begin{aligned}\label{eq:nlp_gen}
& \underset{x\in\mathbb{R}^n}{\min}
& & \phi(x) \triangleq{ f(x)+g(x) }, \\
& \hspace{0.0cm}\textrm{subject to }
& & \hspace{0.05cm} F(x)=0,
\end{aligned}
\end{equation}
where $f : \R^n \to \R$  and $F : \R^n \to \R^m$ are  smooth nonlinear functions, while $g: \R^n \to \bar{\R}$ is a simple proper lsc function (e.g., the indicator function of a simple set). 
Note that problem \eqref{eq:nlp_gen} is very general, as the functions $f, F$ and $g$ are possibly nonconvex, and encompasses a wide spectrum of applications. 
For example, the constrained composite optimization problem:
\begin{align*}
\min_{z \in \mathcal{Z}} \big\{  \ell(z) + h(G(z)) \big\} \quad \text{subject to} \ \; G(z) \in \mathcal{Y},
\end{align*}
that usually appears in optimal control (see \cite{MesBau:21}),  where $\ell : \R^{n_1} \to \R$  and $G : \R^{n_1} \to \R^m$ are  smooth functions,  $h: \R^{n_2} \to \bar{\R}$ is a proper lsc function, while $\mathcal{Z}$ and $\mathcal{Y}$ are simple convex sets, 
can be easily recast in the form of  \eqref{eq:nlp_gen}, by defining \(x=(z,y)\), \(F(x) = G(z) - y\), \(f(x) = \ell(z)\), and \(g(x) = h(y) + \delta_{\mathcal{Z}}(z) + \delta_{\mathcal{Y}}(y)\), where \(\delta_{\Omega}\) denotes the indicator function of a given set \(\Omega\). Hence, our problem \eqref{eq:nlp_gen} can also deal  with inequality constraints.  

\medskip 

To investigate  the optimization model \eqref{eq:nlp_gen}, we impose the following assumptions:

%%% Assumption A.1.
\begin{assumption}\label{assump1}
There exists $\rho_0\geq0$ such that $\Phi_{\rho_0}(x) \triangleq{ f(x) + g(x) + \frac{\rho_0}{2}\Vert F(x) \Vert^2 }$ has  compact level sets, i.e. for any value  $\alpha\in\mathbb{R}$, the following sublevel set is empty or compact:
\begin{equation*}
\mathcal{S}_{\alpha}^0\triangleq{\{x \ : \, \Phi_{\rho_0}(x) \leq \alpha \} }.
\end{equation*}
\end{assumption}
Note that Assumption \ref{assump1} holds provided that the objective function is strongly convex or  $\Phi_{\rho_0}(\cdot)$ is coercive. 
An immediate consequence of Assumption \ref{assump1} is the fact that the function  $\Phi_{\rho_0}(\cdot)$ is  bounded from below, i.e.:
\begin{equation}\label{lem1}
\underline{\Phi}\triangleq{\inf_{x\in\mathbb{R}^n} \Big\{ \Phi_{\rho_0}(x) \triangleq f(x) + g(x) + \tfrac{\rho_0}{2}\Vert F(x)\Vert ^2 \Big\}}>-\infty.   
\end{equation}

%%% Assumption A.2.
\begin{assumption}\label{as:Assum}
For any given compact set $\mathcal{S}\subseteq\mathbb{R}^n$, there exist positive constants $M_f$, $M_F$, $L_f$, $L_F$, and $\sigma$ such that $f$ and $F$ satisfy the following conditions:
\begin{itemize}
  \item[$\mathrm{(i)}$] $\Vert \nabla f(x)\Vert \leq M_f$ and $\Vert \nabla f(x)-\nabla f(y)\Vert \leq L_f\Vert x-y\Vert  \;\text{ for all } x, y\in\mathcal{S}$.
  \item[$\mathrm{(ii)}$] $\Vert {J_F}(x)\Vert \leq M_F$ and  $\Vert {J_F}(x)-{J_F}(y)\Vert \leq L_F\Vert x-y\Vert \;\text{ for all } x, y\in\mathcal{S}$.
  \item[$\mathrm{(iii)}$] $\sigma\Vert  F(x) \Vert \leq \mathrm{dist}\left(-J_F(x)^TF(x),\partial^{\infty}g(x)\right) \;\text{ for all } x\in\mathcal{S}$,
\end{itemize}
where  ${J_F}\in\mathbb{R}^{m\times n}$ denotes  the Jacobian matrix of the functional constraints $F$. 
\end{assumption}

\noindent Assumption \ref{as:Assum} allows us to consider quite general classes of problems. In particular, Condition $\mathrm{(i)}$ holds if $f(\cdot)$ is differentiable and $\nabla f(\cdot)$ is \textit{locally}  Lipschitz continuous on a neighborhood of $\mathcal{S}$. 
Conditions $\mathrm{(ii)}$ hold when $F(\cdot)$ is differentiable on a neighborhood of $\mathcal{S}$ and $J_F(\cdot)$ is \textit{locally} Lipschitz continuous on $\mathcal{S}$. 
Finally, the constraint qualification condition $\mathrm{(iii)}$, newly introduced in this work, guarantees, in particular, the existence of a KKT point  \cite{HalTeb:23},  as it is a generalization of LICQ condition which is commonly used in nonconvex optimization with smooth objective, see e.g.,  \cite{XieWri:21}. 
Moreover, this regularity condition $\mathrm{(iii)}$ combined with a perturbed ascent step of the dual variables  allows us to prove in the next sections the boundedness of the dual iterates. 
More detailed discussion on the  constraint qualification condition $\mathrm{(iii)}$ is given in Section \ref{sec_reg}.  

%%% Assumption A.3.
\begin{assumption}\label{assump3}
There exists a finite constant $\bar{\alpha} \in \R$ such that $\phi(x) \leq \bar{\alpha}$ for all $x \in \Lambda := \{x \in \dom g :\, \Vert F(x) \Vert \leq 1 \}$.
\end{assumption}

\noindent Assumption \ref{assump3} holds provided that e.g.,   the set $\Lambda := \{x \in \mathbb{R}^n : \, \Vert F(x)\Vert \leq 1\}$ is compact. 
However, this assumption is unnecessary in our analysis below if we can select the initial point $x^0$ such that $F(x^0) = 0$.  
In many cases, we can easily find such an $x^0$, e.g., when $m$ is sufficient small.
\medskip 

Next, we provide a definition of an $\epsilon$-first-order optimal solution of  \eqref{eq:nlp_gen}.

%%% Definition 4.
\begin{definition}\label{firstorder}[$\epsilon$ first-order optimal solution]
Let $\epsilon>0$ be a fixed accuracy. 
A vector $x^{*}$ is said to be an $\epsilon$-first-order optimal solution to \eqref{eq:nlp_gen} if  $\exists y^{*} \in\mathbb{R}^m$ such that:
\begin{equation*}
\mathrm{dist}\Big(-\nabla f(x^{*})-\nabla {F(x^{*})}^Ty^{*},\partial g(x^{*}) \Big) \leq \epsilon\hspace{0.3cm} \text{and} \hspace{0.3cm} \Vert F(x^{*})\Vert \leq\epsilon.
\end{equation*} 
\end{definition}

\noindent Definition \ref{firstorder} introduces the notion of an $\epsilon$-first-order optimal solution to \eqref{eq:nlp_gen}, which can be seen as an approximation of the KKT conditions within an accuracy  $\epsilon$. 
Specifically, under some appropriate regularity condition, related to Assumption \ref{as:Assum}-$\mathrm{(iii)}$,  we can  ensure the existence of  KKT points for problem  \eqref{eq:nlp_gen}, i.e. a pair $(\bar{x}, \bar{y})$  satisfying the following relations \cite{HalTeb:23}:
\begin{equation*}
-\nabla f(\bar{x})-\nabla {F(\bar{x})}^T\bar{y}\in\partial g(\bar{x})\quad \text{ and }\quad F(\bar{x})=0.
\end{equation*} 
However, in many practical scenarios we may not be able to find an exact KKT point, but instead we can identify an approximate solution ($\epsilon$-KKT point) as in  Definition \ref{firstorder}. 

%%%%%%%%%%%%%%%%%%%%%%%%%%%%%%%%%%%%%%%%%%%%%%%%%%%%%%%%%%
%%% 3. A Linearized Perturbed Augmented Lagrangian Method
%%%%%%%%%%%%%%%%%%%%%%%%%%%%%%%%%%%%%%%%%%%%%%%%%%%%%%%%%%
\section{A Linearized Perturbed Augmented Lagrangian Method.}\label{sec3}
In this section, we develop a \textit{novel  linearized perturbed augmented Lagrangian method} to solve \eqref{eq:nlp_gen} under the assumptions stated in Section~\ref{sec2}. 
The main idea is as follows.
\begin{itemize}
\item First, a perturbation is introduced at the level of  dual variables in the augmented Lagrangian function through a subunitary parameter. 
\item Then, at each iteration we linearize the smooth part of the objective and the functional constraints within the perturbed augmented Lagrangian  and add a quadratic regularization. 
\item Finally, the dual variables are then updated through a perturbed ascent step. 
\end{itemize}

Let $(x^0,y^0)\in\mathbb{R}^n\times\mathbb{R}^m$ be the starting points of our algorithm and  $\tau \in (0,1]$ be a perturbation parameter. 
The perturbed augmented Lagrangian function $\mathcal{L}^{\tau}_{\rho}$ associated with   \eqref{eq:nlp_gen} is defined as
\begin{equation}\label{eq:pert_AL}
\arraycolsep=0.2em
\begin{array}{lcl}
\psi^{\tau}_{\rho}(x,y;y^0)  & \triangleq &  f(x)+\langle \tau y^0 + (1-\tau)y, F(x) \rangle+\frac{\rho}{2}{\Vert F(x)\Vert ^2}, \vspace{1ex}\\
\mathcal{L}^{\tau}_{\rho}(x,y;y^0) & \triangleq &    g(x) + \psi^{\tau}_{\rho}(x,y;y^0).
\end{array}
\end{equation}
Here, $\psi^{\tau}_{\rho}(\cdot, y; y^0)$ is  differentiable w.r.t. $x$. 

\medskip 

Note that for $\tau=1$, the perturbed augmented Lagrangian function reduces to the standard quadratic penalty function. 
The previous Lagrangian type function for $y^0 = 0$ has been also considered in \cite{HajHon:19} in the context of nonconvex problems with linear constraints. 
The \textbf{anchor vector} $y^0$ is useful in our framework when restarting is incorporated in the implementation.  
The gradient of $\nabla_x\psi^{\tau}_{\rho}$ can be expressed as
\begin{equation}\label{eq:grad_psi}
     \nabla_x\psi^{\tau}_{\rho}(x,y;y^0)=\nabla f(x)+{J_F(x)}^T\left(\tau y^0 + (1-\tau)y + \rho F(x)\right).
\end{equation}
Note that if Assumption \ref{as:Assum} holds on a compact set $\mathcal{S}\subseteq\mathbb{R}^n$, then for any compact set $\mathcal{Y}\subseteq\mathbb{R}^m$ (containing the dual variables), the gradient   $\nabla_x\psi^{\tau}_{\rho}(x, y; y^0)$ is locally smooth, i.e. it is  Lipschitz continuous w.r.t. $x$ on the compact set $\mathcal{S}\times\mathcal{Y}$. 
In the following lemma we prove this statement.  For clarity, we provide the proofs of all the lemmas in Appendix.

%%% Lemma 1.
\begin{lemma}[Smoothness of $\psi^{\tau}_{\rho}$]\label{lemma2} 
If Assumption~\ref{as:Assum} holds on a compact set $\mathcal{S}$, then for any compact set $\mathcal{Y}\subset\mathbb{R}^m$ there exists $L^{\tau}_{\rho}>0$ such that:
\begin{equation*}
\Vert \nabla_x \psi^{\tau}_{\rho}(x,y;y^0)-\nabla_x \psi^{\tau}_{\rho}(x',y;y^0)\Vert \leq {L^{\tau}_{\rho}}\left\Vert x-x'\right\Vert  \hspace{0.5cm}\forall x,x'\in\mathcal{S}, \; \forall y\in\mathcal{Y},
\end{equation*}
where $L^{\tau}_{\rho}  \triangleq \sup_{(x,y)\in\mathcal{S}\times\mathcal{Y}}\left\{L_f+L_F\Vert \tau y^0 + (1-\tau)y +\rho F(x)\Vert +\rho M_F^2\right\}$.
\end{lemma}

%%% Proof of Lemma 1
\proof{Proof.} See Appendix.
\Eproof
\endproof
%%% End of Proof.

\medskip
The following linear functions will be  frequently used in the sequel:
\begin{equation}\label{eq:linearization_funcs}
\ell_f(x;\bar{x}):=f(\bar{x})+\langle\nabla f(\bar{x}),x-\bar{x}\rangle \quad\text{and} \quad    \ell_F(x;\bar{x}):=F(\bar{x})+ J_F(\bar{x})(x-\bar{x}) \hspace{0.5cm}, \forall x,\bar{x} \in \mathcal{S}.
\end{equation}

\noindent\textbf{The proposed algorithm.}
Let us now present the \textbf{new}  \textit{Linearized Perturbed Augmented Lagrangian} (LIPAL) algorithm for solving \eqref{eq:nlp_gen} (see Algorithm \ref{alg:LIPAL}).

%%%% Algorithm 1.
\begin{algorithm}\caption{(Linearized Perturbed Augmented Lagrangian (LIPAL))}\label{alg:LIPAL}
\begin{algorithmic}[1]
\State\textbf{Initialization:} 
Choose $x^0 \in \dom g$ and $y^0\in\mathbb{R}^m$, and parameters $\tau \in(0,1]$, $\rho > 0$, and $\beta > 0$.
\For{$k = 0, 1, \cdots, K$}
\State If an $\epsilon$-KKT condition is met, then \textbf{terminate}.
\State Set $y_{\tau}^k := \tau y^0 + (1 - \tau)y^k$ and update the primal step:
\begin{equation}\label{eq:LIPAL_subprob}
 x^{k+1}\gets\argmin\limits_{x \in \mathbb{R}^n} \Big\{\mathcal{Q}_k(x) \triangleq \ell_f(x; x^k) + g(x) + \iprods{ y^k_{\tau}, \ell_F(x;x^k)} + \tfrac{\rho}{2}\norms{\ell_F(x;x^k)}^2 + \tfrac{\beta}{2}\norms{x - x^k}^2 \Big\}.
\end{equation}
\State Update the dual step: $y^{k+1}\gets  y^k_{\tau} + \rho F(x^{k+1})$.
\EndFor
\end{algorithmic}
\end{algorithm}

\medskip

\noindent Note that if $g$ is a (weakly) convex function (or prox-regular), then the objective function in the subproblem \eqref{eq:LIPAL_subprob} of Algorithm \ref{alg:LIPAL} becomes strongly convex, provided that $\beta$ is chosen appropriately. 
Consequently, the solution of \eqref{eq:LIPAL_subprob} is well-defined and unique. 
The existence and boundedness of the primal iterates $x^{k+1}$ can be also ensured if the domain of $g$ is compact or $g$  is coercive.  
Furthermore, our algorithm is easy to implement as it uses only first-order information and the subproblem \eqref{eq:LIPAL_subprob} is simple provided that $g$ is a simple function (e.g., $g$ weakly convex). 
Thus it allows us to tackle large-scale instances of \eqref{eq:nlp_gen}. 

\medskip 

When it is difficult to compute an exact solution $x^{k+1}$ of \eqref{eq:LIPAL_subprob}, one can consider computing an inexact solution of \eqref{eq:LIPAL_subprob}. 
Under appropriate conditions on the inexactness,  our convergence results derived in the sequel will remain valid.
For instance, if $x^{k+1}$ in \eqref{eq:LIPAL_subprob} satisfies  the following condition:
For $\mathcal{Q}_k(\cdot)$ defined by \eqref{eq:LIPAL_subprob}, there exists $\alpha > 0$ and $s^{k+1}\in\partial\mathcal{Q}_k(x) \big|_{x = x^{k+1}}$ such that
\begin{equation*}
\Vert s^{k+1}\Vert \leq\alpha\Vert x^{k+1}-x^k\Vert \quad \textrm{and} \quad \mathcal{Q}_k(x^{k+1}) \leq \mathcal{L}^{\tau}_{\rho}(x^{k},y^{k};y^0). 
\end{equation*}
Note that any descent algorithm, when initialized at the current iterate $x^k$, can  guarantee the second descent condition. For example,   one can solve the subproblem in step 4 using an accelerated proximal gradient algorithm (possibly with line search for the Lipschitz constant of the gradient) \cite{GhaLan:16, Nes:18}. Similar approximate optimality  conditions have been considered e.g.,  in \cite{XieWri:21}. 

\medskip 

The subproblem \eqref{eq:LIPAL_subprob} can be viewed as a prox-linear operator in Gauss-Newton-type methods (also called prox-linear methods).
However,  we can linearize the entire smooth component of the perturbed augmented Lagrangian using a gradient-based approach as in \cite{HalTeb:23}.
In this case, \eqref{eq:LIPAL_subprob} can be replaced by the following subproblem:
\begin{equation*}
x^{k+1}\gets\argmin_{x \in \mathbb{R}^n} \Big\{ g(x)+ \iprods{\nabla_x\psi^{\tau}_{\rho}(x^k,y^k;y^0), x-x^k} +  \tfrac{\beta}{2}\norms{x - x^k}^2 \Big\}.
\end{equation*}
This modification  leads to  a proximal gradient update at Step 4, while still maintaining the validity of our convergence proofs in the next sections. 
However, in practice, we observed that the prox-linear (or Gauss-Newton) mechanism generally works better than the gradient-based strategy. 
Finally, Step 5 in Algorithm~\ref{alg:LIPAL} is a perturbed ascent step with a given \textbf{anchor point} $y^0$. 

\medskip 

In what follows, we will frequently use the following notations:
\begin{equation}\label{eq:PD_decrements}
\Delta{x}^k := x^{k+1} - x^k \quad \text{ and } \quad \Delta{y}^k := y^{k+1} - y^k.
\end{equation}
These quantities can be seen as the decrement of the primal and dual variables $x$ and $y$, respectively.

%%%%%%%%%%%%%%%%%%%%%%%%%%%%%%%%%%%%%%%
%%% 3. Convergence Analysis
%%%%%%%%%%%%%%%%%%%%%%%%%%%%%%%%%%%%%%%
\section{Convergence Analysis.}\label{sec4}
In this section, we first analyze the global convergence of our LIPAL, then we derive improved local convergence rates under KL property.  
The analysis is relatively technical and contains several lemmas.

%%% 3.1. Technical lemmas
\subsection{Technical lemmas}\label{subsec:technical_lemmas}
Let us first prove an initial bound of $\Vert y^{k}\Vert$ for the dual iterates.

%%% Lemma 2.
 \begin{lemma}\label{le:lambda_bou}
Let $\{ (x^k,y^k)\}$ be generated by Algorithm \ref{alg:LIPAL}. 
Suppose that Assumption \ref{as:Assum} holds on some compact set $\mathcal{S}$  on which the primal sequence $\{x^k\}_{k\geq0}$  belongs to.
Then, we have 
\begin{equation}\label{eq:bound_for_dual}
   \Vert y^{k}\Vert \leq\Vert y^0\Vert +\frac{\rho}{\tau}\Delta \quad \text{and} \quad    \Vert y^{k+1}-y^0\Vert \leq \frac{\rho}{\tau}\Delta,
\end{equation}
where $\Delta:=\sup_{x\in\mathcal{S}} \Vert F({x})\Vert$.
\end{lemma}

%%% Proof of Lemma 2
\proof{Proof.} See Appendix.
\Eproof
\endproof
%%% End of Proof.

\medskip 

%%% Remark 1.
\begin{remark}\label{re:Lipschitz_constant}
As a consequence of Lemma~\ref{le:lambda_bou}, when the primal iterates $\{x^k\}$ belong to a compact set $\mathcal{S}$, it follows that the ball centered at the origin of radius $R \triangleq \Vert y^0\Vert + \frac{\rho}{\tau}\Delta$ (denoted $\mathbb{B}_R$) contains all the dual iterates $\{y^k \}$. 
Hence, the Lipschitz constant of $\nabla_x\psi^{\tau}_{\rho}$ is independent of the dual iterates of Algorithm \ref{alg:LIPAL}, i.e.:
\begin{equation*}
\arraycolsep=0.2em
\begin{array}{lcl}
L^{\tau}_{\rho}  & := & L_f + L_F\sup_{(x, y)\in\mathcal{S}\times\mathbb{B}_R}\left\{ \norms{ \tau y^0 + (1-\tau)y   + \rho F(x)} \right\} +   \rho M_F^2 \vspace{1ex} \\
& \leq & L_f + (1-\tau)L_F\Vert y^0\Vert +\frac{\rho}{\tau}L_F\Delta +\tau L_F \Vert y^0\Vert  +   \rho M_F^2 \vspace{1ex} \\
& = & L_f + L_F\Vert y^0\Vert +\frac{\rho}{\tau}L_F\Delta  +   \rho M_F^2.
\end{array}
\end{equation*}
\end{remark}

Our next step is to prove a bound for the dual iterates $\{y^k\}$ generated by Algorithm~\ref{alg:LIPAL}.
This result is new and important (see our discussion in the ``Introduction" and also in \cite{HalTeb:23}).

%%% Lemma 3.
 \begin{lemma}[Bound for $\Vert y^{k}\Vert $]\label{le:lambda_bou1}
Let $\{ (x^k,y^k) \}$ be generated by Algorithm \ref{alg:LIPAL}. 
Suppose that Assumption \ref{as:Assum} holds on some compact set $\mathcal{S}$  on which the primal sequence $\{x^k\}_{k\geq0}$  belongs to, and $\beta\geq 2L^{\tau}_{\rho}$.
Then,  we have
\begin{equation}\label{eq:lambda_1}
 \Vert y^{k+1}-y^0\Vert   \leq  \frac{1}{\sigma} \Big[  M_f+1+\Vert y^0\Vert  +  2\beta \Vert \Delta x^{k}\Vert +\frac{(1-\tau)}{\tau} (M_F + \sigma) \Vert \Delta y^k\Vert \Big].
\end{equation}
\end{lemma}

%%% Proof of Lemma 3
\proof{Proof.} See Appendix.
\Eproof
\endproof
%%% End of Proof.

\medskip 

Now, we establish a bound on  $\Vert \Delta y^{k}\Vert ^2$ using a similar technique as in \cite{Lu:22}. 
This bound will be useful to construct a Lyapunov function later.

%%% Lemma 4.
\begin{lemma}[Bound for $\Vert \Delta y^{k}\Vert $]\label{le:lambda_bound}
Let  $\{ (x^k,y^k)\}$ be generated by Algorithm \ref{alg:LIPAL}. I
Suppose that Assumption \ref{as:Assum} holds on some compact set $\mathcal{S}$  on which $\{x^k\}_{k\geq0}$  belongs to.
Then, we have
\begin{equation}\label{eq:lambda_squared}
(1-\tau)\norms{\Delta y^k}^2\leq (1-\tau)\norms{\Delta y^{k-1}}^2 +\frac{\rho^2M_F^2}{2\tau}\norms{\Delta x^k}^2 - \frac{\tau}{2}\Vert \Delta y^k\Vert ^2.
\end{equation}
\end{lemma}

%%% Proof of Lemma 4
\proof{Proof.} See Appendix.
\Eproof
\endproof
%%% End of Proof.

\medskip 

Our key step is to prove the following descent property of $\mathcal{L}_{\rho}^{\tau}$ defined by \eqref{eq:pert_AL}.

%%% Lemma 5.
\begin{lemma}[Descent property of $\mathcal{L}^{\tau}_{\rho}$]\label{lemma3} 
Let $\{ (x^k,y^k) \}$ be generated by Algorithm \ref{alg:LIPAL}. 
Suppose that Assumption \ref{as:Assum} holds on some compact set $\mathcal{S}$  on which the primal sequence $\{x^k\}_{k\geq0}$  belongs to, and $\beta\geq 2L^{\tau}_{\rho}$.
Then, for all $k\geq 0$, we have 
\begin{equation}\label{eq:descent_of_PAL}
\mathcal{L}^{\tau}_{\rho}(x^{k},y^{k};y^0) - \mathcal{L}^{\tau}_{\rho}(x^{k+1},y^{k};y^0) \geq  \frac{\beta}  {4}\Vert x^{k+1}-x^{k}\Vert ^2 \geq \frac{L^{\tau}_{\rho}}{2}\Vert x^{k+1}-x^{k}\Vert ^2. 
\end{equation}
\end{lemma}

%%% Proof of Lemma 5
\proof{Proof.} See Appendix.
\Eproof
\endproof
%%% End of Proof.

\medskip 
\noindent Note that in practice the regularization parameter \( \beta \) can be generated dynamically (through a line search procedure) to ensure the decrease in \eqref{eq:descent_of_PAL} at each iteration.   To obtain an approximate solution, we need to upper bound both the optimality and feasibility measures through the primal and dual decrements $\Delta{x}^k$ and $\Delta{y}^k$ defined by \eqref{eq:PD_decrements}, as in the following lemma.

%%%% Lemma 6.
\begin{lemma}[Bound of optimality measures]\label{le:measure_critic}
Let $\{ (x^k,y^k) \}$ be generated by Algorithm \ref{alg:LIPAL}. 
Suppose that Assumption \ref{as:Assum} holds on some compact set $\mathcal{S}$  on which the primal sequence $\{x^k\}_{k\geq0}$  belongs to, and $\beta \geq 2L^{\tau}_{\rho}$.
Then, we have 
\begin{equation}\label{eq:measure_critic}
\arraycolsep=0.2em
\left\{\begin{array}{ll}
	& \mathrm{dist}\big( -\nabla f(x^{k+1}) - {{J_F}(x^{k+1})}^Ty^{k+1},\partial g(x^{k+1}) \big)  \leq    2\beta\Vert \Delta x^{k}\Vert, \vspace{1ex}\\ 
	& \Vert F(x^{k+1})\Vert  \leq   \frac{2\beta\tau}{\rho\sigma}\Vert \Delta x^{k}\Vert +\frac{1-\tau}{\rho}\left(2+\frac{M_F}{\sigma}\right)\Vert \Delta y^{k}\Vert +\frac{\tau}{\rho\sigma} \left( M_f+1+\Vert y^0\Vert \right).
\end{array}\right.
\end{equation}
\end{lemma}

%%% Proof of Lemma 6
\proof{Proof.} See Appendix.
\Eproof
\endproof
%%% End of Proof.

%%%% The Lyapunov function.
\medskip

\noindent\textbf{The Lyapunov function.}
We introduce the following Lyapunov function (similar to  \cite{Lu:22}):
\begin{equation}\label{eq:Lyapunov_func}
P(x,y,y';y^0)=\mathcal{L}^{\tau}_{\rho}(x,y;y^0)-\frac{\tau(1-\tau)}{2\rho}\Vert y-y^0\Vert ^2+\frac{2(1-\tau)^2}{\tau\rho}\Vert y-y'\Vert ^2.
\end{equation}
The value of this function at each triple $(x^k,y^k,y^{k-1})$ for $k\geq 1$ is denoted by
\begin{equation}\label{eq:lyapunov_function}
 \Pc_{k}=P(x^k,y^k,y^{k-1};y^0).
\end{equation}
Our next step is to prove that the sequence $\{\Pc_{k}\}_{k\geq1}$ is decreasing  and bounded from bellow. Let us first prove that  $\{ \Pc_k \}_{k\geq 1}$ is decreasing.
 
%%% Lemma 7.
 \begin{lemma}[Decrease of Lyapunov function]\label{le:decrease_of_P}  
Let $\{ (x^k,y^k) \}$ be generated by Algorithm \ref{alg:LIPAL}. 
Suppose that Assumption \ref{as:Assum} holds on some compact set $\mathcal{S}$  on which the primal sequence $\{x^k\}_{k\geq0}$  belongs to. 
Suppose further that $\beta$ is chosen such that
\begin{equation}\label{eq:gamma_rho}
\beta \geq \max\left\{ 2L_\rho^\tau, \ \tfrac{8(1-\tau)\rho M_F^2}{\tau^2} \right\}.
\end{equation}
Then, for all $k \geq 1$, $\Pc_k$ defined by \eqref{eq:lyapunov_function} satisfies 
\begin{equation}\label{eq:decrease_Lyapunov}
	\Pc_{k+1}-\Pc_{k} \leq  -\frac{\beta}{8}\Vert \Delta x^{k}\Vert ^2-\frac{\tau(1-\tau)}{2\rho}\Vert \Delta y^{k}\Vert^2. 
\end{equation}
\end{lemma}

%%% Proof of Lemma 7
\proof{Proof.} See Appendix.
\Eproof
\endproof
%%% End of Proof.

\medskip 

In the sequel, we assume that $x^0\in \dom g$ is chosen such that
\begin{equation}\label{eq:choice_of_x0}
	\Vert F(x^0)\Vert ^2\leq\min\left\{1,\frac{c_0}{\rho}\right\} \hspace{0.4cm} \textrm{for some}~ c_0 \geq 0.
\end{equation}
Then, from Assumption \ref{assump3}, it follows that  $\phi(x^0) \triangleq f(x^0) + g(x^0) \leq \bar{\alpha}$. 
Let us also define
\begin{equation}\label{eq:P_upper_bound}
	\bar{\Pc} := \left(\frac{18\left(1-\tau\right)^2}{\tau}+1\right) \left(\bar{\alpha}+c_0+2\Vert y^0\Vert ^2\right)-\frac{18\left(1-\tau\right)^2}{\tau}\underline{\Phi}. 
\end{equation}
Moreover, for $\rho_0 > 0$ given in Assumption \ref{assump1}, we also choose $\rho$ such that
\begin{equation}\label{eq:choice_of_rho_b1}
\rho \geq \max\left\{\!{1}, 3\rho_0 \right\}.
\end{equation}
By the definition \eqref{eq:pert_AL} of $\mathcal{L}^{\tau}_{\rho}$, we have 
\begin{equation}\label{eq:pert_AL_at_x0}
\arraycolsep=0.2em
\begin{array}{lcl} 
	\mathcal{L}^{\tau}_{\rho}(x^0,y^0;y^0) &= & \phi(x^0) + \langle y^0  ,F(x^0)\rangle+\frac{\rho}{2}\Vert F(x^0)\Vert ^2  \vspace{1ex} \\
	& \leq & \phi(x^0) + \frac{\Vert y^0\Vert ^2}{2\rho}+\frac{\rho}{2}\Vert F(x^0)\Vert ^2+\frac{\rho}{2}\Vert F(x^0)\Vert ^2 \vspace{1ex}\\
	& \overset{\tiny\eqref{eq:choice_of_x0}}{\leq} & \bar{\alpha}+\frac{1}{2\rho}\Vert y^0\Vert ^2+c_0.
\end{array}    
\end{equation}
This inequality implies that
\begin{equation}\label{eq:line11}
\arraycolsep=0.1em
\begin{array}{lcl} 
	\bar{\alpha}+c_0 - \underline{\Phi} & \geq & \phi(x^0) + \langle y^0  ,F(x^0)\rangle+\frac{\rho}{2}\Vert F(x^0)\Vert ^2-\frac{1}{2\rho}\Vert y^0\Vert ^2-\underline{\Phi} \vspace{1ex}\\
	& \overset{\tiny{\rho\geq3\rho_0}}{\geq} & \phi(x^0) + \frac{\rho_0}{2}\Vert F(x_0)\Vert ^2-\underline{\Phi}+\langle y^0 ,F(x^0)\rangle+\frac{\rho}{3}\Vert F(x^0)\Vert ^2-\frac{\Vert y^0\Vert ^2}{2\rho} \vspace{1ex} \\
	& {\overset{{\eqref{lem1}}}{\geq}} & \frac{\rho}{3}\Vert F(x^0)+\frac{3y^0}{2\rho}\Vert ^2-\frac{3}{4\rho}\Vert y^0\Vert ^2 \!-\frac{\Vert y^0\Vert ^2}{2\rho}  {\overset{{(\rho\geq1)}}{\geq}}-\frac{5}{4}\Vert y^0\Vert ^2.
\end{array}    
\end{equation}
%%%
The following lemma shows that if $\{x_{k}\}_{k\geq0}$ generated by  Algorithm \ref{alg:LIPAL} is bounded, then $\{P_{k}\}_{k\geq0}$ is also bounded.

%%% Lemma 8.
\begin{lemma}\label{bbound} 
Let $(x^k,y^k)$ be generated by Algorithm \ref{alg:LIPAL} and  $\{\Pc_{k}\}_{k\geq0}$ be defined by \eqref{eq:lyapunov_function}. 
Suppose that Assumptions \ref{assump1}, \ref{as:Assum}, and \ref{assump3} hold on some compact set $\mathcal{S}$  on which the primal sequence $\{x^k\}_{k\geq0}$  belongs to.
Moreover, suppose that $\rho$ is chosen as in \eqref{eq:choice_of_rho_b1} and $x^0$ is chosen as in \eqref{eq:choice_of_x0}.
Then, there exist $\underline{\Pc} > -\infty$ and $\bar{\Pc}$ defined in \eqref{eq:P_upper_bound} such that for all $k \geq 0$, we have 
\begin{equation}\label{eq:important}
	\underline{\Pc} \leq \Pc_k \leq \bar{\Pc}.
\end{equation}
\end{lemma}

%%% Proof of Lemma 8
\proof{Proof.} See Appendix.
\Eproof
\endproof
%%% End of Proof.

\medskip 

Note that Assumptions \ref{assump1} and \ref{assump3} are used solely for establishing the upper bound $\Pc_u$ on $\{\Pc_k\}$.
They play no role in the derivation of the lower bound $\underline{\Pc}$ or in the proofs of other lemmas. 

%%%%%%%%%%%%%%%%%%%%%%%%%%%%%%%%%%%%%%%%%
%%% 4.2. Global Convergence
%%%%%%%%%%%%%%%%%%%%%%%%%%%%%%%%%%%%%%%%%
\subsection{Global convergence.}\label{subsec:globa_conver}
Now, we are ready to present our first main convergence  results of this paper.
We  derive the  complexity of  Algorithm \ref{alg:LIPAL} to find an $\epsilon$-first-order optimal solution to \eqref{eq:nlp_gen} (in terms of the number of evaluations of the objective function, gradient, functional constraints, and Jacobian matrix of functional constraints).

%%% Theorem 1.
\begin{theorem}[Iteration complexity]\label{th:complex_bound1} 
Let $\{(x^k,y^k)\}_{k\geq0}$ be generated by Algorithm \ref{alg:LIPAL} to solve \eqref{eq:nlp_gen}. 
Suppose that Assumptions \ref{assump1}, \ref{as:Assum}, and \ref{assump3} hold on some compact set $\mathcal{S}$  on which the primal sequence $\{x^k\}_{k\geq0}$  belongs to.
Suppose further that $\rho$ is chosen as in \eqref{eq:choice_of_rho_b1}, $x^0$ is chosen as in \eqref{eq:choice_of_x0}, and $\beta$ is selected such that
\begin{equation}\label{eq:choice_of_beta}
	\beta\geq\max\left\{2L_\rho^\tau, \frac{8(1-\tau)\rho M_F^2}{\tau^2}\right\}.
\end{equation}
Then, we have
\begin{equation}\label{eq:FO_convergence_limits}
\lim_{k\to\infty}\Vert \Delta x^k\Vert =0 \quad \textrm{and} \quad \lim_{k\to\infty}\Vert \Delta y^k\Vert =0. 
\end{equation}
Additionally, for any $\epsilon>0$, if  $\rho$ is chosen such that
\begin{equation}\label{eq:choice_of_rho}
 \rho \geq \max\left\{1, 3\rho_0, \frac{2\tau (M_f+1+\Vert y^0\Vert )}{\sigma\epsilon}\right\}, 
\end{equation}
 then Algorithm \ref{alg:LIPAL} yields an $\epsilon$-first-order optimal solution to \eqref{eq:nlp_gen} after $K := \mathcal{O}\left(\frac{1}{\epsilon^3}\right)$ iterations. 
\end{theorem}

%%% Proof of Theorem 1.
\proof{Proof.}
From  \eqref{eq:decrease_Lyapunov}, for all $k\geq 1$, we have
\begin{equation*}
	\frac{\beta}{8}\Vert \Delta x^{k}\Vert ^2+\frac{\tau(1-\tau)}{2\rho}\Vert \Delta y^{k}\Vert ^2\leq \Pc_{k}-\Pc_{k+1}.
\end{equation*}
Summing up this inequality from $k := 1$ to $k := K$, and noting that $\Pc_1 \leq \bar{\Pc}$ and $\Pc_{K+1} \geq \underline{\Pc}$, we get
\begin{equation}\label{eq:th1_limit}
\sum_{k=1}^{K}{\left( \frac{\beta}{8}\Vert \Delta x^k \Vert ^2+\frac{\tau(1-\tau)}{2\rho}\Vert \Delta y^k \Vert ^2\right)} \leq \Pc_{1} - \Pc_{K+1}\overset{\tiny\eqref{eq:important}}{\leq} \bar{\Pc} - \underline{\Pc} < \infty. 
\end{equation}
Since \eqref{eq:th1_limit} holds for any $K \geq 1$, passing to the limit as $K\to\infty$, we obtain
\begin{equation*}
\sum_{k = 1}^{\infty}{\left(\frac{\beta}{8}\Vert \Delta x^k \Vert ^2+\frac{\tau(1-\tau)}{2\rho}\Vert \Delta y^k \Vert ^2\right)}  < +\infty.
\end{equation*}
This summable result  together with the facts that $\beta>0$ and $\tau\in(0,1)$ yield 
\begin{equation*} 
	\lim_{k\to\infty}{\Vert \Delta x^{k}\Vert }=0 \quad \textrm{and}\quad  \lim_{k\to\infty}{\Vert \Delta y^{k}\Vert }=0,
\end{equation*}
which exactly proves \eqref{eq:FO_convergence_limits}.  Since $\{x^k\} \subset \mathcal{S}$, where $\mathcal{S}$ is compact, and $\{y^k\}$ is bounded as proved in Lemma~\ref{le:lambda_bou1},  we conclude that $\{(x^{k},y^{k})\}_{k\geq0}$ is bounded.
Therefore, there exists a convergent subsequence $\{(x^{k},y^{k})\}_{k \in \mathcal{K}}$ converging to a limit point $(x^{*}, y^{*})$.
Since $F$, $\nabla{f}$, and $J_F$ are continuous, and $\partial{g}$ is closed, passing the limit as $k \in \mathcal{K}$ tends to $+\infty$ in the two measurements \eqref{eq:measure_critic} of Lemma \ref{le:measure_critic}, we can show that
\begin{equation*}
\mathrm{dist}\big( -\nabla f(x^{*})-{J_F(x^{*})}^Ty^{*},\partial g(x^{*}) \big) = 0 \quad \textrm{ and } \quad \Vert F(x^{*})\Vert \leq\frac{\tau}{\rho}\Vert y^{*}-y^0\Vert \leq \frac{\tau (M_f+1+\Vert y^0\Vert )}{\rho\sigma}.
\end{equation*}
Let $\epsilon>0$ and $\rho\geq\frac{2\tau (M_f+1+\Vert y^0\Vert )}{\sigma\epsilon} = \mathcal{O}\left(\frac{\tau}{\epsilon}\right)$. 
Then, for some fixed $c_0$, if $x_0$ is chosen such that $x^0 \in \dom g$ and $\Vert F(x^0)\Vert ^2 \leq \frac{c_0}{\rho}$,  then from \eqref{eq:th1_limit}, we have (recall that $\bar{\Pc}$ is independent of $\rho$):
\begin{equation*}
\sum_{k = 1}^{K}{\left( \frac{\beta}{8}\Vert \Delta x^k \Vert ^2+\frac{\tau(1-\tau)}{2\rho}\Vert \Delta y^k \Vert ^2\right)} \leq \bar{\Pc} - \underline{\Pc}. 
\end{equation*}
Hence, there exists $k^{*}\in\{1, \ldots ,K\}$ such that
\begin{equation*}
\frac{\beta}{8}\Vert \Delta x^{k^{*}}\Vert ^2+\frac{\tau(1-\tau)}{2\rho}\Vert \Delta y^{k^{*}}\Vert ^2\leq\frac{\bar{\Pc} - \underline{\Pc}}{K}
\end{equation*}
This expression implies that 
\begin{equation*}
\Vert \Delta x^{k^{*}}\Vert \leq\sqrt{\frac{8(\bar{\Pc} - \underline{\Pc})}{\beta K}} \quad \textrm{and} \quad
\Vert \Delta y^{k^{*}}\Vert \leq\sqrt{\frac{2\rho(\bar{\Pc} - \underline{\Pc} )}{\tau(1-\tau)K}}. 
\end{equation*}             
Now, substituting these bounds into  \eqref{eq:measure_critic} of Lemma \ref{le:measure_critic} and noting that $\bar{\Pc} - \underline{\Pc} =\mathcal{O}\left( \frac{1}{\tau} \right)$ (see \eqref{eq:P_upper_bound} and \eqref{eq:decrease_Lyapunov}), we can derive that
\begin{equation*}
\arraycolsep=0.2em
\begin{array}{lcl}
	\mathrm{dist} \big( -\nabla f(x^{k^{*}+1})-{{J_F}(x^{k^{*}+1})}^Ty^{k^{*}+1},\partial g(x^{k^{*}+1}) \big) & \leq & 2\beta\Vert \Delta x^{k^{*}}\Vert \vspace{1ex}\\
	& \leq &  2\beta\sqrt{\frac{8( \bar{\Pc} - \underline{\Pc} )}{\beta K}} \vspace{1ex}\\
	& \leq & \mathcal{O}\left(\sqrt{\frac{\beta}{\tau}}\right)\frac{1}{\sqrt{K}}.
\end{array}
\end{equation*}
This estimate shows that for any $ \epsilon>0$ if $K \geq \mathcal{O}\left(\frac{\beta}{\tau\epsilon^2}\right)$, then we obtain $\mathrm{dist} \big(-\nabla f(x^{k^{*}+1})-{{J_F}(x_{k^{*}+1})}^Ty^{k^{*}+1},\partial g(x^{k^{*}+1}) \big)\leq\epsilon$.  Now, from the second line of \eqref{eq:measure_critic} in Lemma \ref{le:measure_critic} and the choice of $\rho$ in \eqref{eq:choice_of_rho}, we can also show that
\begin{equation*}
\arraycolsep=0.2em
\begin{array}{lcl}
	\Vert F(x^{k^{*}+1})\Vert & \leq & \frac{2\beta\tau}{\rho\sigma}\Vert \Delta x^{k^{*}}\Vert +\frac{1-\tau}{\rho}\left(2+\frac{M_F}{\sigma}\right)\Vert \Delta y^{k^{*}}\Vert +\frac{\tau }{\rho\sigma} \big( M_f+1+\Vert y^0\Vert \big) \vspace{1ex}\\
	& \leq & \frac{2\beta\tau}{\rho\sigma}\sqrt{\frac{8( \bar{\Pc} - \underline{\Pc} )}{\beta K}} + \frac{1-\tau}{\rho}\left(2+\frac{M_F}{\sigma}\right)\sqrt{\frac{2\rho( \bar{\Pc} - \underline{\Pc} )}{\tau(1-\tau)K}} + \frac{\tau}{\rho\sigma}  \big( M_f+1+\Vert y^0\Vert \big) \vspace{1ex}\\
	& \leq & \mathcal{O}\left(\frac{\sqrt{\tau\beta}}{\rho}\right)\frac{1}{\sqrt{K}}+ \mathcal{O}\left(\sqrt{\frac{1-\tau}{\tau^2\rho}}\right)\frac{1}{\sqrt{K}} + \dfrac{\epsilon}{2}.
\end{array}
\end{equation*}
This estimate shows that  for any $ \epsilon>0 $,  if $K\geq \max\big\{ \frac{\tau\beta}{\rho^2}, \sqrt{\frac{1-\tau}{\tau^2\rho}} \big\} \mathcal{O}\left(\frac{1}{\epsilon^2}\right)$, then we have  $\Vert F(x^{k^{*}+1})\Vert \leq \epsilon$.
Combining both conditions on $K$, we conclude that if
\begin{equation}\label{eq:K_bound}
K \geq  \max \left\{ \frac{\beta}{\tau},  \frac{\tau\beta}{\rho^2}, \sqrt{\frac{1-\tau}{\tau^2\rho}} \right\} \cdot \mathcal{O}\left(\frac{1}{\epsilon^2}\right),
\end{equation}
then $(x^{k^{*}+1}, y^{k^{*}+1})$ is an $\epsilon$-first-order optimal solution to \eqref{eq:nlp_gen}. Since $\rho = \mathcal{O}\left(\frac{\tau}{\epsilon}\right)$ from \eqref{eq:choice_of_rho}, by \eqref{eq:choice_of_beta}, we have $\beta \geq \mathcal{O}\left(\frac{\rho}{\tau^2}\right) = \mathcal{O}\big( \frac{1}{\tau\epsilon} \big)$.
Therefore, we can easily show that $\max \left\{ \frac{\beta}{\tau},  \frac{\tau\beta}{\rho^2}, \sqrt{\frac{1-\tau}{\tau^2\rho}} \right\} = \mathcal{O}\big( \frac{1}{\tau^2\epsilon}\big)$.
Substituting this into \eqref{eq:K_bound}, we can show that $K \geq \mathcal{O}\big( \frac{1}{\tau^2\epsilon^{3} }\big)$.
Overall, for a fixed value $\tau \in (0, 1]$, we can conclude that after $K = \mathcal{O}\big( \frac{1}{\epsilon^{3} }\big)$ iterations, $(x^{k^{*}+1}, y^{k^{*}+1})$ is an $\epsilon$-first-order optimal solution for problem \eqref{eq:nlp_gen}.
\Eproof
\endproof 
%%% End of Proof.

\medskip 

%%%% Remark 2.
\begin{remark}\label{re:compare_to_Lu22}
In contrast to \cite{Lu:22} where $\tau$ needs to be fixed a priori (larger than some  threshold), Theorem \ref{th:complex_bound1} shows that we are free to choose $\tau \in (0,1]$.  More specifically, if we set $\tau = \mathcal{O}(\epsilon^{\eta}) \in (0,1]$ with $\eta \in [0,1]$, then
\begin{equation*}
\rho \geq \frac{2\tau}{\sigma \epsilon}  \big( M_f + 1 + \Vert y^0\Vert \big) = \mathcal{O}\left(\frac{1}{\epsilon^{1 - \eta}}\right).
\end{equation*}
Therefore, for \(\eta = 0\) (hence, $\tau$ is independent of the desired accuracy $\epsilon$), the complexity of  Algorithm~\ref{alg:LIPAL} is $\mathcal{O}\left(\frac{1}{\epsilon^3}\right)$, which matches the one derived in \cite{Lu:22} for perturbed augmented Lagrangian methods and in \cite{BirGar:16} for quadratic penalty methods. 
For $\eta \in (0,1]$, the complexity increases to $\mathcal{O}\left(\frac{1}{\epsilon^{3 + 2\eta}}\right)$. 
In particular, when $\eta = 1$ (hence,  $\tau = \mathcal{O}(\epsilon)$), we find that $\rho = \mathcal{O}(1)$. 
Moreover, if $\bar{\Pc} - \underline{\Pc}$ is assumed to be independent of $1/\tau$, as usually considered in the literature, then our complexity becomes $\mathcal{O}\left(\frac{1}{\epsilon^{3 + \eta}}\right)$, which is better than the one in \cite{SahEft:19}.
\end{remark}
 
%%%%%%%%%%%%%%%%%%%%%%%%%%%%%%%%%%%%
%%%% 4.3. Improved rates under KL condition.
\subsection{Improved rates under KL condition.} \label{subsec:KL_rates}
In this subsection, we derive improved local convergence rates for Algorithm~\ref{alg:LIPAL} provided that the Lyapunov function defined in \eqref{eq:Lyapunov_func} satisfies a KL property. 
Let us first bound a subgradient $v \in\partial P(\cdot)$ of the Lyapunov function $P(\cdot)$.

%%% Lemma 9.  
\begin{lemma}[Boundedness of $v\in\partial P$]\label{le:bounded_grad} 
Let $\{ z^{k} \triangleq (x^{k},y^k) \}$ be  generated by  Algorithm \ref{alg:LIPAL} and  $P(\cdot)$ be the Lyapunov function defined  in \eqref{eq:Lyapunov_func}.  
If  $\{x^k\}$  belongs to a compact set $\mathcal{S}$ on which Assumption \ref{as:Assum} holds and $\beta$ is chosen as in Theorem \ref{th:complex_bound1}, then there exists a subgradient $v^{k+1}\in\partial P(x^{k+1},y^{k+1},y^{k};y^0)$ of $P$ such that for all $k \geq 0$, we have
\begin{equation}\label{eq:bounded_subgrad_of_P}
	\Vert v^{k+1}\Vert \leq c_1\Vert \Delta x^{k}\Vert + c_2\Vert \Delta y^{k}\Vert,
\end{equation}
where $c_1 \triangleq L^{\tau}_{\rho}+\beta+\rho M_F^2$ and $ c_2 \triangleq \frac{(8+\tau)(1-\tau)^2}{\tau\rho}+(1-\tau)M_F$.
\end{lemma}

%%% Proof of Lemma 9
\proof{Proof.} See Appendix.
\Eproof
\endproof
%%% End of Proof.

\medskip 

Let $\tau \in (0,1]$ be fixed and $\rho$ be chosen as in Theorem \ref{th:complex_bound1}.
Then, Lemma~\ref{le:bounded_grad} implies that
\begin{equation*} 
	\Vert v^{k+1}\Vert ^2\leq 2c_1^2\left(\Vert \Delta x^{k}\Vert ^2+\Vert \Delta y^{k}\Vert ^2\right). 
\end{equation*}
In addition, utilizing  \eqref{eq:decrease_Lyapunov}, we get
\begin{equation}\label{eq:lm10_decrease_Lyapunov_KL}
	\Pc_{k+1} - \Pc_{k} \leq -  \underline{\gamma} \left(\Vert \Delta x^{k}\Vert ^2+\Vert \Delta y^{k}\Vert ^2\right), \quad \textrm{where} \quad \underline{\gamma} := \frac{\tau(1-\tau)}{2\rho}.
\end{equation}
Combining the last two inequalities, we eventually get 
 \begin{equation}\label{eq:lm10_rate_of_v}
	\Pc_{k+1}-\Pc_{k}\leq-\frac{\ubar{\gamma}}{2c_1^2}\left\Vert v^{k+1}\right\Vert ^2.
 \end{equation}
Recall that we denoted $z^{k}=(x^{k},y^{k})$ and  $u^{k}=(x^{k},y^{k},y^{k-1})$. 
Moreover,  $\mathrm{crit}(P)$ denotes the set of  critical points of the Lyapunov function $P$ defined in \eqref{eq:Lyapunov_func} (i.e. $\mathrm{crit}(P) \triangleq \{ (x,y,y') \; : \;  0 \in \partial{P}(x,y,y';y^0) \}$). 
Furthermore, we denote $\mathcal{E}_{k}=\Pc_{k}-\Pc^{*}$, where $\Pc^{*}=\lim_{k\to\infty}{\Pc_{k}}$ (recall that since $\{\Pc_k\}_{k\geq1}$ is decreasing and bounded from below, it has the limit $P^{*}$). 
We also denote the set of limit points of $\{u^{k}\}_{k\geq1}$ by
\begin{equation*}
	\Omega:=\{ u^{*}   \, : \,  \textrm{there exists a subsequence  $\{u^{k}\}_{k \in \mathcal{K}}$ of $\{u^k\}$ such that $\lim_{k \in \mathcal{K},k \to \infty}{u^{k}} = u^{*}$} \}.
\end{equation*}
We prove the following result.

%%% Lemma 10.
\begin{lemma}\label{le:added_lemma}
Let $\{z^{k}:=(x^{k},y^{k})\} $ be  generated by  Algorithm \ref{alg:LIPAL} and $ P(\cdot)$ and $\{\Pc_{k}\}$ be defined as in \eqref{eq:Lyapunov_func}  and \eqref{eq:lyapunov_function}, respectively. 
Suppose that  $\{x^{k}\}$  belongs to a compact set $\mathcal{S}$ on which  Assumptions \ref{assump1}, \ref{as:Assum} and \ref{assump3} hold, $\rho$ is chosen as in \eqref{eq:choice_of_rho_b1}, $x^0$ is chosen as in \eqref{eq:choice_of_x0}, and $\beta$ is chosen as in Theorem \ref{th:complex_bound1}.
Then, the following  statements are valid:
\begin{enumerate} 
	\item[]$\mathrm{(i)}$~$\Omega$  is a compact subset of $\mathrm{crit}(P)$ and   $ \lim_{k\to\infty}{\mathrm{dist}(u^{k}, \mathrm{crit}(P))}=0$. \label{lem_item1}
	\item[]$\mathrm{(ii)}$~For any $u^{*} \in \Omega,$ we have $P(u^{*}) = \Pc^{*}$. \label{lem_item2}
	\item[]$\mathrm{(iii)}$~For any $(x^{*}, y^{*}, \hat{y}^{*} )\in \mathrm{crit}(P)$,   $(x^{*}, y^{*})$ is an $\epsilon$ - first-order optimal solution to \eqref{eq:nlp_gen} for $\epsilon := \frac{\tau}{\rho}\Vert y^{*} - y^0\Vert$. \label{lem_item3}
\end{enumerate}
\end{lemma}

%%% Proof of Lemma 10
\proof{Proof.} See Appendix.
\Eproof
\endproof
%%% End of Proof.

\medskip 

Our next step is to  prove that $\{ \Vert \Delta x^{k}\Vert +\Vert \Delta y^{k}\Vert  \}$ has finite length under the KL condition.

%%%% Lemma 11.
\begin{lemma}\label{le:finite_length}
Let $\{ z^{k}:=(x^{k},y^{k}) \}$ be  generated by  Algorithm \ref{alg:LIPAL} and let $P(\cdot)$ and $\{\Pc_{k}\}_{k\geq1}$ be defined as in \eqref{eq:Lyapunov_func} and \eqref{eq:lyapunov_function}, respectively.  
Suppose that $\{x^{k}\}$  belongs to a compact set $\mathcal{S}$ on which  Assumption \ref{as:Assum} holds,  $\beta$ is chosen as in Theorem \ref{th:complex_bound1}, and additionally, $P(\cdot)$ satisfies the KL property on $\Omega$.
Then, $\{z^{k}\}$  satisfies the following finite length property and $\{z^k\}$ converges to $z^{\star}$:
\begin{equation}\label{eq:lm11_finite_length}
\sum_{k=0}^{\infty}{\Vert \Delta x^{k}\Vert +\Vert \Delta y^{k}\Vert }<\infty \quad \textrm{and} \quad \lim_{k \to \infty} z^k=z^{*}:=(x^{*},y^{*}),
\end{equation}
such that $(x^{*},y^{*})$ is an $\epsilon$-first-order optimal solution to \eqref{eq:nlp_gen} with $\epsilon := \frac{\tau}{\rho}\Vert y^{*} - y^0\Vert$.
\end{lemma}

%%% Proof of Lemma 11
\proof{Proof.} See Appendix.
\Eproof
\endproof
%%% End of Proof.

\medskip 

Lemma \ref{le:finite_length} shows  that the set of  limit points of  $\{z^k=(x^{k},y^{k})\}_{k\geq1}$ is a singleton, denoted   $z^{*} =(x^{*},y^{*})$. 
Finally, let us bound $\Vert z^k-z^{*}\Vert$. 

%%% Lemma 12.
 \begin{lemma}\label{le:bound_of_diff_zk} 
 Let  $\{z^k:=(x^{k},y^{k})\}$ be  generated by Algorithm \ref{alg:LIPAL} and let $ P(\cdot)$ and $\{\Pc_{k}\}_{k\geq1}$ be defined by \eqref{eq:Lyapunov_func}  and \eqref{eq:lyapunov_function}, respectively. 
 Suppose that $\{x^{k}\}$  belongs  to a compact set $\mathcal{S}$ on which  Assumptions \ref{assump1}, \ref{as:Assum} and \ref{assump3} hold, $\rho$ is chosen as in \eqref{eq:choice_of_rho_b1}, $x^0$ is chosen as in \eqref{eq:choice_of_x0},  and $\beta$ is chosen as in Theorem \ref{th:complex_bound1}.
 Additionally,  assume that $P(\cdot)$  satisfies the KL property on $\Omega=\{u^{*} : u^{*}:=(x^{*},y^{*},y^{*})\}$.  
Then, there exists  $k_1\geq1$ such that for all $k\geq k_1$ we have 
\begin{equation}\label{eq:bound_of_diff_zk} 
	\Vert z^{k}-z^{*}\Vert \leq C\max\{\varphi(\mathcal{E}_{k}),\sqrt{\mathcal{E}_{k-1}}\},
\end{equation}
 where $C>0$ is a given constant  and  $\varphi\in\Psi_{\varepsilon}$ $($with $\varepsilon>0$$)$ denotes a desingularizing function. 
 \end{lemma}
 
%%% Proof of Lemma 12
\proof{Proof.} See Appendix.
\Eproof
\endproof
%%% End of Proof.

\medskip 

Finally, we prove the main result of this subsection, by showing local convergence rates $\{ z^k\}$ generated by Algorithm~\ref{alg:LIPAL} when the Lyapunov function $P$ satisfies the KL property with the following special desingularizing function:
\begin{equation}\label{eq:desingu_func}
	\varphi:[0,\varepsilon)\to[0,+\infty), \quad \varphi(s)=s^{1-\nu}, \quad \textrm{where}\quad  \nu\in[0,1).
\end{equation}
Note that this particular KL condition holds when $P(\cdot)$ is semi-algebraic, see, e.g., \cite{BolDan:07}.

%%% Theorem 2. 
\begin{theorem}[Convergence rates of $\{(x^{k},y^{k})\}$]\label{th:LIPAL_convergence_under_KL}  
Let  $\{z^k:=(x^{k},y^{k})\}$ be generated by Algorithm \ref{alg:LIPAL}  and let $ P(\cdot)$ and $\{\Pc_{k}\}_{k\geq1}$ be defined as in \eqref{eq:Lyapunov_func}  and \eqref{eq:lyapunov_function}, respectively. 
Suppose that $\{x^{k}\}$  belongs  to a compact set $\mathcal{S}$ on which  Assumptions \ref{assump1}, \ref{as:Assum} and \ref{assump3} hold, $\rho$ is chosen as in \eqref{eq:choice_of_rho_b1}, $x^0$ is chosen as in \eqref{eq:choice_of_x0}, and $ \beta$ is chosen as in Theorem \ref{th:complex_bound1}.
Suppose additionally that $P(\cdot)$ satisfies the KL property on $\Omega=\{ u^{*} \; : \; u^{*}:=(x^{*},y^{*},y^{*}) \}$ with  the desingularizing function \eqref{eq:desingu_func}.
Then, the following local convergence rates hold:
\begin{enumerate}
	\item[]$\mathrm{(i)}$~If $\nu=0$, then $\{ z^{k} \}$ converges to $z^{*}$ in a finite number of iterations.

	\item[]$\mathrm{(ii)}$~If $\nu\in(0,\frac{1}{2}]$, then there exists $k_1 \geq 1$ such that for all $k \geq k_1$, we have the following linear convergence rate:
	\begin{equation}\label{eq:th2_linear_rate}
         \Vert z^{k}-z^{*}\Vert \leq\frac{\sqrt{\mathcal{E}_{k_1}}}{(1+\bar{c}\mathcal{E}_{k_1}^{2\nu-1})^{\frac{k-k_1}{2}}},  \quad  \textrm{where} \  \bar{c} := \frac{\ubar{\gamma}}{2c_1^2}.
	\end{equation}
	
	\item[]$\mathrm{(iii)}$~If $\nu\in(\frac{1}{2},1)$, then there exists $k_1 \geq 1$ such that for all $k >  k_1 $, we have the following sublinear convergence rate:
         \begin{equation}\label{eq:th2_sublinear_rate}
         \Vert z^{k}-z^{*}\Vert \leq \Big(\frac{1}{\mu(k-k_1)+\mathcal{E}_{k_1}^{1-2\nu}} \Big)^{\frac{1-\nu}{2\nu-1}}.
         \end{equation}
       \end{enumerate}
 \end{theorem}
 
 %%% Proof of Theorem 2.
 \proof{Proof.} 
 Let $ \nu\in[0,1)$.
 For all  $s\in [0, \tau)$, since $\varphi(s)=s^{1-\nu}$, we have $\varphi'(s)=(1-\nu)s^{-\nu}$.  
 From Lemma \ref{le:bound_of_diff_zk}, it follows that for all $k \geq k_1$, we have
\begin{equation}\label{eq:th2_rate_point1}
	\Vert z^{k}-z^{*}\Vert \leq C \cdot \max\{\mathcal{E}_{k}^{1-\nu},\sqrt{\mathcal{E}_{k-1}}\}.
\end{equation}
Furthermore, for $\nabla{P}(x^{k},y^{k},y^{k-1};y^0) \in \partial{P}(x^{k},y^{k},y^{k-1};y^0)$, \eqref{eq:lm12_KL2} yields
\begin{equation*}
\mathcal{E}_k^{\nu}\leq \Vert \nabla{P}(x^{k},y^{k},y^{k-1};y^0)\Vert, \quad \forall k\geq k_1.
\end{equation*}
From \eqref{eq:lm10_rate_of_v} and Lemma \ref{le:bounded_grad}, there exists $v^k := \nabla{P}(x^{k},y^k,y^{k-1};y^0) \in \partial{P}(x^{k},y^k,y^{k-1};y^0)$ such that for all $k\geq 1$, we have
\begin{equation*}
        \Vert v^k\Vert ^2 = \Vert \nabla{P}(x^{k},y^k,y^{k-1};y^0)\Vert^2 \leq \frac{2c_1^2}{\ubar{\gamma}}(\mathcal{E}_{k-1}-\mathcal{E}_{k}).
\end{equation*}
Combining the last two relations, for all $k \geq k_1$, one has
\begin{equation}\label{eq:th2_key_rate_est}
	\bar{c} \mathcal{E}_k^{2\nu}\leq \mathcal{E}_{k-1}-\mathcal{E}_{k}, \quad \textrm{where} \ \bar{c} := \frac{\ubar{\gamma}}{2c_1^2}.
\end{equation}
We consider the following cases:
\begin{enumerate}
\item[$\mathrm{(i)}$] \textit{The case $\nu=0$.} 
If $\mathcal{E}_k>0$ for any $k \geq k_1$, then we have $\bar{c}\leq \mathcal{E}_{k-1}-\mathcal{E}_{k}$. 
As $k$ goes to infinity, the right hand side approaches zero. 
Then, $0<\bar{c}\leq0$, which is a contradiction. 
Hence,  there exists $k \geq k_1 $ such that $ \mathcal{E}_k=0$.
Then, $\mathcal{E}_k\to 0$ in a finite number of steps and  from \eqref{eq:th2_rate_point1}, $z^k\to z^{*}$ in a finite number of steps.

\medskip 

\item[$\mathrm{(ii)}$] \textit{The case $\nu\in(0,\frac{1}{2}]$.}
We have $2\nu-1\leq0$.
For $k \geq k_1$, since $\{\mathcal{E}_i\}_{i\geq k_1}$ is monotonically decreasing, $\mathcal{E}_i\leq\mathcal{E}_{k_1}$ for any $i\in\{k_1+1, k_1+2, \ldots, k\}$ and 
\begin{equation*}
	\bar{c}\mathcal{E}_{k_1}^{2\nu-1}\mathcal{E}_k\leq\mathcal{E}_{k-1}-\mathcal{E}_{k}.
\end{equation*}
Rearranging this expression, for all $k\geq k_1$, we get 
\begin{equation*}
\mathcal{E}_k\leq \frac{\mathcal{E}_{k-1}}{1+\bar{c}\mathcal{E}_{k_1}^{2\nu-1}}\leq\frac{\mathcal{E}_{k-2}}{(1+\bar{c}\mathcal{E}_{k_1}^{2\nu-1})^2}\leq \cdots  \leq\frac{\mathcal{E}_{k_1}}{(1+\bar{c}\mathcal{E}_{k_1}^{2\nu-1})^{k-k_1}}.
\end{equation*}
Therefore, we have $\max\{\mathcal{E}_k^{1-\nu},\sqrt{\mathcal{E}_{k-1}}\}=\sqrt{\mathcal{E}_{k-1}}$.
It then follows that
\begin{equation*}
         \Vert z^{k}-z^{*}\Vert \leq\frac{\sqrt{\mathcal{E}_{k_1}}}{\sqrt{(1+\bar{c}\mathcal{E}_{k_1}^{2\nu-1}})^{k-k_1}}.
\end{equation*}

\item[$\mathrm{(iii)}$] \textit{The case $\nu\in(1/2,1)$.}
For all $k \geq k_1$, from \eqref{eq:th2_key_rate_est} we have
\begin{equation}\label{eq:eqqq}
	\bar{c}\leq(\mathcal{E}_{k-1}-\mathcal{E}_k)\mathcal{E}_k^{-2\nu}.
\end{equation}
Let $h:\mathbb{R}_{+}\to\mathbb{R}$ be defined by $h(s)=s^{-2\nu}$ for any $s\in\mathbb{R}_{+}$. 
It is clear that $h$ is monotonically decreasing and for all $s\in\mathbb{R}_+$, we have $h'(s)=-2\nu s^{-(1+2\nu)}<0$. 
Since $\mathcal{E}_k\leq\mathcal{E}_{k-1}$ for all $k \geq k_1$, $h(\mathcal{E}_{k-1})\leq h(\mathcal{E}_k)$ for all $k \geq k_1$. 
We consider two subcases:

%%% Case 1.
\noindent\textbf{Case 1.} 
Let $r_0\in(1,+\infty)$ such that $ h(\mathcal{E}_k)\leq r_0h(\mathcal{E}_{k-1})$ for all $k \geq k_1$.
Then, from  \eqref{eq:eqqq} we get
\begin{equation*}
\arraycolsep=0.2em
\begin{array}{lcl}
	\bar{c} & \leq & r_0(\mathcal{E}_{k-1}-\mathcal{E}_k)h(\mathcal{E}_{k-1}) \leq  r_0h(\mathcal{E}_{k-1})\int_{\mathcal{E}_k}^{\mathcal{E}_{k-1}}{1\,ds} \vspace{1ex} \\
	& \leq & r_0\int_{\mathcal{E}_k}^{\mathcal{E}_{k-1}}{h(s)\,ds}= r_0\int_{\mathcal{E}_k}^{\mathcal{E}_{k-1}}{s^{-2\nu}\,ds}=\frac{r_0}{1-2\nu}(\mathcal{E}_{k-1}^{1-2\nu}-\mathcal{E}_{k}^{1-2\nu}).
\end{array}    
\end{equation*}
Since $\nu>\frac{1}{2}$, the last inequality leads to
\begin{equation*}
	0<\frac{\bar{c}(2\nu-1)}{r_0}\leq \mathcal{E}_{k}^{1-2\nu}-\mathcal{E}_{k-1}^{1-2\nu}.
\end{equation*}
Let us define $\hat{c}=\frac{\bar{c}(2\nu-1)}{r_0}$ and $\hat{\nu}=1-2\nu<0$. 
We get
\begin{equation}\label{eq:need1}
	0<\hat{c}\leq \mathcal{E}_{k}^{\hat{\nu}}-\mathcal{E}_{k-1}^{\hat{\nu}}, \quad \forall k geq k_1.
\end{equation}

%%% Case 2.
\noindent\textbf{Case 2.} 
Let $r_0\in(1,+\infty)$ such that $h(\mathcal{E}_k) >  r_0h(\mathcal{E}_{k-1})$ for all $k \geq k_1$. 
Then, we have $\mathcal{E}_k^{-2\nu}\geq r_0 \mathcal{E}_{k-1}^{-2\nu}$. 
This inequality leads to 
\begin{equation*}
q\mathcal{E}_{k-1}\geq\mathcal{E}_k, \quad \textrm{where} \ q={r_0}^{-\frac{1}{2\nu}}\in(0,1).
\end{equation*}
Since $\hat{\nu}=1-2\nu<0$, we have $ q^{\hat{\nu}} \mathcal{E}_{k-1}^{\hat{\nu}}\leq\mathcal{E}_k^{\hat{\nu}}$, which is equivalent to
\begin{equation*}
(q^{\hat{\nu}}-1)\mathcal{E}_{k-1}^{\hat{\nu}}\leq\mathcal{E}_{k-1}^{\hat{\nu}}-\mathcal{E}_{k}^{\hat{\nu}}.
\end{equation*}
Since $q^{\hat{\nu}}-1>0$ and $\mathcal{E}_k \to 0^+$ as $k\to\infty$, there exists $\Tilde{c}$ such that $(q^{\hat{\nu}}-1)\mathcal{E}_{k-1}^{\hat{\nu}}\geq\tilde{c}$ for all $k \geq k_1$. 
Therefore, we can show that
\begin{equation}\label{eq:need2}
        0<\tilde{c}\leq \mathcal{E}_{k}^{\hat{\nu}}-\mathcal{E}_{k-1}^{\hat{\nu}}, \quad \forall k \geq k_1.
\end{equation}
By choosing $\mu :=\min\{\hat{c},\tilde{c}\}>0$, one can combine \eqref{eq:need1} and \eqref{eq:need2} to obtain
\begin{equation*}
0<{\mu}\leq \mathcal{E}_{k}^{\hat{\nu}}-\mathcal{E}_{k-1}^{\hat{\nu}}, \quad \forall k \geq k_1.
\end{equation*}
Summing this inequality from $k_1 + 1$ to some $k > k_1$, we get
\begin{equation*}
	\mu(k-k_1) + \mathcal{E}_{k_1}^{\hat{\nu}}\leq  \mathcal{E}_{k}^{\hat{\nu}}.
\end{equation*}
Hence, we can easily show that
\begin{equation*}
\mathcal{E}_k\leq (\mu(k-k_1) + \mathcal{E}_{k_1}^{\hat{\nu}})^{\frac{1}{\hat{\nu}}}=(\mu(k-k_1) + \mathcal{E}_{k_1}^{1-2{\nu}})^{\frac{1}{1-2{\nu}}}.
\end{equation*}
Since $\nu\in(\frac{1}{2},1)$, we have $\max\{{\mathcal{E}_{k-1}}^{1-\nu},\sqrt{\mathcal{E}_{k-1}}\}={\mathcal{E}_{k-1}}^{1-\nu}$.
Therefore, \eqref{eq:th2_rate_point1} becomes 
\begin{equation*}
\Vert z^{k}-z^{*}\Vert \leq\left(\frac{1}{\mu(k-k_1)+\mathcal{E}_{k_1}^{1-2\nu}}\right)^{\frac{1-\nu}{2\nu-1}}, \quad \forall k \geq k_1,
\end{equation*}
which shows a sublinear convergence rate of $\{ \Vert z^k - z^{\star}\Vert\}$.
\end{enumerate} 
Putting three cases together, we complete our proof.
\Eproof
\endproof 
%%%% End of Proof.

\medskip 

Note that $z^{*}=(x^{*},y^{*})$ from Theorem \ref{th:LIPAL_convergence_under_KL}  is $\epsilon$-first-order optimal solution to \eqref{eq:nlp_gen} with $\epsilon = \frac{\tau}{\rho}\Vert y^{*}-y^0\Vert$.
Therefore, this point is also an $\epsilon$-first-order optimal solution to \eqref{eq:nlp_gen} if we choose $\rho$ as in Theorem \ref{th:complex_bound1} (i.e., $\rho = \mathcal{O}\big( \frac{\tau}{\epsilon}\big)$ is of order $\frac{\tau}{\epsilon}$).

%%%%%%%%%%%%%%%%%%%%%%%%%%%%%%%%%%%%%%%%%%%%%%
%%% 4.4. Selection of the penalty parameter $\rho$.
%%%%%%%%%%%%%%%%%%%%%%%%%%%%%%%%%%%%%%%%%%%%%%
\subsection{Selection of the penalty parameter $\rho$.}\label{subsec:selection_of_rho}
The convergence results from Sections \ref{subsec:globa_conver} and \ref{subsec:KL_rates} estimate the total number of iterations required to compute an $\epsilon$-first-order optimal solution to \eqref{eq:nlp_gen}.
These complexity results rely on the assumption that the penalty parameter $\rho$ exceeds a certain threshold. 
In fact, for any $\rho > 0$, we can always guarantee the optimality condition.
However, to guarantee the feasibility, we need to choose the penalty parameter $\rho$ sufficiently large (as seen in Theorem \ref{th:complex_bound1}).
Clearly, determining this threshold of $\rho$ beforehand poses some challenges since it depends on unknown parameters rendered from \eqref{eq:nlp_gen} as well as the algorithm's parameters.  

\medskip 

To overcome this challenge, we propose a variant of Algorithm~\ref{alg:LIPAL} which allows us to determine a sufficiently large value of $\rho$ without knowning specific parameter information. 
Inspired by Algorithm 3 in \cite{XieWri:21}, our method repeatedly employs Algorithm \ref{alg:LIPAL} as an inner loop, and calls it $S$ iterations (called stages). 
If Algorithm \ref{alg:LIPAL} fails to converge within a predetermined number of iterations, we gradually increase the penalty parameter $\rho$ by multiplying it by a constant at each outer iteration. 
The detailed implementation of this scheme is outlined in Algorithm \ref{alg:LIPAL_v2}. 

%%% Algorithm 2.
\begin{algorithm}
\caption{(LIPAL with trial values of $\rho$)}\label{alg:LIPAL_v2}
\begin{algorithmic}[1]
\State \textbf{Initialization:} Choose $x^0 \in \dom g$ and $y^0\in\mathbb{R}^m$.
\State Choose $\delta_1>1$, $0 < \delta_2 < 1$, $\epsilon>0$, $\tau_0\in(0,1]$, and $\rho_0 > 0$.
\For{$s = 0, 1, \cdots, S$}
\State If \textrm{\textbf{Feasible}} then \textbf{terminate}.
\State Call \textbf{Algorithm \ref{alg:LIPAL}} with $\tau_s$, $\rho_s$ and $\beta_s$, where $\beta_s$ is determined by a backtracking procedure, 
\Statex {~~~~} and  using warm start  $(x_s^0, y_s^0)\gets(x_{s-1}^{*}, y_{s-1}^{*})$ until  $ \mathrm{dist} \big(-\nabla f(x_s^{*})-{J_F(x_s^{*})}^Ty_s^{*},\partial g(x_s^{*}) \big) \leq \epsilon$.
\State $\rho_{s+1}\gets \delta_1 \rho_s$
\State $\tau_{s+1}\gets \delta_2 \tau_s$
\EndFor
\end{algorithmic}
\end{algorithm}

Note that in Algorithm~\ref{alg:LIPAL_v2}, $(x_s^0, y_s^0)$ denotes the initial primal and dual iterate for Algorithm \ref{alg:LIPAL} at the $s^{\text{th}}$ stage of Algorithm \ref{alg:LIPAL_v2}, while $(x_{s-1}^{*},y_{s-1}^{*})$ denotes the last primal and dual iterate (or the output) yielded by Algorithm \ref{alg:LIPAL} at the $(s-1)^{\text{th}}$ stage of Algorithm \ref{alg:LIPAL_v2}. 
At $s=0$ (the first stage), $(x_0^0,y_0^0)$ can be chosen arbitrarily. 

\medskip 

Algorithm \ref{alg:LIPAL_v2} is well-defined and is terminated after a finite number of iterations. 
In fact, during the $s$-th stage of Algorithm \ref{alg:LIPAL_v2}, we have $\rho_s = \delta_1^s \rho_0$ and $\tau_s = \delta_2^s \tau_0$. 
Consequently,  from Theorem \ref{th:complex_bound1}, we have $\Vert F(x_s^{*})\Vert  \leq \frac{\tau_{s}(M_f + 1 + \Vert y_s^0\Vert )}{\rho_{s}\sigma} $, leading to
\begin{equation*}
\Vert F(x_s^{*})\Vert  \leq \frac{\tau_s (M_f+1+\Vert y_s^0\Vert )}{\rho_s \sigma} \leq \frac{\tau_0(M_f+1+\Vert \bar{y}^0\Vert )}{\rho_0}\left(\frac{\delta_2}{\delta_1}\right)^s. 
\end{equation*}
Hence, if $\frac{\tau_0(M_f+1+\Vert \bar{y}^0\Vert )}{\rho_0}\left(\frac{\delta_2}{\delta_1}\right)^s \leq \epsilon$, then $\Vert F(x_s^{*})\Vert  \leq \epsilon$.
This condition holds if we choose $s$ such that
\begin{equation*}
s \log\left(\frac{\delta_1}{\delta_2}\right) \geq \log\frac{\tau_0(M_f+1+\Vert \bar{y}^0\Vert )}{\rho_0} + \log\frac{1}{\epsilon},
\end{equation*}
which leads to
\begin{equation*}
s \geq S := \left\lfloor \frac{\log(\tau_0) + \log(M_f + 1 + \Vert \bar{y}^0\Vert ) - \log(\rho_0) - \log(\epsilon) }{ \log(\delta_1) - \log(\delta_2) } \right\rfloor + 1.
\end{equation*}
Therefore, we conclude that one has to increase $\rho$ (and/or decrease $\tau$) at most $S$ times (for $S$ defined above) to reach a value that guarantees the desired precision for feasibility and consequently to yield an $\epsilon$-first-order optimal solution to \eqref{eq:nlp_gen}.

\section{Discussion on Regularity Conditions.}\label{sec_reg}
In this section, we examine our new  regularity condition stated in Assumption \ref{as:Assum}-\textrm{(iii)}, which played  a crucial role in our convergence analysis, in particular for  bounding the dual iterates $\{y^k\}$. 
We first discuss some special cases of our assumption. Then, we explore an algorithmic perspective to understand the necessity of  such regularity condition, highlighting how the absence of this condition  leads to infeasibility.

%%% 3.1. Regularity condition for special cases
\subsection{Special cases of the regularity condition.}
In this subsection, we provide several examples of problems for which   Assumption \ref{as:Assum}-\textrm{(iii)} is valid.

\medskip 

\noindent\textbf{Case 1 -- Lipschitz function $g$.}
If the nonsmooth term $g$ of \eqref{eq:nlp_gen} is Lipchitz continuous on $\mathcal{S}$ (e.g., $g$ is a norm or an affine max-type function), then we have $\partial^{\infty}g(x)=\{0\}$ for all $x\in \mathcal{S}$, see  \cite{RocWet:98}. Hence, Assumption \ref{as:Assum}-\textrm{(iii)} becomes
\begin{equation*}
\sigma\left\Vert F(x)\right\Vert \leq \Vert J_F(x)^TF(x)\Vert,  \quad  \forall x\in \mathcal{S}.
\end{equation*}
Note that this condition always holds provided that  $J_F(x)$ has full row rank for all $x\in\mathcal{S}$,  equivalently, LICQ holds on $\mathcal{S}$. 
Clearly, LICQ  is frequently used in the literature when the objective function is smooth or Lipchitz continuous \cite{XieWri:21,Ber:96}.
Consequently, our assumption recovers this particular setting.

\medskip 

\noindent\textbf{Case 2 -- Indicator function $g$.}
If the nonsmooth term $g$ in \eqref{eq:nlp_gen} is the indicator function of a convex set $\Omega$, then we have $\partial^{\infty}g(x)=N_{\Omega}(x)$ for any $x \in \Omega$. Hence, Assumption \ref{as:Assum}-\textrm{(iii)} reduces to
\begin{equation*}
\sigma\Vert F(x)\Vert \leq\mathrm{dist}\left(-J_F(x)^TF(x),N_{\Omega}(x)\right), \quad \text{for all } x\in\mathcal{S}\subseteq \Omega.
\end{equation*}
This condition holds e.g., when $F$ is an affine map and $\Omega$ is a polyhedral set (see \cite{LiChe:21}). 
Recall that for any proper, lsc, convex function $g$, we have $\partial^{\infty}g(x)=N_{\dom g}(x)$. 
Furthermore, a similar constraint qualification condition was considered in \cite{Lu:22} in the context of inequality constraints.

\medskip 

\noindent\textbf{Case 3 -- Separable objective and constraints.}
If $g(x) = h(x_1)$ and $F(x) = G(x_1) + H(x_2)$ in \eqref{eq:nlp_gen},  where $x = (x_1, x_2) \in \mathbb{R}^{n_1} \times \mathbb{R}^{n_2}$, with $J_H$ having full row rank on a set $\mathcal{S}_2\subseteq\mathbb{R}^{n_2}$ and $h$ being a proper lsc function, then Assumption \ref{as:Assum}-\textrm{(iii)} holds on $\dom h \times \mathcal{S}_2$.
Indeed, we have
\begin{equation*}
\begin{array}{lcl}
	\mathrm{dist}^2\left(-J_F(x)^TF(x), \partial^{\infty} g(x)\right) 
	 & = &  \mathrm{dist}^2\left(-\begin{bmatrix} J_G(x_1)^T \\ J_H(x_2)^T \end{bmatrix} F(x), \begin{bmatrix}  \partial^{\infty} h(x_1) \\ 0_{n_2} \end{bmatrix}\right) \vspace{1ex}\\
    	& = & \mathrm{dist}^2\left(-J_G(x_1)^TF(x), \partial^{\infty} h(x_1)\right) + \mathrm{dist}^2\left(-J_H(x_2)^TF(x), 0_{n_2} \right) \vspace{1ex} \\
	& \geq & \left\Vert -J_H(x_2)^TF(x)\right\Vert ^2 \vspace{1ex}\\
	& \geq & \sigma^2 \Vert F(x)\Vert ^2,
\end{array}    
\end{equation*}
where the last inequality follows from the full-row rank of $J_H(x_2)$ on $\mathcal{S}_2$, or equivalently,  there exists $\sigma>0$ such that $\Vert J_H(x_2)^T y\Vert \geq \sigma \Vert y\Vert $ for any $y\in\mathbb{R}^m$ and $x_2 \in \mathcal{S}_2$.

\medskip 

\noindent\textbf{Comparison with \cite{SahEft:19}.} 
Next, we compare  Assumption \ref{as:Assum}-\textrm{(iii)} with the regularity condition introduced in \cite{SahEft:19}  for the problem \eqref{eq:nlp_gen}, but assuming $g$ to be convex, which is formulated as follows:
\begin{equation}\label{eq:reg_Volkan}
\sigma\Vert F(x^k)\Vert \leq\mathrm{dist}\left(-J_F(x^k)^TF(x^k),\ \rho_k^{-1} \partial g(x^k) \right), \quad \forall k\geq 0,
\end{equation}
where $\rho_k$ denotes the penalty parameter in a pure augmented Lagrangian setting, which should increase along the iterations, and $\partial g(x^k)$ denotes the  subdifferential of the convex function $g$ at the iterate $x^k$. It was argued in \cite{SahEft:19} that when $g=0$, this condition simplifies to:
\[
\sigma\Vert F(x^k)\Vert \leq\Vert J_F(x^k)^TF(x^k)\Vert ,
\]
which is always true provided  $J_F(x^k)$ has full row rank on the level set where the iterates belong, see also \textit{Case 1} above. 
However, the regularity condition \eqref{eq:reg_Volkan}  has some drawbacks: (i) First, this condition is given in terms of  the iterates $x_k$ and the penalty parameter $\rho_k$ of the algorithm instead of relying solely on the problem's data.  Hence, this condition cannot be checked a priori. Moreover, we argue that the horizon subdifferential of the nonsmooth objective function is the proper tool to be used instead of the  subdifferential. Indeed, with the horizon subdifferential, their condition would no longer depend on the iterates or the parameter $\rho_k$, as the horizon subdifferential is a  cone, rendering the penalty parameter $\rho_k$ obsolete in \eqref{eq:reg_Volkan}.
(ii) Second, we argue that \eqref{eq:reg_Volkan} does not have any meaning  when $g$ is locally Lipschitz or  locally smooth. Indeed, in this scenarios an adequate constraint qualification should not involve $g$. However, our constraint qualification condition still makes sense in this setting, see \textit{Case 1}.

\medskip 

\noindent\textbf{Constraint qualifications of \eqref{eq:nlp_gen} and their relation to Assumption~\ref{as:Assum}-\textrm{(iii)}.}  
\noindent Finally, we introduce a constraint qualification condition for problem \eqref{eq:nlp_gen} at a feasible point $\bar{x}$,  commonly used in the literature, such as  \cite{GuoYe:18, HalTeb:23, RocWet:98} and make some connections  with  Assumption \ref{as:Assum}-\textrm{(iii)}.  

\medskip 

%%% Definition 7.
\begin{definition} \label{de:regularity_def}
Let $\bar{x} \in \mathbb{R}^n$ be a feasible point to \eqref{eq:nlp_gen}. 
We say that problem \eqref{eq:nlp_gen} is regular at $\bar{x}$ if the following constraint qualification condition holds:
\begin{equation}\label{CQ_regularity}
	\mathrm{dist}\left(-J_F(\bar{x})^T y,\partial^{\infty} g(\bar{x})\right)=0, \quad y\in N_{\{0_m\}}(F(\bar{x}))=\mathbb{R}^m \implies y=0.
\end{equation}
\end{definition}

\medskip 

This regularity notion extends LICQ at a point $\bar{x}$ to problems of the form \eqref{eq:nlp_gen} with non-Lipschitzian objective. 
If  $\bar{x}$ is a local minimizer of \eqref{eq:nlp_gen}, which satisfies the constraint qualification condition \eqref{CQ_regularity}, then $\bar{x}$ is a KKT point of  \eqref{eq:nlp_gen} (see Proposition 6.9 in \cite{KruMeh:22}).  

\medskip 

Now, let us define a class of functions that satisfies the constraint qualification \eqref{CQ_regularity} at a point or on a set. 
This class of functions are called amenable functions, see \cite{RocWet:98}.

\medskip 

%%% Definition 8.
\begin{definition}[Amenable  functions]\label{def:amenable}
 A function $\varphi: \mathbb{R}^n \to \mathbb{R}$ is \textit{amenable} at $\bar{x}$ if $\varphi$ is finite at $\bar{x}$ and there is an open neighborhood $V$ of $\bar{x}$ on which $\varphi$ can be represented as $\varphi = h \circ G$, where
\begin{itemize}
\item[]$\mathrm{(i)}$ $G : V\to \mathbb{R}^m\) is a $C^1$ mapping;  
\item[]$\mathrm{(ii)}$ $h : \mathbb{R}^m \to \mathbb{R}$ is a proper, lsc, convex function; 
\item[]$\mathrm{(iii)}$ given $D = \operatorname{cl}(\dom h)$, the only vector $y \in N_D\left( G(\bar{x})\right)$ with $ J_G(\bar{x})^T y = 0$ is $y = 0$.
\end{itemize}
Moreover, $\varphi$ is said to be amenable on a set $\mathcal{S}$ if it is amenable at all points of $\mathcal{S}$.
\end{definition}

\medskip 

From Definition~\ref{def:amenable}, it follows that if $g$ in \eqref{eq:nlp_gen} is convex and the constraint qualification \eqref{CQ_regularity} holds at a feasible point $\bar{x}$, then the function  
\begin{equation*}
\varphi(x) := f(x) + g(x) + \delta_{\{0_m\}}(F(x)) = f(x) + h(\tilde{F}(x))
\end{equation*}
is amenable at $\bar{x}$, where $\tilde{F}(x) = \begin{pmatrix} F(x) \\ x \end{pmatrix}$ is a $C^1$-mapping from $\mathbb{R}^n$ to $\mathbb{R}^{m+n}$ and  
$h(y,x) = \delta_{\{0_m\}}(y) + g(x)$ is a proper, lsc, and convex function.    
Moreover, since $f$ is a $C^1$ function, if $h \circ \tilde{F}$ is amenable at $\bar{x}$, then $\varphi$ is amenable at $\bar{x}$ (by the calculus of amenability). 
Hence, it remains to verify that $h \circ \tilde{F}$ is amenable at $\bar{x}$, i.e.,  to show that for $D = \text{cl}(\dom  h) = \{0_m \times \text{cl}(\dom  g)\}$, the only vector  
\begin{equation*}
w = \begin{pmatrix} y \\ z \end{pmatrix} \in N_D(\tilde{F}(\bar{x})) = \mathbb{R}^m \times N_{\dom g} (\bar{x})
\end{equation*}
satisfying $J_{\tilde{F}}(\bar{x})^T w = 0$ is $w = 0$.  
Indeed, let $w \in \mathbb{R}^m \times N_{\dom  g}(\bar{x})$ such that  
\begin{equation} \label{reg_amenable}
J_{\tilde{F}}(\bar{x})^T w = 0 \iff 
J_{\tilde{F}}(\bar{x})^T w = \begin{pmatrix} J_F(\bar{x}) \\ I_n \end{pmatrix}^T \begin{pmatrix} y \\ z \end{pmatrix} = J_F(\bar{x})^T y + z=0.
\end{equation}
Since $g$ is convex and $z\in N_{\dom  g}(\bar{x})$, we also have
\begin{equation*}
	\mathrm{dist}\left(-J_F(\bar{x})^T y, \partial^\infty g(\bar{x})\right) =\mathrm{dist}\left(-J_F(\bar{x})^T y, N_{\dom g}(\bar{x})\right) \leq \Vert  -J_F(\bar{x})^T y - z \Vert  = 0.
\end{equation*}
Furthermore, since the constraint qualification \eqref{CQ_regularity} holds at $\bar{x}$, it follows that $y = 0$. Substituting $y = 0$ in \eqref{reg_amenable}, we find  $z = 0$.  
Therefore, $\varphi$ is amenable at $\bar{x}$ when \eqref{CQ_regularity} holds.

\medskip 

It is easy to see that the constraint qualification \eqref{CQ_regularity} is actually equivalent to the constraint qualification defining amenablility of $\varphi$ (see condition (iii) in Definition \ref{def:amenable}). 
Moreover, since $\varphi$ is amenable at $\bar{x}$, it follows that there exists an open neighborhood $V$ of $\bar{x}$ such that $\varphi$ is amenable at every point $x \in V \cap \dom \varphi$. 
Hence, we have
\begin{equation}\label{CQ_regularity_neighboor}
	\mathrm{dist}\left(-J_F(x)^T y,\partial^{\infty} g(x)\right)=0, \quad y\in\mathbb{R}^m \implies y=0 \quad \forall x \in V \cap \dom  \varphi.
\end{equation}
Note that if there exists $\sigma > 0$ such that, for any $y\in \mathbb{R}^m$, we have
\begin{equation}\label{ours_local}
	\sigma\Vert y\Vert  \leq \mathrm{dist}\left(-J_F(x)^T y,\partial^{\infty} g(x)\right) \quad \forall x \in V \cap \dom  \varphi,
\end{equation}
then \eqref{CQ_regularity_neighboor} follows. Moreover, the equivalence between \eqref{CQ_regularity_neighboor} and  \eqref{ours_local} holds if additionally $J_F(x)$ has full row rank for any $x \in V \cap \dom  \varphi$. 
Note also that the  full row rank condition is a necessary condition for \eqref{ours_local} to be satisfied. 
This is proved in the next lemma.  

\medskip 

%%% Lemma 18.
\begin{lemma}\label{lemma18}
Assume that $g(\cdot)$ is convex and $J_F(\cdot)$ is continuous and has full row rank on $V \cap \dom \varphi$.
Then, \eqref{CQ_regularity_neighboor} is valid if and only if \eqref{ours_local} holds for some $\sigma>0$.
\end{lemma}

%%% Proof of Lemma 18
\proof{Proof.} See Appendix.
\Eproof
\endproof
%%% End of Proof.

\medskip

Note that \eqref{CQ_regularity_neighboor} (which is equivalent to \eqref{ours_local} under a full row rank condition) is assumed to hold for feasible points (i.e., $F(x)=0$). 
By drawing an analogy with the extension of LICQ  for problems with smooth and/or Lipschitz continuous objective functions to infeasible points, one would require the following condition for problems with non-Lipschitz objective functions:
\begin{equation}\label{extension_LICQ}
\sigma\Vert y\Vert  \leq \mathrm{dist}\left(-J_F(x)^T y, \partial^{\infty} g(x)\right) \quad \forall x \in \mathcal{S}, \; y \in \mathbb{R}^m,
\end{equation}
where usually $\mathcal{S}$ should be larger than $V\cap \dom \varphi$ in order to ensure global convergence for an algorithm. 
However, requiring \eqref{extension_LICQ} on such larger set $\mathcal{S}$ may be too restrictive, as there may exist $y \neq 0$ such that $-J_F(x)^T y \in \partial^{\infty} g(x)$. 
On the other hand, assuming that \eqref{extension_LICQ} holds not for all $y \in \mathbb{R}^m$ but only for a specific choice of $y$, i.e. $y = F(x)$, then  \eqref{extension_LICQ} becomes Assumption \ref{as:Assum}-\textrm{(iii)} and this appears to be less restrictive.
Indeed, if $g$ is convex,  $F(x) = Ax - b$, and the problem \eqref{eq:nlp_gen} is well-posed in the sense that there exists $v \in \dom g$ such that $Av - b = 0$, then for any infeasible point $x \in \dom g$ (i.e. $Ax - b \neq 0$), we have
\begin{equation*}
	\langle -J_F(x)^TF(x), v - x \rangle =  \langle -A^T(Ax-b), v - x \rangle = \langle Ax - b, Ax - Av \rangle = \Vert Ax - b\Vert ^2 > 0.
\end{equation*}
Hence, for an infeasible point $x$, we always have $-A^T(Ax - b) \notin \partial^{\infty} g(x)= N_{\dom g}(x)$. 
Consequently,  there may exist a $\sigma>0$ for which Assumption \ref{as:Assum}-\textrm{(iii)} holds, while \eqref{extension_LICQ} fails to hold.

%%%%%%%%%%%%%%%%%%%%%%%%
%%% 6.2. Algorithmic perspective.
%%%%%%%%%%%%%%%%%%%%%%%%
\subsection{Algorithmic perspective.}
From an algorithmic standpoint, augmented Lagrangian and penalty-type methods usually guarantee that any limit point of the iterate sequences  satisfies the stationarity (optimality) condition and the KKT condition corresponding to the feasible problem, see, e.g., \cite{DemJia:23}. 
However, the feasibility at this limit point may not hold. 
Therefore, one needs to impose a regularity condition to enforce the feasibility at such limit point. For nonlinear programs with a smooth or Lipschitz continuous objective function,  LICQ condition is sufficient to ensure that any KKT point of the feasible problem is feasible for the original problem. However, when the objective function is  non-Lipschitz continuous, the horizon subdifferential of such an objective function is not reduced to the singleton $\{0\}$. 
Therefore, the horizon subdifferential must be considered in constructing an appropriate constraint qualification condition to enforce the feasibility. 
Consequently, our regularity condition given in Assumption \ref{as:Assum}-\textrm{(iii)} represents a generalization of LICQ  constraint qualification, extending it to accommodate non-Lipschitz objective functions from an algorithmic perspective.  

\medskip 

In the following lemma, we prove that in the absence of  Assumption \ref{as:Assum}-\textrm{(iii)}, any limit point of our LIPAL algorithm is  only an $\epsilon$-first-order solution of the associated feasibility problem.

%%% Lemma 19.
\begin{lemma}[Limit points are $\epsilon$-KKT points of feasible problem]\label{le:limit_feasibility}
Let  $\{(x^k,y^k)\}_{k\geq0}$ be generated by Algorithm \ref{alg:LIPAL}. 
Suppose that Assumptions \ref{assump1}, \ref{as:Assum}-$\mathrm{(i)}$, \ref{as:Assum}-$\mathrm{(ii)}$ and \ref{assump3} hold on some compact set $\mathcal{S}$  on which the primal sequence $\{x^k\}_{k\geq0}$  belongs to,  $x^0$ is chosen as in \eqref{eq:choice_of_x0}, $\tau\in(0,1)$, $\epsilon>0$,  $\rho$, and $\beta$ are chosen as
\begin{equation*}
	\rho \geq \max\left\{1, 3\rho_0, \frac{2\tau (M_f+1+\Vert y^0\Vert )}{\sigma\epsilon}\right\} \quad\textrm{and} \quad  \beta\geq\max\left\{2L_\rho^\tau, \frac{8(1-\tau)\rho M_F^2}{\tau^2}\right\}.
\end{equation*}
Moreover, let  $\{x^k\}_{k\in\mathcal{K}}$ be a subsequence  satisfying  $x_k \xrightarrow[]{g} x^{*}$. 
Then,  $x^{*}$ is an $\epsilon$-first-order solution of the following feasible problem:
\begin{equation*}
	\min_{(x, \alpha) \in \emph{epi}(g) }   \tfrac{1}{2} \Vert F(x) \Vert ^2.
\end{equation*}
Additionally, if $g$ is locally Lipschitz continuous at $x^{*}$, then $x^{*}$ is an $\epsilon$-first-order solution of the following feasible problem:
\begin{equation*}
	\min_{x \in \mathbb{R}^n } \tfrac{1}{2} \Vert F(x) \Vert ^2.
\end{equation*}
\end{lemma}

%%% Proof of Lemma 19
\proof{Proof.} See Appendix.
\Eproof
\endproof
%%% End of Proof.

\medskip 

\noindent Note that Lemma \ref{le:limit_feasibility} shows that when the regularity  Assumption \ref{as:Assum}-\textrm{(iii)} is not used, Algorithm \ref{alg:LIPAL} can only ensure converge to an $\epsilon$-first-order solution  of the feasible problem, but  we cannot ensure that $\Vert F(x^{*})\Vert \leq \epsilon$. However, when Assumption \ref{as:Assum}-\textrm{(iii)} holds, then we proved in Theorem \ref{th:complex_bound1} that any limit point $x^{*}$ is an $\epsilon$-first-order solution of the original problem \eqref{eq:nlp_gen}, i.e., we can ensure  $\Vert F(x^{*})\Vert \leq \epsilon$.

%%%%%%%%%%%%%%%%%%%%%%%%%%%%%%%%%%%%%%%%%%%%%%%
%%% 7. Numerical Experiments.
%%%%%%%%%%%%%%%%%%%%%%%%%%%%%%%%%%%%%%%%%%%%%%% 
\section{Numerical experiments.}\label{sec5}
%In this section, we provide detailed numerical experiments to verify our algorithm and compare it with existing methods.

%\qcmt{
%\subsection{Mathematical program with complementarity constraints.}\label{subsec:MPCC}
%In this experiment, we consider the following mathematical program with complementarity constraint:
%\begin{equation}\label{eq:MPCC} \arraycolsep=0.2em
%\left\{\begin{array}{ll}
%{\displaystyle\min_{u, v}} &f(u, v), \vspace{1ex}\\
%\mathrm{s.t.} & (u, v) \in \mathcal{S}, \vspace{1ex}\\
%& u \geq 0, \ Au + Bv + c \geq 0, \ u^{\top}(Au + Bv + c) = 0,
%\end{array}\right.
%\end{equation}
%where $f : \R^{n_1}\times\R^{n_2} \to \R$ is convex, $\mathcal{S}$ is a nonempty, closed, and convex set in $\R^n$ for $n:=n_1+n_2$, $A$ and $B$ are given matrices with consistent dimensions, and $c$ is a given vector. Clearly, problem \eqref{eq:MPCC} fits our model \eqref{eq:nlp_gen} with $x = (u, v)$, $f(x) = f(u, v)$, $F(x) = u^{\top}(Au + Bv + c)$, and $g(x) = \delta_{\mathcal{X}}(x)$ the indicator of $\mathcal{X} = \mathcal{S} \cap \{(u,v): u \geq 0, \ Au + Bv + c \geq 0\}$.The main reference is "Facchinei, F., Jiang, H. and Qi, L.: A smoothing method for mathematical programs with equilibrium constraints. Math. Program., 85, 107–134 (1999)." We can use this problem as a synthetic example to verify the algorithm. This paper also has concrete examples for our test.}

%%%%%%%%%%%%%%%%%%%%%%%%%

%\subsection{Minimum Sum-of-Squares Clustering.}
\label{subsec:MSSC}
In this section we analyze the practical performance  of our algorithm  by comparing it with state-of-the art existing methods and software on several numerical experiments.
Our mathematical model is a semidefinite programming relaxation of large-scale clustering problem. 
By utilizing the well-known Burer-Monteiro factorization, we can reformulate such a problem into \eqref{eq:nlp_gen}. 
First, let us describe our problem as follows.

Given $m$ data points $\{a_1, \dots, a_m\} \subset \mathbb{R}^d$, we want to partition them into $k\) clusters (usually $k \ll m$). To achieve this we consider  minimizing the sum of squared Euclidean distances between the points and their cluster centroids, yielding the following  optimization problem:
\[
\min_{S_i,\mu_i} \sum_{i=1}^{k} \sum_{a_j \in S_i} \Vert  a_j - \mu_i \Vert ^2,
 \]
where $S_i \subseteq \{a_1, \cdots, a_m\}\) denotes the set of data points assigned to $i$th cluster  and $\mu_i \in \mathbb{R}^d$ is the corresponding centroid, with  $i \in \{1, \cdots,k\}\). The above minimum sum-of-squares clustering (MSSC) problem can be re-written with the help of an assignment matrix $Y = [y_{ij}] \in \mathbb{R}^{m \times k}$  defined as:
\[
y_{ij} =
\begin{cases} 
1 & \text{if } a_i \text{ is assigned to cluster } S_j \\
0 & \text{otherwise}.
\end{cases}
\]
The centroid  $\mu_j$ of cluster $S_j$ is the mean of all points in that cluster:
\[
\mu_j = \frac{\sum_{l=1}^{m} y_{lj} a_l}{\sum_{l=1}^{m} y_{lj}}.
\]
Using this fact, we can  represent  MSSC problem as \cite{PenWei:07}:
\begin{align}
    \min_{Y = [y_{ij}] \in \{0, 1\}^{m \times k}} &  \sum_{j=1}^{k} \sum_{i=1}^{m} y_{ij} \left\Vert  a_i - \frac{\sum_{l=1}^{m} y_{lj} a_l}{\sum_{l=1}^{m} y_{lj}} \right\Vert ^2  \\
   & \text{s.t.} \;\; \sum_{j=1}^{k} y_{ij} = 1 \quad \forall i = 1:m,  \label{sum_row}\\
    & \qquad  \sum_{i=1}^{m} y_{ij} \geq 1 \quad \forall j = 1:k.  \label{sum_col}
%& \qquad  y_{ij} \in \{0, 1\} \quad \forall \,  i = 1:m,\;   j = 1:k. \label{0-1}
\end{align}
The constraint \eqref{sum_row} ensures that each point $a_i$ is assigned to one and only one cluster (i.e., $Y\mathbf{1}_k=\mathbf{1}_m$, where  $\mathbf{1}_k$  denotes  the vector of all ones in $\mathbb{R}^k$) and the constraint \eqref{sum_col} ensures that each cluster is nonempty.  Consequently, we have  $\text{rank}(Y)=k$. Following  \cite{PenWei:07}, we can define:
\[
Z := [z_{ij}] = Y (Y^T Y)^{-1} Y^T \in \mathbb{R}^{m\times m}.
\]
It follows that:
\[
\text{rank}(Z)= \text{rank}(Y)=k \quad \text{ and } \quad Z\mathbf{1}_m=ZY\mathbf{1}_k=Y\mathbf{1}_k=\mathbf{1}_m.
\]
Moreover, $Z$ is a symmetric projection matrix with nonnegative entries (i.e., $Z^2=Z, Z=Z^T$ and $Z\geq 0$) and  replaces the $0-1$ constraint (i.e., $Y = [y_{ij}] \in \{0, 1\}^{m \times k}$), thus leading to the following equivalent formulation (here $A \in \mathbb{R}^{m \times d}$ is the data matrix, i.e., $i$th row is data $a_i$):
\begin{align} \label{0-1_SDP}
    \min_{Z\in \mathbb{R}^{m\times m}} & \quad \text{Tr}(A A^T (I - Z)) \\
    \text{s.t.} & \quad Z \mathbf{1}_m= \mathbf{1}_m, \;\; \text{Tr}(Z) = k,  \nonumber\\
    & \quad Z \geq 0, \; Z^2=Z, \; Z=Z^T. \nonumber
\end{align}
Since the constraints   $Z^2=Z, \; Z=Z^T,\;  Z \geq 0$  and $\text{Tr}(Z)=k$ can be equivalently replaced by the semidefinite condition $ Z \succeq 0,\;  Z \geq 0$ and   $\text{rank}(Z) = \text{Tr}(Z) = k$, we end up with the following   semidefinite programming (SDP) problem, which is  equivalent to \eqref{0-1_SDP}:
\begin{align} \label{SDP_original}
    \min_{Z\in \mathbb{R}^{m\times m}} & \quad \text{Tr}(A A^T (I - Z)) \\
    \text{s.t.} & \quad Z \mathbf{1}_m= \mathbf{1}_m, \;\; \text{Tr}(Z) = k,  \;\; \text{rank}(Z) = k, \nonumber\\
    & \quad  Z \geq 0, \; Z \succeq 0. \nonumber
\end{align}
Note that to get a convex SDP problem, we drop the rank constraint in \eqref{SDP_original}. Standard convex SDP solvers struggle to scale with the number of data points $m$, as they typically need to deal with a large number of decision variables proportional to $\mathcal{O}(m^2)$  and require projections onto the semidefinite cone, which incurs a computational cost per iteration of order $\mathcal{O}(m^3)$. To address this issue, we employ Burer-Monteiro (BM) factorization \cite{BurMon:03} for  the optimization problem \eqref{SDP_original}, trading off convexity for reduced number of variables ($\mathcal{O}(mr)$ with $k \leq r \ll m$,  usually $r$ is proportional to $k$)  and consequently less computational burden. Specifically, we aim to solve the  nonconvex~problem:
\begin{align} \label{SDP_BM}
    \min_{X \in\mathbb{R}^{m\times r}} & \quad \text{Tr}(A A^T (I - XX^T)) \\
    \text{s.t.} & \quad XX^T \mathbf{1}_m = \mathbf{1}_m, \quad \text{Tr}(X X^T) \leq r, \quad X \geq 0, \nonumber
\end{align}

\noindent In \eqref{SDP_BM}, the constraint \( \text{Tr}(Z) = k \) is relaxed to \( \text{Tr}(XX^T) \leq r \), as $k \leq r$. Moreover, the constraint \( \text{rank}(Z) = k \) is relaxed to \( \text{rank}(X) = \text{rank}(Z)  \leq r \).   It is also important to note that the constraint $Z \geq 0$ in \eqref{SDP_original} has been substituted in \eqref{SDP_BM} with the stronger, yet easier to implement, constraint $X \geq 0$. Obviously, problem \eqref{SDP_BM} is in the form of our problem  (1). Specifically,  let $x_i \in \mathbb{R}^r$ denote the $i$th row of $X$. We then construct $x \in \mathbb{R}^n$, with $n = m r$, as $x := [x_1^T, \ldots, x_m^T]^T \in \mathbb{R}^n$, and define:

\vspace{-0.5cm}

\begin{gather*}
    f(x) = \text{Tr}(A A^T)-\sum_{i,j=1}^{m} [A A^T]_{i,j}  x_i^T x_j, \quad g(x) = \delta_C(x), \quad \\
    F(x) = [x_1^T (\sum_{j=1}^{m} x_j - \mathbf{1}_r),\ldots, x_m^T (\sum_{j=1}^{m} x_j - \mathbf{1}_r)]^T\in\mathbb{R}^{m},
\end{gather*}

\vspace{-0.2cm}

\noindent where \( C \) is the  convex set defined as the intersection between the positive orthant in \( \mathbb{R}^n \) and the ball centered at \( 0_n \) with radius \( \sqrt{r} \) (since \( \text{Tr}(X X^T) = \sum_{j=1}^{m} x_j^T x_j = \|x\|^2 \)). Note that the projection onto  $C$ can be performed in closed form in $\mathcal{O}(n)$ operations.

\medskip 

\noindent We compare our algorithm, Algorithm \ref{alg:LIPAL} (LIPAL), with the Adaptive Lagrangian Minimization Scheme  (ALMS)  from \cite{HalTeb:23},  which is also an augmented Lagrangian type algorithm, on solving the nonconvex  problem \eqref{SDP_BM} after its vectorization and  another augmented Lagrangian based solver (SDPNAL+) from \cite{YanSun:15}  that  solves directly the convex SDP obtained by removing the rank constraint in  \eqref{SDP_original} and relaxing the trace constraint to \(\text{Tr}(Z) \leq r\).  We stop  LIPAL and ALMS algorithms at an \(\epsilon\)-first-order optimal solution (see Definition \ref{firstorder}), where the tolerance for stationarity condition violation is \(\epsilon_1 = 10^{-1}\) and the tolerance for constraints violation is \(\epsilon_2 = 10^{-3}\). We initialize LIPAL and ALMS with  the same point selected randomly.

\medskip 

\noindent All the codes were implemented in Matlab, and executed on a computer with (i9, CPU 3.50GHz, 64GB RAM). For our Algorithm \ref{alg:LIPAL} (LIPAL) we use the accelerated proximal gradient method  \cite{Nes:18}  for solving the  strongly convex quadratic subproblem with a simple feasible set in Step 4,  which is terminated when norm of the gradient mapping is less than the tolerance $\epsilon_{\text{sub}}= 10^{-3}$.  Note that the subproblem in the algorithm ALMS  from \cite{HalTeb:23} reduces to a projected gradient step. Finally,  SDPNAL+ computes the  projection onto positive semidefine cone using partial eigenvalue decomposition whenever it is expected to be more economical than a full eigenvalue decomposition. 

\medskip 

\noindent The numerical results are summarized in Table 1.  In our numerical tests, we consider both synthetic (first part of Table 1) and real (second part of Table 1) datasets. The synthetic data consists of $m$ randomly generated points in $\mathbb{R}^d$, with $m$ ranging from $50$ to $2000$ and $d\in \{30, 100\}$, contained in  $k=10$ separable unitary balls (clusters) (with a minimum pairwise distance of at least $3$).  The real data were taken from \cite{Par:23}, having  the number of features  \( d \) ranging from $4$ to $18$, different numbers of clusters \( k \in \{2, 3, 4\}\) and different numbers of data points \( m \) ranging from $187$ to $1372$.  For each algorithm, we report the number of iterations (number of evaluations of the first derivatives of problem's functions), CPU time (in seconds), optimal function values for problems \eqref{SDP_original} and \eqref{SDP_BM}, and feasibility violation \( \|F\| \). For SDPNAL+, we specifically consider the primal feasibility violation. If an algorithm does not solve a specific problem  in one hour we consider that problem unsolved by that algorithm and we mark the corresponding entry in Table 1 as "-".  To assess the sensitivity of our algorithm to parameters choice \( \rho \), \( \tau \), and \( \beta \), we test several values for these parameters. We chose for $\tau$ two different values and then $\rho$  is fixed so that $\frac{\tau}{\rho}$ is proportional to feasibility violation tolerance $\epsilon_2$ (see \eqref{eq:measure_critic} and \eqref{eq:choice_of_rho}). Moreover,  $\beta$ is chosen dynamically through a line search procedure to ensure the decrease in \eqref{eq:descent_of_PAL}. For ALMS, the penalty parameter is selected dynamically according to the procedure in Step 5 of the algorithm and the regularization parameter $\beta$ is  selected dynamically to ensure the descent (4.2) in \cite{HalTeb:23}.

\medskip

%%%%%%%%%%%%%
\begin{table}
\small
{
\centering   
\begin{adjustbox}{width=\columnwidth,center}
\begin{tabular}{|c|c|c|cc|cc|cc|}
    \hline
   \multirow{1}{*}{Data : $(m,d,k)$} 
   &\multicolumn{2}{c|}{\backslashbox{LIPAL params}{Algs}}
     &
     \multicolumn{2}{c|}{LIPAL} &
      \multicolumn{2}{c|}{ALMS \cite{HalTeb:23}}&
      \multicolumn{2}{c|}{SDPNAL+} \\ \cline{2-9}
       \multirow{2}{*}{Dim-SDP / Dim-BM}
        & \multirow{2}{*}{$\tau$} & \multirow{2}{*}{$\rho$} 
               
          & \# iter    & $f^{*}$  &
           \# iter     & $f^{*}$ &
           \# iter     & $f^{*}$     \\ %\cdashline{2-9}
        & & 
          
          &  cpu   & $\Vert F\Vert $ &
            cpu   & $\Vert F\Vert $  &
           cpu    & prim-feas\\
    \hline

             \multirow{2}{*}{Random: $(50, 30, 10)$}
    
      & \multirow{2}{*}{$10^{-5}$} & \multirow{2}{*}{$10$}  
    & 50 &  1.73 & 
       &  & 
      &  \\% \cdashline{2-9}
    &
     & 
      & \textbf{0.43} & 9 e-4 & 
      2770 & 1.77   & 
     369  & 1.15  \\ \cline{2-5}
\multirow{2}{*}{$1.25\times 10^3$/ $10^3$}
     & \multirow{2}{*}{$10^{-2}$} & \multirow{2}{*}{$ 10^4$} 
      &  72 &  1.73 & 
      15.85 & 9 e-4 & 
     1.35 & 8 e-5 \\% \cdashline{2-9}
    &
     & &  
      0.93 &  2 e-4 &
       &  &
       &   \\\hline

             \multirow{2}{*}{Random: $(100, 30, 10)$}

             & \multirow{2}{*}{$10^{-5}$} & \multirow{2}{*}{$10$} 
    & 124 &  4.24 & 
       &  & 
      &  \\% \cdashline{2-9}
    &
     &  
      & \textbf{3.96} & 6 e-4 & 
      1231 & 4.33   & 
     442  & 3.49  \\ \cline{2-5}
\multirow{2}{*}{$5\times 10^3$/ $2\times10^3$}
     & \multirow{2}{*}{$10^{-2}$} & \multirow{2}{*}{$ 10^4$} 
      &  79 &  4.24 & 
      29.98 & 9 e-4 & 
     4.83 & 8 e-5 \\% \cdashline{2-9}
    &
     & &  
      6.50 &  3 e-4 &
       &  &
       &   \\\hline  

                    \multirow{3}{*}{Random: $(200, 100, 10)$}

             & \multirow{2}{*}{$10^{-5}$} & \multirow{2}{*}{$10^2$} 
    & 21 &  11.06 & 
       &  & 
      &  \\% \cdashline{2-9}
    &
     &  
      & \textbf{5.38} & 3 e-4 & 
       &    & 
       &   \\ \cline{2-5}

       \multirow{3}{*}{$2\times 10^4$/ $4\times10^3$}
                    & \multirow{2}{*}{$10^{-3}$} & \multirow{2}{*}{$2\times10^3$} 
    & 102 &  11.06 & 
      119 & 11.27  & 
     446 & 10.41 \\% \cdashline{2-9}
    &
     &  
      & 33.51 & 7 e-4 & 
      14.63 &  9 e-4  & 
      21.68 & 5 e-5  \\ \cline{2-5}

     & \multirow{2}{*}{$10^{-2}$} & \multirow{2}{*}{$5\times 10^4$}  
      &  353 &  6396.06 & 
       &  & 
      &  \\% \cdashline{2-9}
    &
     & &  
      69.59 &  2 e-4 &
       &  &
       &   \\\hline  

                    \multirow{2}{*}{Random: $(500, 100, 10)$}

             & \multirow{2}{*}{$10^{-5}$} & \multirow{2}{*}{$5\times10$} 
    & 86 &  35.44 & 
       &  & 
      &  \\% \cdashline{2-9}
    &
     &  
      & \textbf{97.57} & 4 e-4 & 
      143 &35.45   & 
     443  & 31.07  \\ \cline{2-5}
\multirow{2}{*}{$1.25\times 10^5$/ $10^4$}
     & \multirow{2}{*}{$10^{-3}$} & \multirow{2}{*}{$ 2\times10^3$} 
      &  179 &  35.44 & 
      120.37 & 8 e-4 & 
     290.04 & 2 e-5 \\% \cdashline{2-9}
    &
     & &  
      264.69 &  5 e-4 &
       &  &
       &   \\\hline

                           \multirow{2}{*}{Random: $(1000, 100, 10)$}

             & \multirow{2}{*}{$10^{-5}$} & \multirow{2}{*}{$5\times 10$} 
    & 56 &  69.68 & 
       &  & 
      &  \\% \cdashline{2-9}
    &
     &  
      & \textbf{310.05} & 6 e-4 & 
      127 & 69.70   & 
     534  & 62.53  \\ \cline{2-5}
\multirow{2}{*}{$5\times 10^5$/ $10^4$}
     & \multirow{2}{*}{$10^{-3}$} & \multirow{2}{*}{$ 2\times10^4$}  
      &  154 &  69.69 & 
     465.73 & 4 e-4 & 
     3548.42 & 2 e-5 \\% \cdashline{2-9}
    &
     & &  
     1293.91 &  5 e-4 &
       &  &
       &   \\\hline

                                  \multirow{2}{*}{Random: $(2000, 100, 10)$}

             & \multirow{2}{*}{$10^{-5}$} & \multirow{2}{*}{$5\times 10$} 
    & 39 &  115.05 & 
       &  & 
      &  \\% \cdashline{2-9}
    &
     &  
      & \textbf{1069.65} & 9 e-4 & 
      112 & 115.10   & 
     -  & -  \\ \cline{2-5}
\multirow{2}{*}{$2\times 10^6$/ $4\times 10^4$}
     & \multirow{2}{*}{$10^{-3}$} & \multirow{2}{*}{$ 2\times10^4$}  
      &  - &  - & 
     1803.01 & 5 e-4 & 
     - &  - \\% \cdashline{2-9}
    &
     & &  
     - &  - &
       &  &
       &   \\\hline 
       \hline

                        \multirow{2}{*}{Wine: $(187, 13, 3)$}

             & \multirow{2}{*}{$10^{-5}$} & \multirow{2}{*}{$10$} 
    & 68 &  1062.58 & 
       &  & 
      &  \\% \cdashline{2-9}
    &
     &  
      & \textbf{2.22} & 4 e-4 & 
      851 & 1064.38 & 
     225  & 985.3 \\ \cline{2-5}
\multirow{2}{*}{$46056$/ $1212$}
     & \multirow{2}{*}{$10^{-2}$} & \multirow{2}{*}{$ 2\times10^3$}  
      &  238 &  1058.57 & 
     12.17 & 8 e-4 & 
     15.21 &  9 e-5 \\% \cdashline{2-9}
    &
     & &  
     23.11 &  5 e-4 &
       &  &
       &   \\\hline

        \multirow{2}{*}{Breast: $(277, 9, 2)$}

             & \multirow{2}{*}{$10^{-5}$} & \multirow{2}{*}{$10$} 
    & 122 &  1562.7 & 
       &  & 
      &  \\% \cdashline{2-9}
    &
     &  
      & \textbf{6.4} & 7 e-4 & 
      3039 & 1563.01 & 
     213  &1504.2  \\ \cline{2-5}
\multirow{2}{*}{$38503$/ $1108$}
     & \multirow{2}{*}{$10^{-2}$} & \multirow{2}{*}{$ 2\times10^3$}  
      &  562 &  1563.21 & 
     81.99 & 8 e-4 & 
     37.93 &  5 e-5 \\% \cdashline{2-9}
    &
     & &  
     48.51 &  4 e-4 &
       &  &
       &   \\\hline 

                                                    \multirow{2}{*}{Heart: $(303, 13, 2)$}

             & \multirow{2}{*}{$10^{-5}$} & \multirow{2}{*}{$10$}  
    & 112 &  2843.06 & 
       &  & 
      &  \\% \cdashline{2-9}
    &
     &  
      & \textbf{8.95} & 3 e-4 & 
      2871 & 2843.12 & 
     253  & 2719.0 \\ \cline{2-5}
\multirow{2}{*}{$46056$/ $1212$}
     & \multirow{2}{*}{$10^{-2}$} & \multirow{2}{*}{$ 2\times10^3$}  
      &  730 &  2840.92 & 
     85.91 & 5 e-4 & 
     48.43 &  7 e-5 \\% \cdashline{2-9}
    &
     & &  
     124.7 &  4 e-4 &
       &  &
       &   \\\hline

                                         \multirow{2}{*}{Vehicle: $(818, 18, 4)$}

             & \multirow{2}{*}{$10^{-5}$} & \multirow{2}{*}{$10$} 
    & 285 &  4554.7 & 
       &  & 
      &  \\% \cdashline{2-9}
    &
     &  
      & \textbf{249.62} & 6 e-4 & 
      2145 & 4555.3 & 
     668  & 4174.1  \\ \cline{2-5}
\multirow{2}{*}{$334562$/ $6544$}
     & \multirow{2}{*}{$10^{-3}$} & \multirow{2}{*}{$ 10^4$}  
      &  869 &  4558.12 & 
     739.34 & 7 e-4 & 
     1946.24 &  8 e-6 \\% \cdashline{2-9}
    &
     & &  
     2134.63 &  6 e-4 &
       &  &
       &   \\\hline

                                                \multirow{2}{*}{Banknote: $(1372, 4, 2)$}

             & \multirow{2}{*}{$10^{-5}$} & \multirow{2}{*}{$10$}  
    & 36 &  1984.22 & 
       &  & 
      &  \\% \cdashline{2-9}
    &
     &  
      & \textbf{68.86} & 4 e-4 & 
      498 & 1984.34.3 & 
     -  & - \\ \cline{2-5}
\multirow{2}{*}{$941192$/ $5488$}
     & \multirow{2}{*}{$10^{-2}$} & \multirow{2}{*}{$ 2\times10^3$} 
      &  544 &  1984.07 & 
     273.49 & 8 e-4 & 
     - &  - \\% \cdashline{2-9}
    &
     & &  
     1410.73 &  5 e-4 &
       &  &
       &   \\\hline

          \multirow{2}{*}{Spambase: $(2000, 57, 2)$ }
    
      & \multirow{2}{*}{$10^{-5}$} & \multirow{2}{*}{$10^2$} 
    & 221 &  97304.2 & 
       &  & 
      &  \\% \cdashline{2-9}
    &
     &  
      & \textbf{1867.9} & 4 e-4 & 
      - & -   & 
     -  & -  \\ \cline{2-5}
\multirow{2}{*}{$2\times 10^6$/ $4\times 10^3$}
     & \multirow{2}{*}{$10^{-2}$} & \multirow{2}{*}{$2\times 10^2$}  
      &  - &  - & 
      - & - & 
     - & - \\% \cdashline{2-9}
    &
     & &  
      - &  - &
       &  &
       &   \\\hline

\end{tabular}%
\end{adjustbox}
}
\caption{Performance comparison between LIPAL, ALMS and SDPNAL+ on clustering using  synthetic (top)  and real (bottom) datasets with different sizes. For each algorithm we provide the number of iterations, cpu time in sec, optimal value function and feasibility violation. The best  time achieved by an algorithm to solve a given clustering problem is written in bold.}
\label{tab1}
\end{table}

\medskip

\noindent From Table 1, we observe that a smaller perturbation parameter  \( \tau \) leads to faster convergence for LIPAL. We attribute this to the fact that smaller values of $\tau$ makes our method to resemble more to a standard augmented Lagrangian scheme rather than a quadratic penalty scheme.  Moreover, for small values of the perturbation parameter \( \tau \), our algorithm outperforms ALMS and SDPNAL+ in terms of computational time (sometimes even 10 times faster). For SDPNAL+, the slower performance can be attributed to the fact that the problem dimension it solves is always larger than that handled by LIPAL and ALMS (see the first column of Table 1). In contrast, the slower performance of ALMS is primarily due to its high number of iterations, which we believe it is due to the  choice of the parameters and  the large approximation error of the model generated by the linearization of the smooth part of the augmented Lagrangian function compared to a linearization in a Gauss-Newton type setting used in LIPAL. %We believe this is caused by the sensitivity of its parameters, while ALMS consistently finds approximately feasible points, it requires significantly more time to satisfy approximate stationarity conditions.  
Additionally, we observe that the objective values obtained by SDPNAL+ are better than those produced by LIPAL and ALMS. This can be explained by the fact that the feasible set of the problem solved by SDPNAL+ is much larger than the one given in \eqref{SDP_BM}. Indeed, in the convex SDP, we have the nonnegativity constraint on the entries of \( Z \) (i.e., \( Z \geq 0 \)), while in \eqref{SDP_BM} a stronger constraint is required (i.e., \( X \geq 0 \));  additionally, removing $\text{rank}(Z) =k$ in \eqref{SDP_original} and considering only $\text{Tr}(Z) \leq r$ may yield solutions $Z^*$ with rank larger than $r$, while \eqref{SDP_BM} forces rank of the matrix $Z$ to be always less than $r$.

\medskip

\noindent Finally, we evaluate the clustering performance of the compared methods on a synthetic dataset consisting of \( m = 150 \) data points, generated from \( 10 \) randomly created (possibly) overlapping unitary balls in a \( d = 2 \) dimensional space. We set \( r = 12 \) and the results are presented in Figure 1, which displays the data distribution and the sparsity pattern of the matrix \( Z \) yielded by the three algorithms.  Since the data points are ordered in the data matrix,  $A$ i.e., data points belonging to the same cluster are grouped together in \( A \)), the matrix $Z$ should be block diagonal. We observe that all the three algorithms successfully identify  \( 9 \) clusters. Since the balls 7 and 8 are overlapping, ALMS attributes the data points in these two balls to only one cluster,  while LIPAL and SDPNAL+ although detect some  correlation between the two balls,  they still cluster them.  We attribute this to the fact that the dataset is not fully clusterable. In the  experiments from Table 1 on random data where unitary balls were non-overlapping, we observed that the clustering was successful. This suggests that the presence of overlapping regions between clusters may affect the quality of the clustering results for the two relaxations \eqref{SDP_original} (removing the rank constraints) and \eqref{SDP_BM}. Additionally, the choice of \( r \) may also influence the quality of the solution.  

\vspace{-0.4cm}

\begin{figure}[htp]
   \begin{center}
       \includegraphics[width=\textwidth,height=5cm]{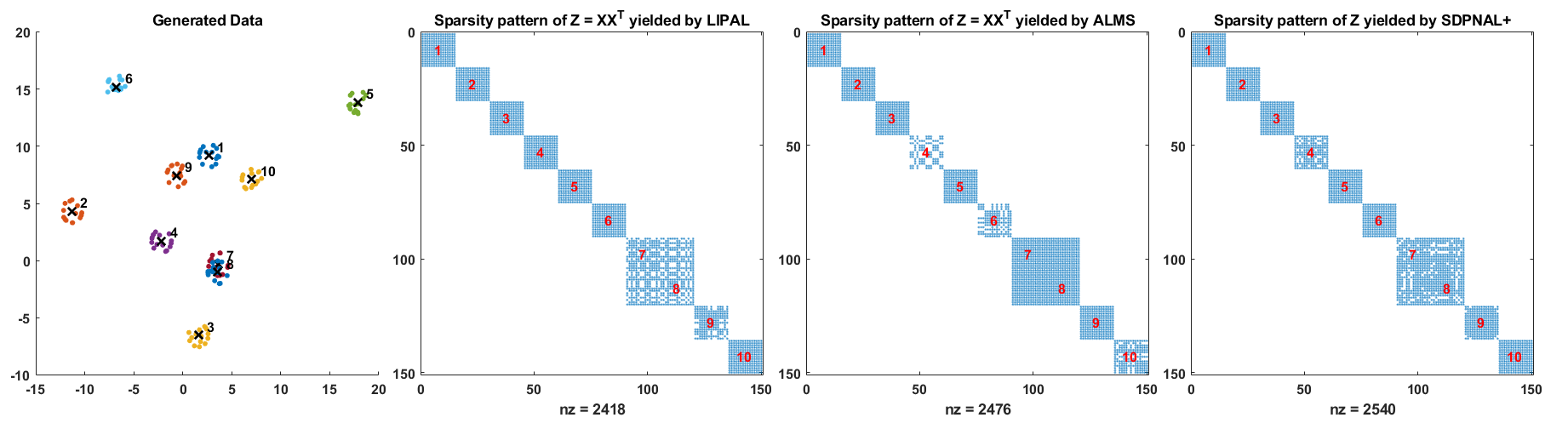}
       \vspace*{-0.7cm}
    \caption{Clustering of $m=150$ random generated data of dimension $d=2$ using LIPAL, ALMS and SDPNAL+.}
    \label{fig2}   
   \end{center}
\end{figure}

\vspace{-0.5cm}

\section{Conclusions}\label{sec6}
In this paper we have studied a class of general  composite optimization problems with nonlinear equality constraints and possibly nonsmooth and nonconvex objective function.  
We have proposed a linearized perturbed augmented Lagrangian method to solve this problem class, where we have linearized the smooth part of a perturbed augmented Lagrangian function in a Gauss-Newton fashion and added a quadratic regularization.  We have also  introduced a new constraint qualification condition that allows us to bound the dual iterates. Consequently,  we have derived global sublinear convergence rates for the iterates of our method and  improved local  convergence results were obtained under the KL condition.  Finally, the numerical experiments have demonstrated the effectiveness of our proposed algorithm in solving large-scale clustering problems.

%%%%%%%%%%%%%%%%%

\section*{Acknowledgments.}
\noindent 
The research of Necoara and El Bourkhissi was  supported by the European Union's Horizon 2020 research and innovation programme under the Marie Sklodowska-Curie grant agreement No. 953348. The research of Tran-Dinh was partly supported by the National Science Foundation (NSF), grant no. NSF-RTG DMS-2134107 and the Office of Naval Research (ONR), grant No. N00014-23-1-2588 (2023-2026).

%%%%%%%%%%%%%%%%%%%%%%%%%%%%%%%%%%%%%%%

\section*{Appendix}
This appendix provides the full proof of the technical results in the main text.

\medskip 

%%% Proof of Lemma 1.
\proof{\textbf{Proof of Lemma \ref{lemma2}.}} 
Denote $y_{\tau} := \tau y^0 + (1-\tau)y$.
By the definition of $\psi^{\tau}_{\rho}$ given in \eqref{eq:pert_AL} and of $\nabla_x{\psi^{\tau}_{\rho}}$ given in \eqref{eq:grad_psi}, we  have 
\begin{align*}
 \Vert \nabla_x & \psi^{\tau}_{\rho}(x,y;y^0)  -  \nabla_x \psi^{\tau}_{\rho}(x',y;y^0) \Vert \\
 & =   \big\Vert  \nabla f(x) + {J_F(x)}^T ( y_{\tau} + \rho F(x)) - \nabla f(x') - {J_F(x')}^T ( y_{\tau}  + \rho F(x') ) \big\Vert   \\
  &=  \big\Vert  ( \nabla f(x)-\nabla f(x')) + ({{J_F}(x)}-{{J_F}(x')} )^T(y_{\tau} + \rho F(x) ) + \rho {{J_F}(x')}^T ( F(x) - F(x') )\big\Vert   \\
  &\leq  \left\Vert \nabla f(x)-\nabla f(x')\right\Vert +\left\Vert {{J_F}(x)}-{{J_F}(x')}\right\Vert \left\Vert  y_{\tau} +\rho F(x)\right\Vert +\rho \left\Vert J_F(x')\right\Vert \left\Vert F(x)-F(x')\right\Vert  \\
  &{\overset{\text{Ass. } \ref{as:Assum}}{\leq}}  \left(L_f+L_F\left\Vert  y_{\tau} + \rho F(x)\right\Vert +\rho M_F^2\right)\left\Vert x-x'\right\Vert  \leq   L^{\tau}_{\rho}\left\Vert x-x'\right\Vert ,
\end{align*}
which proves our statement with $L^{\tau}_{\rho} \triangleq \sup_{(x,y)\in\mathcal{S}\times\mathcal{Y}}\left\{L_f + L_F\Vert  y_{\tau} + \rho F(x)\Vert +\rho M_F^2\right\}$.  
\Eproof
\endproof   
%%% End of Proof.

\medskip

%%% Proof of Lemma 2.
\proof{\textbf{Proof of Lemma \ref{le:lambda_bou}.}} 
Recall that $y^k_{\tau} := \tau y^0 + (1 - \tau)y^k$.
Using the dual update $y^{k+1} = y^k_{\tau} + \rho F(x^{k+1})$, we have
\begin{equation*}
\arraycolsep=0.2em
\begin{array}{lcl}
y^{k+1} &= &  \tau y^0 + (1 - \tau)y^k +  \rho F(x^{k+1}) \vspace{1ex} \\
& = & (1-\tau)^{k+1}y^0 + \tau\sum_{i=0}^{k}{(1-\tau)^{i}y^0}+\sum_{i=1}^{k+1}{(1-\tau)^{k+1-i}\rho F(x^{i})} \vspace{1ex} \\
& = & y^0+\sum_{i=1}^{k+1}{(1-\tau)^{k+1-i}\rho F(x^{i})}.
\end{array}
\end{equation*}
Therefore, by the triangle inequality, we can show that 
\begin{equation*}
\begin{array}{lcl}
\Vert y^{k+1}\Vert &\leq\Vert y^0\Vert +\rho\sum_{i=1}^{k+1}{(1-\tau)^{k+1-i}\Vert F(x^{i})\Vert } \leq\Vert y^0\Vert +\frac{\rho}{\tau}\Delta. 
\end{array}
\end{equation*}
Similarly, we can also get
\begin{equation*}
\begin{array}{lcl}
 \Vert y^{k+1}-y^0\Vert &\leq\rho\sum_{i=1}^{k+1}{(1-\tau)^{k+1-i}\Vert F(x^{i})\Vert } \leq\frac{\rho}{\tau}\Delta.
\end{array}
\end{equation*}
These statements prove \eqref{eq:bound_for_dual}.
\Eproof
\endproof 
%%% End of Proof.

\medskip

%%% Proof of Lemma 3.
\proof{\textbf{Proof of Lemma \ref{le:lambda_bou1}}.} 
Using the optimality condition of \eqref{eq:LIPAL_subprob} as $0 \in \partial{\mathcal{Q}_k}(x^{k+1})$ and the dual update $y^{k+1} = y^k_{\tau} + \rho F(x^{k+1})$, we can show that
\begin{equation*}
  \left(-\nabla f(x^{k})-J_F(x^{k})^Ty^{k+1}+\rho J_F(x^k)^T\left(F(x^{k+1}) - \ell_F(x^{k+1};x^k)\right)-\beta\Delta x^{k},-1\right)\in N^{\text{lim}}_{\epi g}(x^{k+1},g(x^{k+1})).
\end{equation*}
Since $N^{\text{lim}}_{\epi g}(x^{k+1},g(x^{k+1}))$ is a cone, it follows that
\begin{equation*}
	\frac{\tau}{\rho} \! \left(-\nabla f(x^{k})-J_F(x^{k})^Ty^{k+1} \!+\rho J_F(x^k)^T\left(F(x^{k+1}) - \ell_F(x^{k+1};x^k)\right) \!-\! \beta\Delta x^{k},-1\right) \!\in\!  N^{\text{lim}}_{\epi g}(x^{k+1},g(x^{k+1})).
\end{equation*}
By the triangle inequality, one can show that
\begin{equation*}
\arraycolsep=0.2em
\begin{array}{lcl}
	\Tc_{[1]} & := & \mathrm{dist}\big( -\frac{\tau}{\rho}J_F(x^{k})^Ty^{k+1},\partial^{\infty}g(x^{k+1}) \big) \vspace{1ex}\\
	& = & \mathrm{dist}\big( \big(-\frac{\tau}{\rho}J_F(x^{k})^Ty^{k+1},0 \big), N^{\text{lim}}_{\epi g}(x^{k+1},g(x^{k+1})) \big) \vspace{1ex} \\
	& \leq & \frac{\tau}{\rho}\left(1+\Vert \nabla f(x^{k})\Vert +\rho\Vert J_F(x^k)\Vert \Vert F(x^{k+1})-l_F(x^{k+1};x^k)\Vert +\beta\Vert \Delta x^{k}\Vert \right) \vspace{1ex} \\
	& \leq & \frac{\tau}{\rho} \big[ M_f+1+\left(\beta+2\rho M_F^2\right)\Vert \Delta x^{k}\Vert \big]  \vspace{1ex}\\
	& \leq & \frac{\tau}{\rho} \big[ M_f+1+\left(\beta+2L^{\tau}_{\rho}\right)\Vert \Delta x^{k}\Vert \big] \vspace{1ex}\\
	& \leq & \frac{\tau}{\rho} \big( M_f+1+2\beta\Vert \Delta x^{k}\Vert \big).
\end{array} 
\end{equation*}
Again, by  the triangle inequality, we also have
\begin{equation*}
\arraycolsep=0.2em
\begin{array}{lcl}
	\mathrm{dist} \big(-\frac{\tau}{\rho}J_F(x^{k})^T\left(y^{k+1}-y^0\right),\partial^{\infty}g(x^{k+1}) \big) & \leq & \mathrm{dist} \big(-\frac{\tau}{\rho}J_F(x^{k})^Ty^{k+1},\partial^{\infty}g(x^{k+1}) \big) + \frac{\tau}{\rho}\Vert y^0\Vert \vspace{1ex} \\
	& \leq &\frac{\tau}{\rho} \big[  M_f+1+\Vert y^0\Vert +\left(\beta+2\rho M_F^2\right)\Vert \Delta x^{k}\Vert \big]. 
\end{array} 
\end{equation*}
Therefore, we can show that
\begin{equation*}
\arraycolsep=0.2em
\begin{array}{lcl}
	\mathrm{dist} \big( -\frac{\tau}{\rho}J_F(x^{k+1})^T\left(y^{k+1}-y^0\right),\partial^{\infty}g(x^{k+1}) \big) & \leq & \frac{\tau}{\rho} \big[ M_f+1+\Vert y^0\Vert +\left(\beta+2\rho M_F^2\right)\Vert \Delta x^{k}\Vert \big] \vspace{1ex}\\
	&& + {~} \frac{\tau}{\rho}\Vert J_F(x^{k+1})-J_F(x^{k})\Vert \Vert y^{k+1}-y^0\Vert \vspace{1ex} \\
	& \leq & \frac{\tau}{\rho}\big[ M_f+1+\Vert y^0\Vert +\big( \beta + 2\rho M_F^2 + \frac{\rho \Delta}{\tau}L_F \big) \Vert \Delta x^{k}\Vert \big] \vspace{1ex}\\
	& \leq & \frac{\tau}{\rho}\left( M_f+1+\Vert y^0\Vert +2\beta\Vert \Delta x^{k}\Vert \right).
\end{array} 
\end{equation*}
Furthermore, on the one hand, we can also prove that
\begin{equation}\label{feasibility_measure}
\arraycolsep=0.2em
\begin{array}{lcl}
	\mathrm{dist} \big(-J_F(x^{k+1})^TF(x^{k+1}),\partial^{\infty}g(x^{k+1}) \big) & = & \mathrm{dist}\big(-J_F(x^{k+1})^T\frac{(1-\tau)\Delta y^k+\tau\left(y^{k+1}- y^0\right)}{\rho},\partial^{\infty}g(x^{k+1}) \big) \vspace{1ex}\\
	& \leq & \mathrm{dist} \big( -\frac{\tau} {\rho}J_F(x^{k+1})^T\left(y^{k+1}-y^0\right),\partial^{\infty}g(x^{k+1}) \big) \vspace{1ex}\\
	&& + {~} \frac{1-\tau}{\rho}\Vert J_F(x^{k+1})^T\Delta y^k\Vert \vspace{1ex}\\
	&  \leq & \frac{\tau}{\rho}\big[ (M_f+1+\Vert y^0\Vert ) + 2\beta\Vert \Delta x^{k}\Vert \big] + \frac{(1-\tau)M_F}{\rho}\Vert \Delta y^k\Vert.
\end{array} 
\end{equation}
On the other hand, by the triangle inequality, one can easily show that
\begin{equation*}
\arraycolsep=0.2em
\begin{array}{lcl}
	\Vert F(x^{k+1})\Vert & = & \frac{1}{\rho}\Vert (1-\tau)\Delta y^k+\tau\left( y^{k+1}-y^0\right)\Vert \vspace{1ex}\\
	& \geq & \frac{\tau}{\rho}\Vert y^{k+1}-y^0\Vert  -\frac{1-\tau}{\rho}\Vert \Delta y^k\Vert.
\end{array} 
\end{equation*}
Hence, by using  Assumption \ref{as:Assum}-\textrm{(iii)}, we eventually get
\begin{equation*}
\arraycolsep=0.2em
\begin{array}{lcl}
\frac{\sigma\tau}{\rho}\Vert y^{k+1}-y^0\Vert  & \leq &\frac{\tau}{\rho} \big[ (M_f+1+\Vert y^0\Vert )+2\beta\Vert \Delta x^{k}\Vert \big] + \frac{(1-\tau)}{\rho}\left(M_F+\sigma\right)\Vert \Delta y^k\Vert,
\end{array} 
\end{equation*}
which proves \eqref{eq:lambda_1}.
\Eproof
\endproof 
%%%% End of Proof.

\medskip

%%% Proof of Lemma 4.
\proof{\textbf{Proof of Lemma \ref{le:lambda_bound}.}} 
Using the dual update twice at iterations $k$ and $k-1$,  we have
\begin{equation*} 
\begin{array}{lclcl}
(1-\tau)\Delta{y}^k &= & (1-\tau)\left(y^{k+1}-y^k\right) & = & \rho F(x^{k+1})-\tau \left(y^{k+1}-y^0\right), \vspace{1ex}\\
(1-\tau)\Delta{y}^{k-1} & = & (1-\tau)\left(y^{k}-y^{k-1}\right) & = & \rho F(x^{k})-\tau \left(y^{k}-y^0\right).
\end{array}
\end{equation*}
Computing the difference between two lines, and taking the inner product with $\Delta{y}^k$, we get
\begin{equation}\label{eq:dual_proof1} 
\begin{array}{lclcl}
(1-\tau)\iprods{\Delta y^k-\Delta y^{k-1},\Delta y^k} = \rho \iprods{ F(x^{k+1})-F(x^k), \Delta y^k} - \tau\norms{\Delta y^k}^2.
\end{array}
\end{equation}
On the one hand, we expand and lower bound the left-hand side of \eqref{eq:dual_proof1} as
\begin{equation*} 
\begin{array}{lclcl}
\iprods{\Delta y^k-\Delta y^{k-1},\Delta y^k} = \frac{1}{2}\left(\norms{\Delta y^k-\Delta y^{k-1}}^2 + \norms{\Delta y^k}^2-\norms{\Delta y^{k-1}}^2\right)\geq\frac{1}{2}\left(\norms{\Delta y^k}^2-\norms{\Delta y^{k-1}}^2\right).
\end{array}
\end{equation*}
On the other hand, we upper bound $\rho \iprods{ F(x^{k+1})-F(x^k), \Delta y^k}$ on the right-hand side of \eqref{eq:dual_proof1} as
\begin{equation*} 
\begin{array}{lclcl}
 \rho \iprods{ F(x^{k+1})-F(x^k), \Delta y^k} &\leq & \frac{\rho^2}{2\tau}\Vert F(x^{k+1})-F(x^k)\Vert ^2+\frac{\tau}{2}\Vert \Delta y^k\Vert ^2  \vspace{1ex} \\
 & \leq &  \frac{\rho^2M_F^2}{2\tau}\norms{\Delta x^k}^2+\frac{\tau}{2}\Vert \Delta y^k\Vert ^2. 
\end{array}
\end{equation*}
Substituting the last two inequalities into \eqref{eq:dual_proof1}, and rearranging the result, we obtain \eqref{eq:lambda_squared}.
\Eproof
\endproof 
%%% End of Proof.

\medskip

%%% Proof of Lemma 5.
\proof{\textbf{Proof of Lemma \ref{lemma3}.}}
Recall that $y^k_{\tau} := \tau y^0 + (1-\tau)y^k$.
Then, since $x^{k+1}$ solves \eqref{eq:LIPAL_subprob}, we have
\begin{equation*} 
\begin{array}{lclcl}
	\mathcal{Q}_k(x^{k+1}) &= &  \ell_f(x^{k+1};x^{k})+g(x^{k+1}) +\langle y^k_{\tau}, \ell_F(x^{k+1};x^{k})\rangle+\frac{\rho}{2}{\Vert  \ell_F(x^{k+1};x^{k})\Vert ^2} \vspace{1ex} \\
	& \leq & \mathcal{Q}_k(x^k) - \frac{\beta}{2}\norms{x^{k+1} - x^k}^2 \vspace{1ex}\\
	& = &  f(x^k)+g(x^{k})+\langle(1-\tau) y^{k}+\tau y^0 ,F(x^{k}) \rangle+\frac{\rho}{2}{\Vert F(x^{k})\Vert ^2}-\frac{\beta}{2}\Vert \Delta x^{k}\Vert^2.
\end{array}
\end{equation*}
Rearranging this inequality yields
\begin{equation}\label{eq:use_bellow} 
\begin{array}{lclcl}
	g(x^{k+1}) - g(x^{k}) & \leq & -\langle \nabla f(x^{k}) ,\Delta x^{k} \rangle -\langle  J_F(x^{k})\Delta x^{k}, y^k_{\tau} \rangle-\frac{\rho}{2}\langle J_F(x^{k})\Delta x^{k} ,2F(x^{k})\rangle \vspace{1ex}\\
	&& {~}  -\frac{\rho}{2}\langle J_F(x^{k})\Delta x^{k} ,J_F(x^{k})\Delta x^{k} \rangle-\frac{\beta}{2}\Vert \Delta x^{k}\Vert ^2 \vspace{1ex}\\
	& = & -\langle { \nabla f(x^{k})+J_F(x^{k})}^Ty^k_{\tau} + \rho F(x^{k}) ,\Delta x^{k} \rangle-\frac{\rho}{2}\Vert J_F(x^{k})\Delta x^{k}\Vert ^2-\frac{\beta}{2}\Vert \Delta x^{k}\Vert ^2 \vspace{1ex}\\
	& \leq & -\langle \nabla_x\psi^{\tau}_{\rho}(x^{k},y^{k};y^0) ,\Delta x^{k} \rangle -\frac{\beta}{2}\Vert \Delta x^{k}\Vert ^2. 
\end{array}
\end{equation}
Using the $L_{\rho}^{\tau}$-smoothness of $\psi_{\rho}$ from \eqref{eq:pert_AL}, we get
\begin{equation*}
\psi_{\rho}(x^{k+1},y^k;y^0)- \psi_{\rho}(x^{k},y^k;y^0)- \langle \nabla_x\psi^{\tau}_{\rho}(x^{k},y^{k};y^0) ,\Delta x^{k} \rangle\leq \frac{L^{\tau}_{\rho}}{2}\Vert \Delta x^k\Vert ^2,
\end{equation*}
Adding $f(x^{k+1}) - f(x^k)$ to the last inequality, and using $\mathcal{L}^{\tau}_{\rho}$ from  \eqref{eq:pert_AL} and \eqref{eq:use_bellow}, we can show that
\begin{equation*} 
\begin{array}{lclcl}
	\mathcal{L}^{\tau}_{\rho}(x^{k+1},y^{k};y^0) - \mathcal{L}^{\tau}_{\rho}(x^{k},y^{k};y^0) & = & g(x^{k+1})-g(x^{k})+\psi^{\tau}_{\rho}(x^{k+1},y^{k};y^0)-\psi^{\tau}_{\rho}(x^{k},y^{k};y^0) \vspace{1ex} \\
  & {\overset{{\eqref{eq:use_bellow}}}{\leq}} & -\frac{\beta}{4}\Vert \Delta x^{k}\Vert ^2.   
\end{array}
\end{equation*}
This proves the first inequality of \eqref{eq:descent_of_PAL}.
The second one follows from the condition $\beta\geq 2L^{\tau}_{\rho}$.
\Eproof
\endproof 
%%%% End of Proof.

\medskip

%%%% Proof of Lemma 6.
\proof{\textbf{Proof of Lemma \ref{le:measure_critic}.}}
From the optimality condition of \eqref{eq:LIPAL_subprob}, we have
\begin{equation*}
	\mathrm{dist}\left(-\nabla f(x^{k})-J_F(x^{k})^Ty^{k+1}+\rho J_F(x^{k})^T\left(F(x^{k+1})-\ell_F(x^{k+1};x^k)\right)-\beta\Delta x^{k},\partial g(x^{k+1})\right)=0.
\end{equation*}
Applying the triangle inequality to this expression, and using the $M_F$-boundedness of $J_F$, we get
\begin{equation*}
\arraycolsep=0.2em
\begin{array}{lcl}
	\mathrm{dist}\big(-\nabla f(x^{k})-{{J_F}(x^{k})}^Ty^{k+1},\partial g(x^{k+1}) \big)   & \leq & \rho\Vert J_F(x^{k})\Vert \Vert F(x^{k+1})-\ell_F(x^{k+1};x^k)\Vert +\beta\Vert \Delta x^{k}\Vert \vspace{1ex} \\
	& \leq & \left(\beta+2\rho M_F^2\right)\Vert \Delta x^{k}\Vert. 
\end{array}    
\end{equation*}
Again, by the triangle inequality, the $L_f$-smoothness of $f$, the $L_F$-smoothness of $F$, and $\Vert y^{k+1}\Vert \leq\Vert y^0\Vert +\frac{\rho\Delta}{\tau}$ from Lemma \ref{le:lambda_bou}, we can also prove that
\begin{equation*}
\arraycolsep=0.1em
\begin{array}{lcl}
\Tc_{[2]}  &:= & \mathrm{dist} \big( -\nabla f(x^{k+1})-{{J_F}(x^{k+1})}^Ty^{k+1},\partial g(x^{k+1}) \big) \vspace{1ex} \\
	& = &  \mathrm{dist}\big(-\nabla f(x^{k})-J_F(x^{k})^Ty^{k\!+\!1}+\big(\nabla f(x^{k}) - \nabla f(x^{k \! + \! 1})+\left(J_F(x^{k})-J_F(x^{k\!+\!1})\right)^Ty^{k\!+\!1} \big),\partial g(x^{k\!+\!1}) \big) \vspace{1ex} \\
	& \leq & \mathrm{dist}\big(-\nabla f(x^{k})-{{J_F}(x^{k})}^Ty^{k \!+ \! 1},\partial g(x^{k\!+\!1}) \big) + \Vert \nabla f(x^{k \!+\! 1})-\nabla f(x^{k})\Vert +\Vert {{J_F}(x^{k \!+\! 1})}-{{J_F}(x^{k})}\Vert \Vert y^{k \!+\! 1}\Vert \vspace{1ex}\\
	& \leq & \left(\beta+L_f+\Vert y^0\Vert L_F+\frac{\rho\Delta}{\tau}L_F+2\rho M_F^2\right)\Vert \Delta x^{k}\Vert \vspace{1ex}\\
	& \leq & 2\beta\Vert \Delta x^{k}\Vert,
\end{array}    
\end{equation*}
which proves the first line of \eqref{eq:measure_critic}.  Finally, using the update of the dual variable as $y^{k+1} = \tau y^0 + (1-\tau)y^k + \rho F(x^{k+1})$, we have
\begin{equation*}
\arraycolsep=0.2em
\begin{array}{lcl}
	\Vert F(x^{k+1})\Vert & = & \frac{1}{\rho}\left\Vert (1-\tau)\Delta y^k+\tau \left(y^{k+1}-y^0\right)\right\Vert \vspace{1ex} \\
	& \leq & \frac{1-\tau}{\rho}\Vert \Delta y^{k}\Vert +\frac{\tau}{\rho}\Vert y^{k+1}-y^0\Vert  \vspace{1ex} \\
	& \leq & \frac{2\beta\tau}{\rho\sigma}\Vert \Delta x^{k}\Vert +\frac{1-\tau}{\rho}\left(2+\frac{M_F}{\sigma}\right)\Vert \Delta y^{k}\Vert +\frac{\tau }{\rho\sigma}\left(M_f+1+\Vert y^0\Vert \right),
\end{array}    
\end{equation*}
which proves the second line of \eqref{eq:measure_critic},  where the first inequality follows from the triangle inequality and the last one follows from Lemma \ref{le:lambda_bou1}.
\Eproof
\endproof 
%%% End of Proof

\medskip

%%% Proof of Lemma 7.
\proof{\textbf{Proof of Lemma \ref{le:decrease_of_P}.}}
Combining \eqref{eq:Lyapunov_func}, \eqref{eq:lyapunov_function}, \eqref{eq:pert_AL}, \eqref{eq:descent_of_PAL}, and \eqref{eq:lambda_squared}, we can derive that
\begin{equation*} 
\arraycolsep=0.2em
\begin{array}{lcl}
\Pc_{k+1} & = & \Lc^{\tau}_{\rho}(x^{k+1}, y^{k+1};y^0) - \frac{\tau(1-\tau)}{2\rho}\norms{y^{k+1}-y^0}^2  + \frac{2(1-\tau)^2}{\tau\rho}\norms{\Delta{y}^k}^2 \vspace{1ex} \\
& \overset{\tiny\eqref{eq:pert_AL}}{=} & \Lc^{\tau}_{\rho}(x^{k+1}, y^k; y^0) + (1-\tau)\iprods{\Delta{y}^k, F(x^{k+1})} - \frac{\tau(1-\tau)}{2\rho}\norms{y^{k+1}-y^0}^2  + \frac{2(1-\tau)^2}{\tau\rho}\norms{\Delta{y}^k}^2 \vspace{1ex} \\
& \overset{\tiny\eqref{eq:descent_of_PAL}}{\leq}  &   \Lc^{\tau}_{\rho}(x^{k}, y^k;y^0) - \frac{\tau(1-\tau)}{2\rho} \norms{y^k-y^0}^2 + \frac{2(1-\tau)^2}{\tau\rho}\norms{\Delta{y}^{k-1}}^2 \vspace{1ex} \\
&& + {~} (1-\tau)\iprods{F(x^{k+1}),\Delta y^k} - \frac{\beta}{4}\norms{\Delta x^k}^2 + \frac{2(1-\tau)^2}{\tau\rho}\left(\norms{\Delta y^k}^2-\norms{\Delta y^{k-1}}^2\right) \vspace{1ex} \\
&&  - {~} \frac{\tau(1-\tau)}{2\rho} \left(\norms{y^{k+1}-y^0}^2-\norms{y^{k}-y^0}^2\right)\\
& \overset{\tiny\eqref{eq:lambda_squared}}{\leq}  & \Pc_k  - \big[ \frac{\beta}{4}-\frac{(1-\tau)\rho M_F^2}{\tau^2} \big] \norms{\Delta x^k}^2-\frac{\tau(1-\tau)}{\rho}\iprods{y^{k}-y^0,y^{k}-y^0-\left(y^{k+1}-y^0\right)}\vspace{1ex} \\
&&  - {~} \frac{\tau(1-\tau)}{2\rho} \left[ \norms{y^{k+1}-y^0}^2-\norms{y^{k}-y^0}^2\right] \vspace{1ex}\\
& \leq  & \Pc_k -  \big[ \frac{\beta}{4} - \frac{(1-\tau)\rho M_F^2}{\tau^2} \big] \norms{\Delta x^k}^2-\frac{\tau(1-\tau)}{2\rho}\Vert \Delta y^k\Vert ^2.
\end{array}
\end{equation*}
Using the choice of $\beta$ from \eqref{eq:gamma_rho}, we obtain \eqref{eq:decrease_Lyapunov} from this estimate.
\Eproof
\endproof 
%%%% End of Proof.

\medskip

%%% Proof of Lemma 8.
\proof{\textbf{Proof of Lemma \ref{bbound}.}}
First, let us prove $\Pc_k \leq \bar{\Pc}$ for all $k \geq 1$ by induction. 
Indeed, since $x^0, x^1 \in\mathcal{S}$, from \eqref{eq:descent_of_PAL} of Lemma \ref{lemma3} for $k=0$, it follows that 
\begin{equation}\label{eq:lm8_ine111}
\arraycolsep=0.2em
\begin{array}{lcl} 
	\mathcal{L}_{\tau}^{\rho}(x^1, y^0; y^0) + \frac{\beta}{4}\Vert x^1 - x^0\Vert^2 & = & \phi(x^1) + \langle y^0  ,F(x^1)\rangle+\frac{\rho}{2}\Vert F(x^1)\Vert ^2+\frac{\beta}{4}\Vert x^1-x^0\Vert ^2 \vspace{1ex}\\
	& \leq & \mathcal{L}_{\tau}^{\rho}(x^0, y^0; y^0) = \phi(x^0) + \langle y^0  ,F(x^0)\rangle+\frac{\rho}{2}\Vert F(x^0)\Vert ^2 \vspace{1ex}\\
	& {\overset{{\eqref{eq:pert_AL_at_x0}}}{\leq}} & \bar{\alpha}+\frac{1}{2\rho}\Vert y^0\Vert ^2+c_0.
\end{array}    
\end{equation}
Furthermore, utilizing the definition \eqref{eq:Lyapunov_func} of $\Pc_{k}$ at $k=1$, we have:
\begin{equation}\label{eq:lm8_ine110}
\arraycolsep=0.2em
\begin{array}{lcl} 
	\Pc_1 &= & \mathcal{L}^{\tau}_{\rho}(x^1,y^1;y^0)- \frac{\tau(1-\tau)}{2\rho}\norms{y^{1}-y^0}^2  + \frac{2(1-\tau)^2}{\tau\rho}\norms{y^1-y^0}^2 \vspace{1ex}\\
& = & \mathcal{L}^{\tau}_{\rho}(x^1,y^1;y^0)-\mathcal{L}^{\tau}_{\rho}(x^1,y^0;y^0)+\mathcal{L}^{\tau}_{\rho}(x^1,y^0;y^0)-\mathcal{L}^{\tau}_{\rho}(x^0,y^0;y^0)+\mathcal{L}^{\tau}_{\rho}(x^0,y^0;y^0) \vspace{1ex}\\
&& - \frac{\tau(1-\tau)}{2\rho}\norms{y^{1}-y^0}^2  + \frac{2(1-\tau)^2}{\tau\rho}\norms{y^1-y^0}^2 \vspace{0.5ex}\\
&{\overset{\tiny\eqref{eq:descent_of_PAL}}{\leq}} & \left(1-\tau\right)\langle y^1-y^0  ,F(x^1)\rangle-\frac{\beta}{4}\Vert x^1-x^0\Vert ^2+\mathcal{L}^{\tau}_{\rho}(x^0,y^0;y^0) \vspace{1ex}\\
&& - {~} \frac{\tau(1-\tau)}{2\rho}\norms{y^{1}-y^0}^2  + \frac{2(1-\tau)^2}{\tau\rho}\norms{y^1-y^0}^2 \vspace{1ex}\\
& {\overset{}{\leq}} & \frac{\rho}{2}\Vert F(x^1)\Vert ^2+\frac{\left(1-\tau\right)^2}{2\rho}\left(1+\frac{4}{\tau}\right)\Vert y^1-y^0\Vert ^2+\mathcal{L}^{\tau}_{\rho}(x^0,y^0;y^0) \vspace{1ex}\\
&  {\overset{}{\leq}} & \frac{\rho\left(1-\tau\right)^2\left(2+\tau\right)}{\tau}\Vert F(x^1)\Vert ^2+\mathcal{L}^{\tau}_{\rho}(x^0,y^0;y^0).
\end{array}    
\end{equation}
From \eqref{eq:lm8_ine111}, we can derive that
\begin{equation*} 
\arraycolsep=0.2em
\begin{array}{lcl} 
	\frac{\rho}{6}\Vert F(x^1)\Vert ^2 & \leq & \bar{\alpha}+\frac{1}{2\rho}\Vert y^0\Vert ^2+c_0-\frac{\rho}{6}\Vert F(x^1)\Vert ^2-\langle y^0  ,F(x^1)\rangle-f(x^1)-g(x^1)-\frac{\rho}{6}\Vert F(x^1)\Vert ^2 \vspace{1ex}\\
	& = & \bar{\alpha}+\frac{1}{2\rho}\Vert y^0\Vert ^2+c_0-\frac{\rho}{6}\Vert F(x^1)+\frac{3y^0}{\rho}\Vert ^2+\frac{3}{2\rho}\Vert y^0\Vert ^2-f(x^1)-g(x^1)-\frac{\rho}{6}\Vert F(x^1)\Vert ^2 \vspace{1ex}\\
	& {\overset{\tiny{\rho \geq 3\rho_0}}{\leq}} & \bar{\alpha}+\frac{1}{2\rho}\Vert y^0\Vert ^2+c_0+\frac{3}{2\rho}\Vert y^0\Vert ^2-f(x^1)-g(x^1)-\frac{\rho_0}{2}\Vert F(x^1)\Vert ^2 \vspace{0.5ex}\\
	& {\overset{{\eqref{lem1}}}{\leq}} & \bar{\alpha}+\frac{2}{\rho}\Vert y^0\Vert ^2+c_0-\underline{\Phi} \overset{\tiny{\rho\geq1}}{\leq}  \bar{\alpha}+c_0-\underline{\Phi}+{2}\Vert y^0\Vert ^2.
\end{array}    
\end{equation*}
Substituting this inequality into \eqref{eq:lm8_ine110}, we get
\begin{equation*} 
\arraycolsep=0.2em
\begin{array}{lcl} 
	\Pc_1 & \leq & \frac{6\left(1-\tau\right)^2\left(2+\tau\right)}{\tau}(\bar{\alpha}+c_0-\underline{\Phi}+2\Vert y^0\Vert ^2)+\bar{\alpha}+c_0+\frac{1}{2\rho}\Vert y^0\Vert ^2 \vspace{1ex} \\
    & \leq &  \left(\frac{18\left(1-\tau\right)^2}{\tau}+1\right) \left(\bar{\alpha}+c_0+2\Vert y^0\Vert ^2\right)-\frac{18\left(1-\tau\right)^2}{\tau}\underline{\Phi} \; \overset{\tiny\eqref{eq:P_upper_bound}}{ = }  \; \bar{\Pc}.
\end{array}    
\end{equation*}
This inequality verifies that for $k = 1$, \eqref{eq:important} holds. Now, we assume that \eqref{eq:important} holds for some $k \geq 1$ (induction hypothesis).
We prove that it continues to hold for $k+1$. 
Indeed, since $x_{k}, x_{k+1}\in\mathcal{S}$, by \eqref{eq:decrease_Lyapunov} from Lemma \ref{le:decrease_of_P} together with the induction hypothesis, we easily get
\begin{equation*}
    \Pc_{k+1}\leq \Pc_{k}{\overset{{}}{\leq}}\Pc_u.
\end{equation*}
Hence, we conclude that $\Pc_k \leq \bar{\Pc}$ for all $k\geq 1$.  It remains to prove that $\Pc_{k}\geq \underline{\Pc}$ for all $k\geq1$, where $\underline{\Pc} > -\infty$ is fixed.  
Since $\mathcal{S}$ is compact and $\phi$ is continuous, without loss of generality, we assume $\phi(x) \geq 0$ for all $x \in \mathcal{S}$. 
Then, by the Cauchy-Schwarz inequality and the identity $\langle b, b - a \rangle = \frac{1}{2} \left[ \Vert b - a\Vert ^2 + \Vert b\Vert ^2 - \Vert a\Vert ^2\right]$, we have
\begin{equation*} 
\arraycolsep=0.2em
\begin{array}{lcl} 
	\Tc_{[3]} & := &  \langle (1 - \tau)(y^k - y^0)+y^0, F(x^k) \rangle - \frac{\tau(1 - \tau)}{2\rho}\Vert y^k-y^0\Vert ^2  \vspace{1ex} \\
	& \geq &  \langle (1 - \tau)(y^k - y^0) +y^0, F(x^k) \rangle - \frac{\tau(1 - \tau)}{\rho}\Vert y^k - y^0\Vert ^2  \vspace{1ex} \\
	& = & \langle (1 - \tau)(y^k - y^0), F(x^k) - \frac{\tau}{\rho} (y^k - y^0) \rangle + \langle y^0, F(x^k) \rangle  \vspace{1ex} \\
	&  \geq & \langle (1 - \tau)(y^k - y^0), \frac{y^k - y^{k-1}}{\rho} + \frac{\tau}{\rho} (y^{k-1} - y^0) - \frac{\tau}{\rho} (y^k -y^0) \rangle -\Vert F(x^k)\Vert \Vert y^0\Vert    \vspace{1ex} \\
	& \geq & \frac{(1 - \tau)^2}{\rho} \langle (y^k - y^0), (y^k - y^0) - (y^{k-1} - y^0) \rangle -\Delta \Vert y^0\Vert  \vspace{1ex} \\
	& =  & \frac{(1 - \tau)^2}{2\rho} \left[ \Vert y^k - y^{k-1}\Vert ^2 + \Vert y^k - y^0\Vert ^2 - \Vert y^{k-1} - y^0\Vert ^2 \right] -\Delta \Vert y^0\Vert   \vspace{1ex} \\
	& \geq &  \frac{(1 - \tau)^2}{2\rho} \left[ \Vert y^k - y^0\Vert ^2 - \Vert y^{k-1} - y^0\Vert ^2 \right] - \Delta \Vert y^0\Vert.
\end{array}    
\end{equation*}
This inequality implies that
\begin{equation*} 
	\langle (1 - \tau)(y^k - y^0)+y^0, F(x^k) \rangle - \frac{\tau(1 - \tau)}{2\rho}\Vert y^k-y^0\Vert ^2 + \Delta \Vert y^0\Vert  \geq  \frac{(1 - \tau)^2}{2\rho} \left( \Vert y^k\Vert ^2 - \Vert y^{k-1}\Vert ^2 \right).
\end{equation*}
Summing this inequality from $k=1$ to $k = K$ we get
\begin{equation*} 
\arraycolsep=0.2em
\begin{array}{lcl} 
	\Tc_{[4]} &:= & \sum_{k=1}^{K} \big[  \langle (1 - \tau)(y^k - y^0)+y^0, F(x^k) \rangle - \frac{\tau(1 - \tau)}{2\rho}\Vert y^k-y^0\Vert ^2 + \Delta \Vert y^0\Vert  \big] \vspace{1ex}\\
	& \geq & \frac{(1 - \tau)^2}{2\rho} \sum_{k=1}^{K}\left( \Vert y^k - y^0\Vert ^2 - \Vert y^{k-1} - y^0\Vert ^2 \right) \vspace{1ex} =  \frac{(1 - \tau)^2}{2\rho}  \Vert y^K- y^0\Vert ^2   \geq  0.
\end{array}    
\end{equation*}
By the definition \eqref{eq:Lyapunov_func} of $\Pc_k$, and $\phi(x^k) \geq 0$, we can show that
\begin{equation*} 
\arraycolsep=0.2em
\begin{array}{lcl} 
	\Pc_k + \Delta\Vert y^0\Vert & = & \mathcal{L}_{\tau}^{\rho}(x^k, y^k; y^0) - \frac{\tau(1-\tau)}{2\rho}\Vert y^k - y^0\Vert^2 + \frac{2(1-\tau)^2}{\tau\rho}\Vert \Delta{y}^{k-1}\Vert^2 + \Delta\Vert y^0\Vert \vspace{1ex}\\
	& = & \phi(x^k) + \frac{\rho}{2}\Vert F(x^k) \Vert^2 + \frac{2(1-\tau)^2}{\tau\rho}\Vert \Delta{y}^{k-1}\Vert^2   \vspace{1ex}\\
	&& + {~} \langle (1 - \tau)(y^k - y^0)+y^0, F(x^k) \rangle - \frac{\tau(1 - \tau)}{2\rho}\Vert y^k-y^0\Vert ^2 + \Delta \Vert y^0\Vert \vspace{1ex}\\
	& \geq & \langle (1 - \tau)(y^k - y^0)+y^0, F(x^k) \rangle - \frac{\tau(1 - \tau)}{2\rho}\Vert y^k-y^0\Vert ^2 + \Delta \Vert y^0\Vert.
\end{array}    
\end{equation*}
Summing up this inequality from $k=1$ to $k=K$ and using $\Tc_{[4]}$ above, then passing the limit as $K \to\infty$, we get
\begin{equation*}
\lim_{K \to \infty} \sum_{k=1}^{K} \left(\Pc_k +  \Delta \Vert y^0\Vert \right) = \sum_{k=1}^{\infty} \left(\Pc_k +  \Delta \Vert y^0\Vert \right) \geq 0.
\end{equation*}
This nonnegative sum combined with the fact that $\{\Pc_k\}$ is non-increasing (or equivalently, $\{ \Pc_k +  \Delta \Vert y^0\Vert \}$ is non-increasing), implies the existence of a lower bound $\Pc_L > -\infty\) such that 
\begin{equation*}
\Pc_l \leq \Pc_k + \Delta \Vert y^0\Vert,  \quad \forall k \geq 1.
\end{equation*}
Therefore, if we define $\underline{\Pc} := \Pc_l - \Delta\Vert y^0 \Vert > -\infty$, then we have $\Pc_k \geq \underline{\Pc}$ for all $k \geq 0$.
\Eproof
\endproof 
%%%% End of Proof.

\medskip

%%%% Proof of Lemma 9.
\proof{\textbf{Proof of Lemma \ref{le:bounded_grad}.}} 
Let $v^{k+1} = (v_x^{k+1}, v_y^{k+1}, v_{y'}^{k+1}) \in \partial P(x^{k+1},y^{k+1},y^{k};y^0)$.
Then, from the definition \eqref{eq:Lyapunov_func}, we can show that
\begin{equation}\label{eq:v_subgrad}
\arraycolsep=0.2em
\begin{array}{lclcl}
    v_x^{k+1}    & \in & \partial_x\mathcal{L}^{\tau}_{\rho}(x^{k+1}, y^{k+1};y^0) & = &  \partial g(x^{k+1}) + \nabla_x\psi^{\tau}_{\rho}(x^{k+1}, y^{k+1};y^0), \vspace{1ex}\\ 
    v_y^{k+1}    & = & \nabla_y P(x^{k+1}, y^{k+1},y^k;y^0) & = & \frac{(4+\tau)(1-\tau)^2}{\tau\rho}\left(y^{k+1}-y^k\right), \vspace{1ex} \\ 
    v_{y'}^{k+1} & = & \nabla_{y'} P(x^{k+1}, y^{k+1},y^k;y^0) & = &  \frac{4(1-\tau)^2}{\tau\rho}\left(y^k-y^{k+1}\right). 
\end{array}
\end{equation}
Next, by the optimality condition of \eqref{eq:LIPAL_subprob}, it follows that there exists $s_g^{k+1} \in \partial g(x^{k+1})$ such that
\begin{equation*}
s_g^{k+1} +\nabla_x\psi^{\tau}_{\rho}(x^{k}, y^{k};y^0) + \beta\Delta x^{k} +\rho J_F(x^{k})^TJ_F(x^{k})\Delta x^{k}=0.
\end{equation*}
Combining this expression and the first line of \eqref{eq:v_subgrad}, and using the triangle inequality, we can derive
\begin{equation*}
    \Vert v_x^{k+1}\Vert  \leq \Vert \nabla_x\psi^{\tau}_{\rho}(x^{k+1}, y^{k+1};y^0)-\nabla_x\psi^{\tau}_{\rho}(x^{k}, y^{k};y^0)\Vert  + (\beta+\rho M_F^2)\Vert \Delta x^{k}\Vert .
\end{equation*}
Since $\nabla_x\psi^{\tau}_{\rho}$ is locally Lipschitz continuous and  the sequence $\{ z_k := (x^{k}, y^{k}) \}_{k\geq1}$ is bounded (recall that $\{x^k\} \subset \mathcal{S}$ is bounded  and $\{y^k\}_{k\geq0}$ is bounded by Lemma \ref{le:lambda_bou1}), it follows that
\begin{equation}\label{eq:bound_vx}
	\Vert v_x^{k+1}\Vert  \leq (L^{\tau}_{\rho}+\beta+\rho M_F^2)\Vert x^{k+1}-x^k\Vert +(1-\tau)M_F\left\Vert y^{k+1}-y^k\right\Vert . 
\end{equation}
Next, from the second and third  lines of \eqref{eq:v_subgrad}, we also have
\begin{equation}\label{eq:bound_vy}
\arraycolsep=0.2em
\begin{array}{llcl}
	& \Vert v_y^{k+1}\Vert  & = & \frac{(4+\tau)(1-\tau)^2}{\tau\rho} \Vert y^{k+1} - y^k \Vert, \vspace{1ex}\\
	& \Vert v_{y'}^{k+1}\Vert  & = &  \frac{4(1-\tau)^2}{\tau\rho} \Vert y^{k+1} - y^k \Vert.
\end{array}
\end{equation}
Combining \eqref{eq:bound_vx} and \eqref{eq:bound_vy}, we can show that
\begin{equation*}
	\Vert v^{k+1}\Vert  \leq \big( L^{\tau}_{\rho}+\beta+\rho M_F^2 \big) \Vert x^{k+1} - x^k\Vert  + \big( \tfrac{(8+\tau)(1-\tau)^2}{\tau\rho}+(1-\tau)M_F \big) \Vert y^{k+1} - y^k\Vert,
\end{equation*}
which proves \eqref{eq:bounded_subgrad_of_P}.
\Eproof
\endproof 
%%%% End of Proof.

\medskip

%%% Proof of Lemma 10.
\proof{\textbf{Proof of Lemma \ref{le:added_lemma}.}}
$\mathrm{(i)}$~
Since $\{u^{k}\}$ is bounded, there exists a convergent subsequence  $\{u^{k}\}_{k \in \mathcal{K}}$ satisfying $\lim_{k\in\mathcal{K}, k \to \infty}{u^{k}} = u^{*}$. 
Hence $\Omega$ is nonempty. 
Moreover, $\Omega$ is compact since $\{ u^k \}$ is bounded.

On the other hand, for any $u^{*}\in\Omega$, there exists a $\{u^{k}\}_{k \in \mathcal{K}}$ such that $\lim_{k\in\mathcal{K}, k \to \infty}{u^{k}}=u^{*}$.
By Lemma \ref{le:bounded_grad} and \eqref{eq:FO_convergence_limits}, it follows that there exists $v^{*}\in\partial{P}(u^{*})$ such that
\begin{equation*}
	\Vert v^{*}\Vert =\lim_{ i \to \infty }{\Vert v^{k + 1}\Vert } \overset{\tiny\eqref{eq:FO_convergence_limits}}{=} 0.
\end{equation*}
Hence, $u^{*}\in \mathrm{crit}(P)$ and $0\leq \lim_{k\to\infty}{\mathrm{dist}(u^{k},\Omega)}\leq\lim_{k \in \mathcal{K}, k \to \infty}{\mathrm{dist}(u^{k},\Omega)}=\mathrm{dist}(u^{*},\Omega)=0$. Since $\Omega \subseteq \mathrm{crit}(P)$, it also follows that $\lim_{k\to\infty}{\mathrm{dist}(u^{k}, \mathrm{crit}(P) )} = 0$.

$\mathrm{(ii)}$~
To prove Statement $\mathrm{(ii)}$, we follow similar arguments as in \cite{CohHal:21}. 
Since $\{u^{k}\}$ is bounded, there exists a convergent subsequence $\{u^{k}\}_{k \in \mathcal{K}}$ such that $\lim_{k\in\mathcal{K}, k \to \infty} u^{k} = u^{*}$.
Moreover, since \(P(\cdot)\) is a proper lsc function, we have 
\begin{equation*}
\liminf_{k \in \mathcal{K}, k \to \infty} P(u^{k}) \geq P(u^{*}) = \Pc^{*}.
\end{equation*}
Since $P(\cdot) - g(\cdot)$ is continuous, we get
\begin{equation*}
\lim_{k \in \mathcal{K}, k \to \infty} (P(u^{k}) - g(x^{k})) = P(u^{*}) - g(x^{*}).
\end{equation*}
Thus it remains to prove the following inequality:
\begin{equation*}
\limsup_{k \in \mathcal{K}, k \to \infty} g(x^{k}) \leq g(x^{*}).
\end{equation*}
From \eqref{eq:LIPAL_subprob}, for any $k \in \mathcal{K}$, we have $\mathcal{Q}_{k}(x^k) \leq \mathcal{Q}_{k}(x^{\star})$, leading to 
\begin{equation*} 
\arraycolsep=0.2em
\begin{array}{lclcl}
 g(x^k) - g(x^{*}) & \leq &  \big\langle \nabla_x \psi^{\tau}_{\rho}(x^{k-1}, y^{k-1}; y^0), x^{*} - x^k \big\rangle  \vspace{1ex}\\
 && + {~} \big\langle \frac{1}{2} \left( \rho J_F(x^{k-1})^T J_F(x^{k-1}) + \beta I_n \right) ( x^{*} - x^{k-1} + x^{k} - x^{k-1} ), x^{*} - x^k \big\rangle \vspace{1ex} \\
   & \leq &   \Vert  \nabla_x \psi^{\tau}_{\rho}(x^{k-1}, y^{k-1}; y^0)\Vert \Vert x^k - x^{*}\Vert  \vspace{1ex}\\
   && + {~} \frac{1}{2} \Vert \big( \rho J_F(x^{k-1})^T J_F(x^{k-1}) + \beta I_n \big) ( x^{*} - x^{k-1} + x^{k} - x^{k-1} ) \big\Vert  \Vert x^k - x^{*} \Vert  \vspace{1ex} \\
   &  \leq &  c_k \Vert x^k -  x^{*}\Vert,
\end{array}
\end{equation*}
where  the second inequality follows from the Cauchy-Schwarz inequality, the third inequality is a consequence of the triangle inequality, and
\begin{equation*}
c_k \triangleq \left\Vert  \nabla_x \psi^{\tau}_{\rho}(x^{k-1}, y^{k-1}; y^0) \right\Vert  + \frac{1}{2} \left\Vert  \rho J_F(x^{k-1})^T J_F(x^{k-1}) + \beta I_n \right\Vert  \left( \Vert x^{*} - x^{k-1}\Vert  + \Vert x^{k} - x^{k-1}\Vert  \right).
\end{equation*}
Taking the limit over $\mathcal{K}$ as $k \to \infty$, we obtain
\begin{equation*}
    \limsup_{k\in\mathcal{K}, k \to\infty} g(x^{k}) \leq g(x^{*}) + \limsup_{k\in\mathcal{K}, k \to\infty} c_k \Vert x^k - x^{*} \Vert.
\end{equation*}
Thus we only need to prove that 
\begin{equation*}
    \limsup_{k\in\mathcal{K}, k \to\infty} c_k \Vert x^k - x^{*} \Vert  = 0.
\end{equation*}
Since $\lim_{k\in\mathcal{K}, k \to\infty} x^k = x^{*}$, it suffices to show that 
\begin{equation*}
    \limsup_{k\in\mathcal{K}, k \to\infty} c_k < +\infty.
\end{equation*}
From Theorem \ref{th:complex_bound1}, we have $\lim_{k\in\mathcal{K}, k \to\infty} \Vert x^k - x^{k-1}\Vert  = 0$. 
Therefore, the subsequence $\{ x^{k-1} \}_{k \in \mathcal{K}}$ also converges to $x^{*}$. 
Moreover, since $F$ is $L_F$-smooth and $\rho, \beta < \infty$, we can show that
\begin{equation*}
\lim_{k\in\mathcal{K}, k \to\infty} \left\Vert  \rho J_F(x^{k-1})^T J_F(x^{k-1}) + \beta I_n \right\Vert  = \left\Vert  \rho J_F(x^{*})^T J_F(x^{*}) + \beta I_n \right\Vert  < \infty.
\end{equation*}
Consequently, we obtain
\begin{equation*}
\lim_{k\in\mathcal{K}, k \to\infty} \frac{1}{2} \left\Vert  \rho J_F(x^{k-1})^T J_F(x^{k-1}) + \beta I_n \right\Vert  \left( \Vert x^{*} - x^{k-1}\Vert  + \Vert x^{k} - x^{k-1}\Vert  \right) = 0.
\end{equation*}
Furthermore, since $\nabla_x \psi^{\tau}_{\rho}$ is $L_{\tau}^{\rho}$-Lipschitz continuous, we get
\begin{equation*}
\lim_{k\in\mathcal{K}, k \to\infty} \left\Vert  \nabla_x \psi^{\tau}_{\rho}(x^{k-1}, y^{k-1}; y^0) \right\Vert  = \left\Vert  \nabla_x \psi^{\tau}_{\rho}(x^{*}, y^{*}; y^0) \right\Vert  < \infty.
\end{equation*}
Therefore, we conclude that
\begin{equation*}
	\lim_{k\in\mathcal{K}, k \to\infty} c_k = \limsup_{k\in\mathcal{K}} c_k < \infty.
\end{equation*}
As a result, any converging subsequence $\{P(u^{k}) = \Pc_k\}_{k \in \mathcal{K}}$ must converge to the same limit $\Pc^{*}$.

\medskip 

$\mathrm{(iii)}$~
Let $(x^{*}, y^{*}, \hat{y}^{*}) \in\mathrm{crit}(P)$ be a critical point of $P$. 
Then, it follows that
\begin{equation*} 
\arraycolsep=0.2em
\begin{array}{lclcl}
0 & \in & \partial_x P(x^{*}, y^{*}, \hat{y}^{*}; y^0) & = & \partial_x \mathcal{L}^{\tau}_{\rho}(x^{*}, y^{*}; y^0), \vspace{1ex} \\
0 & = & \nabla_{y}{P}(x^{*}, y^{*}, \hat{y}^{*}; y^0) & = & (1-\tau)\big(F(x^{*})-\frac{\tau}{\rho}(y^{*} - y^0 ) \big)+\frac{4(1-\tau)^2}{\tau\rho}\left(y^{*} - \hat{y}^{*} \right), \vspace{1ex}\\
0 & = & \nabla_{y'}{P}(x^{*}, y^{*}, \hat{y}^{*}; y^0) & = & \frac{4(1-\tau)^2}{\tau\rho} (\hat{y}^{*} - y^{*}).
\end{array}
\end{equation*}
The last line shows that $y^{*} = \hat{y}^{*}$.
Combining this relation and the first two lines, we get
\begin{equation*}
    -\nabla f(x^{*}) - J_F^Ty^{*} \in\partial g(x^{*}) \quad\textrm{and} \quad   F(x^{*}) = \frac{\tau}{\rho}(y^{*} - y^0).
\end{equation*}
Hence, we conclude that $(x^{*}, y^{*})$ is an $\epsilon$-first-order optimal solution to \eqref{eq:nlp_gen}, where $\epsilon =  \frac{\tau}{\rho}\Vert y^{*} - y^0\Vert$.
\Eproof
\endproof 
%%% End of Proof.

\medskip

%%% Proof of Lemma 11.
\proof{\textbf{Proof of Lemma \ref{le:finite_length}.}}
From \eqref{eq:lm10_decrease_Lyapunov_KL}, we have 
\begin{equation}\label{eq:lm11_llyap}
    \Pc_{k+1}-\Pc_{k} {\overset{{\eqref{eq:lm10_decrease_Lyapunov_KL}}}{\leq}} -\ubar{\gamma}\left(\Vert \Delta x^{k}\Vert ^2+\Vert \Delta y^{k}\Vert ^2\right)=-\ubar{\gamma}\Vert z^{k+1}-z^{k}\Vert ^2.
\end{equation}
 Since $ \Pc_{k}\to \Pc^{*}$ and  $\{\Pc_{k} \}$ is monotonically decreasing to $\Pc^{*}$,  it implies that the error sequence $\{\mathcal{E}_{k}\}$ is non-negative, monotonically decreasing and convergent to zero. 
 We consider two cases as follows.

%%% Case 1.
\vspace{0.5ex}
\noindent \textbf{{Case 1.}} 
There exists  $k_1\geq 1$ such that $\mathcal{E}_{k_1}=0$.
Then, it is obvious that $\mathcal{E}_{k}=0$ for all $k \geq k_1$ and using \eqref{eq:lm11_llyap}, we have
\begin{equation} \label{eq:KL_bound}
	\Vert z^{k+1}-z^{k}\Vert ^2\leq\frac{1}{\ubar{\gamma}}(\mathcal{E}_{k}-\mathcal{E}_{k+1}) = 0, \quad\forall k\geq k_1. 
\end{equation}
Since $\{z^{k}\}$ is bounded, one can show that
\begin{equation*}
	\sum_{k=0}^{\infty}{\Vert \Delta x^{k}\Vert +\Vert \Delta y^{k}\Vert } 
	\overset{\eqref{eq:KL_bound}}{=}\sum_{k=0}^{k_1}{\Vert \Delta x^{k}\Vert +\Vert \Delta y^{k}\Vert } <  +\infty.
\end{equation*}
 %%%% Case 2.
\noindent \textbf{{Case 2.}} 
The error $\mathcal{E}_{k} > 0$ for all $k \geq 1$. 
Then,  there exists  $k_1 := k_1(\epsilon, \varepsilon)\geq 1$  such that $\forall k \geq k_1$, we have $\mathrm{dist}(u^{k},\Omega) \leq \epsilon$,  $\Pc^{*} < P(u^{k}) < \Pc^{*} + \varepsilon$ and
\begin{equation}\label{eq:lm11_KL}
	\varphi'(\mathcal{E}_{k})\Vert \nabla{P}(x^{k},y^{k},y^{k-1})\Vert \geq 1, \quad \forall \nabla{P}(x^{k},y^{k},y^{k-1}) \in \partial{P}(x^{k},y^{k},y^{k-1}).
\end{equation}
where $\epsilon>0$, $\varepsilon>0$, and $\varphi\in\Psi_{\varepsilon}$ are  defined from the KL property of  $P$ on $\Omega$ (see Definition \ref{def2}). Since $\varphi$ is concave, we have $\varphi(\mathcal{E}_{k})-\varphi(\mathcal{E}_{k+1}) \geq \varphi'(\mathcal{E}_{k})(\mathcal{E}_{k}-\mathcal{E}_{k+1})$. 
Then, from \eqref{eq:lm11_llyap} and \eqref{eq:lm11_KL}, one can~show:
 \begin{equation*}
 \arraycolsep=0.2em
 \begin{array}{lcl}
 \Vert z^{k+1}-z^{k}\Vert ^2  & {\overset{\eqref{eq:lm11_KL}}{\leq}} & \varphi'(\mathcal{E}_{k})\Vert z^{k+1}-z^{k}\Vert ^2\Vert \nabla{P}(x^{k},y^{k},y^{k-1})\Vert \vspace{1ex} {\overset{\eqref{eq:lm11_llyap}}{\leq}}  \frac{1}{\ubar{\gamma}}\varphi'(\mathcal{E}_{k})(\mathcal{E}_{k}-\mathcal{E}_{k+1})\Vert \nabla{P}(x^{k},y^{k},y^{k-1})\Vert \vspace{1ex}\\
 &\leq & \frac{1}{\ubar{\gamma}}\left(\varphi(\mathcal{E}_{k})-\varphi(\mathcal{E}_{k+1})\right)\Vert \nabla{P}(x^{k},y^{k},y^{k-1})\Vert.
 \end{array}
 \end{equation*}
Since $ \Vert \Delta z^{k}\Vert ^2={\Vert \Delta x_{k}\Vert ^2+\Vert \Delta y_{k}\Vert ^2}$ and using the fact that for any $a,b,c,d\geq0$, if $ {a^2+b^2}\leq c\times d$, then $ (a+b)^2\leq 2(a^2+b^2)\leq 2c\times d\leq(c+d)^2$, for any $\theta > 0$, we can show that
\begin{equation*} 
	\Vert \Delta x^{k}\Vert +\Vert \Delta y^{k}\Vert   \leq \frac{\theta}{\ubar{\gamma}}\big[ \varphi(\mathcal{E}_{k}) - \varphi(\mathcal{E}_{k+1}) \big] +\frac{1}{\theta}\Vert \nabla{P}(x^{k},y^k,y^{k-1})\Vert .
\end{equation*}
Furthermore, from Lemma \ref{le:bounded_grad}, there exists $v^{k} := \nabla{P}(x^{k},y^{k},y^{k-1}) \in\partial P(x^{k},y^{k},y^{k-1})$ such that
\begin{equation*}
	\Vert \nabla{P}(x^{k},y^{k},y^{k-1}) \Vert = \Vert v^{k}\Vert \leq c_1\left(\Vert \Delta x^{k-1}\Vert +\Vert \Delta y^{k-1}\Vert \right).
\end{equation*}
Combining both inequalities, we get
\begin{equation*} 
	\Vert \Delta x^{k}\Vert +\Vert \Delta y^{k}\Vert  \leq \frac{\theta}{\ubar{\gamma}}\left(\varphi(\mathcal{E}_{k})-\varphi(\mathcal{E}_{k+1})\right)+\frac{c_1}{\theta}\left(\Vert \Delta x^{k-1}\Vert +\Vert \Delta y^{k-1}\Vert \right). 
\end{equation*}
Now, let us choose $\theta>0$ such that $0< \frac{c_1}{\theta} <1$ and denote  $\delta_0 := 1-\frac{c_1}{\theta}>0$. 
Then,  summing up the above inequality from $k =  k_1$ to $k=K$ and using the property: $\sum_{k=k_1}^{K}{\Vert \Delta x^{k-1}\Vert +\Vert \Delta y^{k-1}\Vert }\leq\sum_{k=k_1}^{K}{\Vert \Delta x^{k}\Vert +\Vert \Delta y^{k}\Vert }+\Vert \Delta x^{k_1-1}\Vert +\Vert \Delta y^{k_1-1}\Vert $, we get  
\begin{equation}\label{eq:lm11_used_next_lemma}
	\sum_{k=k_1}^{K}{\Vert \Delta x^{k}\Vert +\Vert \Delta y^{k}\Vert }   \leq \frac{\theta}{\ubar{\gamma}\delta_0}\varphi(\mathcal{E}_{k_1})+\frac{c_1}{\eta\delta_0}\Big(\Vert \Delta x^{k_1-1}\Vert +\Vert \Delta y^{k_1-1}\Vert \Big).
\end{equation}
Clearly, the right-hand side of this inequality is bounded for all $K \geq k_1$. 
Letting  $K\to\infty$, we obtain	$\sum_{k=k_1}^{\infty}{\Vert \Delta x^{k}\Vert +\Vert \Delta y^{k}\Vert }<\infty$.  Since $\{(x^{k},y^{k})\}$ is bounded, it follows that
\begin{equation*}
	\sum_{k=0}^{k_1-1}{\Vert \Delta x^{k}\Vert +\Vert \Delta y^{k}\Vert }<\infty.
\end{equation*}
Summing up the last two expressions, we conclude that $ \sum_{k=0}^{\infty}{\Vert \Delta x^{k}\Vert +\Vert \Delta y^{k}\Vert } < \infty$. 
In both cases, we obtain the finite length statement in \eqref{eq:lm11_finite_length}.  Finally, the finite length sum in \eqref{eq:lm11_finite_length} also implies that  $\{z^k\}_{k\geq0}$ is a Cauchy sequence and thus it converges. 
Moreover, from Theorem \ref{th:complex_bound1}, $\{z^k\}$  converges to $z^{*}=(x^{*},y^{*})$ such that it is an $\epsilon$-first-order optimal solution to \eqref{eq:nlp_gen}, where $\epsilon := \frac{\tau}{\rho}\Vert y^{*}-y^0\Vert$.
\Eproof
\endproof 
%%% End of Proof.

\medskip

 %%% Proof of Lemma 12.
 \proof{\textbf{Proof of Lemma \ref{le:bound_of_diff_zk}.}} 
By Lemma \ref{le:decrease_of_P}, $\{\Pc_k\}$ is monotonically decreasing.
Consequently, $\{\mathcal{E}_{k}\}$ is  monotonically decreasing and nonnegative. From \eqref{eq:lm11_llyap} and the nonnegativity of  $\{\mathcal{E}_{k}\}$, we have: 
\begin{equation}\label{eq:lm12_lmit2}
	\Vert \Delta x^{k}\Vert + \Vert \Delta y^{k}\Vert \leq  \frac{\sqrt{2\mathcal{E}_{k}}}{\sqrt{\underline{\gamma}}} \quad \forall k\geq 1.
\end{equation}
 Since $\Pc_{k}\to \Pc^{*}$, ${u^{k}}\to u^{*}$, and $P(\cdot)$ satisfies the KL property at $u^{*}$, there exists $k_1 := k_1(\epsilon,\varepsilon)\geq 1$   such that for all $k \geq k_1$, we have $\Vert u^{k}-u^{*}\Vert \leq \epsilon$ and $\Pc^{*}<\Pc_{k}<\Pc^{*}+\varepsilon$, and 
\begin{equation}\label{eq:lm12_KL2}
	\varphi'(\mathcal{E}_{k})\Vert \nabla{P}(x^{k},y^{k},y^{k-1};y^0)\Vert \geq 1, \quad \forall \nabla{P}(x^{k},y^{k},y^{k-1};y^0) \in \partial P(x^{k},y^{k},y^{k-1};y^0).
\end{equation}
Utilizing the same arguments as \textbf{Case 2} in the proof of Lemma \ref{le:finite_length}, the relation \eqref{eq:lm11_used_next_lemma} follows.
Hence, by the triangle inequality, we get for any $k \geq k_1$, we can show that
\begin{equation*}
	\Vert z^{k}-z^{*}\Vert \leq  \sum_{l\geq k}{\Vert z^{l}-z^{l+1}\Vert }  
 \leq \sum_{l\geq k}{\Vert \Delta x^{l}\Vert +\Vert \Delta y^{l}\Vert } \vspace{1ex} \\
	\leq \frac{\theta}{\ubar{\gamma}\delta_0}\varphi(\mathcal{E}_{{k}}) +\frac{c_1}{\theta\delta_0}\left(\Vert \Delta x^{k-1}\Vert +\Vert \Delta y^{k-1}\Vert \right).
\end{equation*}
Combining this inequality and \eqref{eq:lm12_lmit2}, we get
\begin{equation*}
	\Vert z^{k}-z^{*}\Vert  \leq \frac{\theta}{\ubar{\gamma}\delta_0}\varphi(\mathcal{E}_{{k}}) +\frac{c_1}{\theta\delta_0}\sqrt{\frac{2\mathcal{E}_{k-1}}{\ubar{\gamma}}} 
    \leq C \max\{\varphi(\mathcal{E}_{{k}}),\sqrt{\mathcal{E}_{k-1}}\},
\end{equation*}
which proves \eqref{eq:bound_of_diff_zk}, where $C \triangleq \max \big\{ \frac{\theta}{\ubar{\gamma}\delta_0}, \frac{c_1}{\theta\delta_0}\sqrt{\frac{2}{\ubar{\gamma}}} \big\}$.
 \Eproof
\endproof 
%%% End of Proof.

\medskip

%%% Proof of Lemma 18.
\proof{\textbf{Proof of Lemma \ref{lemma18}.}}
Note that \eqref{ours_local}  trivially implies \eqref{CQ_regularity_neighboor}. 
We now proceed to prove the other  implication by contradiction. 
Assume that $J_F(x)$ has full row rank for any $x \in V \cap \dom \varphi$ and that \eqref{CQ_regularity_neighboor} holds. 
Further, suppose that for any $\sigma > 0$, there exists $\bar y_{\sigma} \neq 0_m$ such that
\begin{equation*}
	\sigma \Vert \bar y_{\sigma}\Vert  > \mathrm{dist}\left(-J_F(x)^T \bar y_{\sigma}, \partial^{\infty} g(x)\right).
\end{equation*}
Since $\partial^{\infty} g(x)$ is a cone and $\bar y_{\sigma} \neq 0_m$,  this implies that for any $\sigma > 0$, there exists $\|y_{\sigma}\|=1$ such that
\begin{equation*}
	\sigma   > \mathrm{dist}\left(-J_F(x)^T y_{\sigma}, \partial^{\infty} g(x)\right).
\end{equation*}
Since $J_F$ has full row rank and it is continuous  on $V \cap \dom \varphi$,  there exists  constant $\gamma > 0$ such that
\begin{equation}\label{row_rank}
	\Vert J_F(x)^T y\Vert  \geq \gamma \Vert y\Vert   \quad \forall y \in \mathbb{R}^m.
\end{equation}
Moreover, since $g$ is convex, it follows that $\partial^{\infty} g(x) = N_{\dom g}(x)$, where $\dom g$ is a convex  \cite{RocWet:98}. 
We distinguish three cases  for $-J_F(x)^T y_{\sigma}$: \\
\noindent\textbf{Case (i)}:  
if $-J_F(x)^T y_{\sigma} \in N_{\dom g}(x)$, then we have
\begin{equation*}
	\mathrm{dist}\left(-J_F(x)^T y_{\sigma}, \partial^{\infty} g(x)\right) = 0,
\end{equation*}
and since  $\|y_{\sigma}\| = 1$,  we get a contradiction of  \eqref{CQ_regularity_neighboor}.

\noindent \textbf{Case (ii)}:  if $-J_F(x)^T y_{\sigma}$ in the tangent cone of $\dom g$, then we have
\begin{equation*}
	\mathrm{dist}\left(-J_F(x)^T y_{\sigma}, \partial^{\infty} g(x)\right) = \Vert -J_F(x)^T y_{\sigma}\Vert  \overset{\eqref{row_rank}}{\geq} \gamma \Vert y_{\sigma}\Vert .
\end{equation*}
This further  yields $\sigma \Vert y_{\sigma} \Vert  > \gamma \Vert y_{\sigma}\Vert$, or equivalently, $\gamma < \sigma$ for all  $\sigma > 0$,  which is impossible (take $\sigma \to 0^{+}$), thus contradicting the full row rank assumption.

\noindent\textbf{Case (iii)}:  if $-J_F(x)^T y_{\sigma}$ forms an acute angle with its projection onto $N_{\dom g}(x)$, then  we have
\begin{equation*}
	\mathrm{dist}\left(-J_F(x)^T y_{\sigma}, \partial^{\infty} g(x)\right) = \Vert -J_F(x)^T y_{\sigma}\Vert  |\sin(\omega_{\sigma})| \overset{\eqref{row_rank}}{\geq} \gamma |\sin(\omega_{\sigma})| \Vert y_{\sigma}\Vert ,
\end{equation*}
where $\omega_{\sigma}$ is the angle between $-J_F(x)^T y_{\sigma}$ and its projection onto $N_{\dom g}(x)$. 
This implies that
\begin{equation*}
	\vert \sin(\omega_{\sigma}) \vert < \frac{\sigma}{\gamma} \quad \forall \sigma > 0.
\end{equation*}
Taking the limit as $\sigma \to 0^+$, we obtain $\vert \sin(\omega_{\sigma}) \vert \to 0$, which implies $\mathrm{dist}\left(-J_F(x)^T y_{\sigma}, \partial^{\infty} g(x)\right) \to 0$,
contradicting \eqref{CQ_regularity_neighboor} since $\|y_{\sigma}\|=1$.
This completes our proof.
\Eproof
\endproof 
%%% End of Proof.

\medskip

%%% Proof of Lemma 19.
\proof{\textbf{Proof of Lemma \ref{le:limit_feasibility}.}}
Our proof  follows a similar reasoning as in  Proposition 3.6   \cite{DemJia:23}. 
Using the optimality of $x^{k+1}$ in \eqref{eq:LIPAL_subprob}, we arrive at \eqref{feasibility_measure}, which we can rewrite as
\begin{equation} \label{eq:lm19_proof1}
\mathrm{dist}\left(-J_F(x^{k+1})^TF(x^{k+1}),\partial^{\infty}g(x^{k+1})\right) \leq  \frac{\tau(M_f+1+\Vert y^0\Vert )}{\rho}  + \frac{2\tau\beta}{\rho}\Vert \Delta x^{k}\Vert  +  \frac{(1-\tau)M_F}{\rho}\Vert \Delta y^k\Vert .
\end{equation}
Moreover, by the decrease of the Lyapunov function $\{\Pc_k\}$ due to Lemma \ref{le:decrease_of_P}, and the fact that $\underline{\Pc} \leq \Pc_k \leq \bar{\Pc}$, we obtain \eqref{eq:FO_convergence_limits}. 
Hence,  passing to the limit in \eqref{eq:lm19_proof1}, we get
\begin{equation*}
	\mathrm{dist}\left(-J_F(x^{*})^TF(x^{*}),\partial^{\infty}g(x^{*})\right)=\lim_{k\in\mathcal{K}}{\mathrm{dist}\left(-J_F(x^{k+1})^TF(x^{k+1}),\partial^{\infty}g(x^{k+1})\right) }\overset{\eqref{eq:FO_convergence_limits}}{\leq}\frac{\tau(M_f+1+\Vert y^0\Vert )}{\rho}.
\end{equation*}
From the choice of $\rho$, it follows that
\begin{equation}\label{eq:lm19_epsilon_feasiblility}
	\mathrm{dist}\left(-J_F(x^{*})^TF(x^{*}),\partial^{\infty}g(x^{*})\right) \leq\epsilon.
\end{equation}
On the other hand, the feasible problem can be rewritten as
\begin{equation*}
	\min_{(x, \alpha) \in\epi g } \frac{1}{2} \Vert F(x) \Vert ^2 =\min_{(x,\alpha)} \frac{1}{2}\Vert F(x) \Vert ^2 + \delta_{\epi g}(x,\alpha),
\end{equation*}
and the corresponding first-order optimality condition reads as
\begin{equation*}
0 \in \begin{pmatrix}  J_F(x)^TF(x) \\ 0 \end{pmatrix} + N^{\text{lim}}_{\text{epi} g} (x,g(x)) \quad \iff \quad  \mathrm{dist}\left(-J_F(x)^TF(x),\partial^{\infty}g(x)\right)=0.
\end{equation*}
Hence, $x^{*}$ in \eqref{eq:lm19_epsilon_feasiblility} is an $\epsilon$-first-order solution of the feasible problem.
Moreover, if $g$ is locally Lipschitz continuous at $x^{*}$, then $\partial^{\infty}g(x^{*})=\{0\}$, and thus \eqref{eq:lm19_epsilon_feasiblility} becomes
\begin{equation*} 
\left\Vert J_F(x^{*})^TF(x^{*})\right\Vert  \leq\epsilon,
\end{equation*}
which corresponds to the $\epsilon$-first-order optimal solution to $ \min_{x \in \mathbb{R}^n } \frac{1}{2} \Vert F(x) \Vert ^2$.  
\Eproof
\endproof 
%%% End of Proof.

%%%%%%%%%%%%%%%%%%%%%%%%%%%%%%%%%%%%%%%%%%%%


\begin{thebibliography}{plain}

\bibitem{AttBol:13}{ H. Attouch,  J. Bolte and  B. Svaiter, \emph{Convergence of descent methods for semi-algebraic and tame
problems: proximal algorithms, forward-backward splitting, and regularized Gauss–Seidel methods}, Mathematical Programming, 137: 91–129, 2013.}

\bibitem{Ber:96}{ D.P. Bertsekas, \emph{Constrained Optimization and Lagrange Multiplier Methods}, Athena Scientific,  1996.}

\bibitem{BirGar:16} 
E.G. Birgin, J.L. Gardenghi, J.M. Martínez, S.A. Santos and Ph.L. Toint, \emph{Evaluation complexity for nonlinear constrained optimization using unscaled KKT conditions and high-order models}, SIAM Journal on Optimization, 26(2): 951-967, 2016.

\bibitem{BolDan:07}
J. Bolte, A. Daniilidis, A. Lewis and M. Shiota, \emph{Clarke subgradients of stratifiable functions},  SIAM Journal on Optimization, 18(2): 556--572, 2007.

%\bibitem{BolDan:08} J. Bolte, A. Daniilidis, O. Ley, and L. Mazet, Laurent, \emph{Characterizations of Lojasiewicz inequalities and applications}, ?????, 2008.

%\bibitem{BolSab:14}{ J. Bolte, S. Sabach and M. Teboulle,  \emph{Proximal alternating linearized minimization for nonconvex and nonsmooth problems}, Mathematical Programming, 146(1–2):  459–494, 2014.}

%\bibitem{BotNgu:20}{ R.I. Bot and D.K. Nguyen, \emph{The proximal alternating direction method of multipliers in the nonconvex setting: convergence analysis and rates}, Mathematics of Operations Research, 45(2): 682–712, 2020.}

\bibitem{BoyPar:11}{ S. Boyd,  N. Parikh, E. Chu, B. Peleato and J. Eckstein,  \emph{Distributed optimization and statistical learning via the alternating direction method of multipliers}, Foundations and Trends in Machine Learning, 3(1): 1–122, 2011.}

\bibitem{BurMon:03}
S. Burer and R. D. Monteiro, \emph{A nonlinear programming algorithm for solving semidefinite programs via low-rank factorization}, Mathematical Programming, 95(2):   329–357, 2003.
 
\bibitem{CohHal:21}{ E. Cohen,  N. Hallak and M. Teboulle, \emph{A Dynamic Alternating Direction of Multipliers for Nonconvex Minimization with Nonlinear Functional Equality Constraints}, Journal of Optimization Theory and Applications, 193: 324–353, 2022.}


%\bibitem{Dan:89}{ Y. M. Danilin, \emph{Quadratic penalty methods based on linear approximation}, USSR Computational Mathematics and Mathematical Physics, 29, 132-139, 1989.}

%\bibitem{Dan:02}{ Y. M. Danilin, \emph{Linearization and Penalty Functions}, Cybernetics and Systems Analysis, 38, 691–702, 2002.}


\bibitem{DemJia:23} A. De Marchi, X. Jia,  C. Kanzow and   P. Mehlitz,  \emph{Constrained composite optimization and augmented Lagrangian methods}, Mathematical Programming, 201: 863–896, 2023.

%\bibitem{DonRoc:14} A. Dontchev and R. Rockafellar,  \emph{ Implicit functions and solution mappings. A view from variational analysis}, 2nd edition, 2014. 


\bibitem{ElbNec:23}  L. El Bourkhissi and I. Necoara, \emph{Convergence rates for an inexact linearized ADMM for nonsmooth optimization with nonlinear equality constraints}, provisionally accepted in Computational Optimization and Applications, 2024.

\bibitem{Fan:97}{ J. Fan, \emph{Comments on wavelets in statistics: a review}, Journal of the Italian Statistical Society, 6: 131–138, 1997.}

%\bibitem{Fes:20}{ J. A. Fessler, \emph{Optimization methods for magnetic resonance image reconstruction: key models and optimization algorithms}, IEEE Signal Processing Magazine, 37(1): 33–40, 2020.}

%\bibitem{GabMer:76}{ G. Gabay, and B. Mercier, \emph{A dual algorithm for the solution of nonlinear variational problems via finite element approximation}, Comput. Math. Appl. 2(1), 17–40, 1976.}

\bibitem{GhaLan:16}{ S. Ghadimi and G. Lan,  \emph{Accelerated gradient methods for nonconvex nonlinear and stochastic
programming},  Mathematical Programming, 156(1-2): 59-99, 2016. }

\bibitem{GloTal:89}{ R. Glowinski and P. Le Tallec, \emph{Augmented Lagrangian and Operator-Splitting Methods in Nonlinear Mechanics}, SIAM, 9, 1989.}

%\bibitem{GoyRoy:24} \red{F. Goyens ?????} and C.W. Royer, \emph{Riemannian trust-region methods for strict saddle functions with complexity guarantees},  Mathematical Programming, 2024. https://doi.org/10.1007/s10107-024-02156-2.


%\bibitem{GraBoy:14}{ M. Grant and S. Boyd, \emph{CVX: Matlab Software for Disciplined Convex Programming}, version 2.1, 2014. [Online]. Available:  \url{http://cvxr.com/cvx, Mar. 2014}}


\bibitem{GuoYe:18}
L. Guo and J. Ye, \emph{Necessary optimality conditions and exact penalization for non-Lipschitz nonlinear programs}, Mathematical Programming, 168: 571–598, 2018.

%\bibitem{HaaBuy:70}{ P. C. Haarhoff, J. D. Buys, \emph{A new method for the optimization of a nonlinear function subject to nonlinear constraints}, Computer Journal, 13: 178–184, 1970.}

\bibitem{HajHon:19}{ D. Hajinezhad and M. Hong, \emph{Perturbed proximal primal-dual algorithm for nonconvex nonsmooth optimization}, Mathematical Programming, 176(1-2): 207-245, 2019.}


\bibitem{HalTeb:23}
N. Hallak and M. Teboulle, \emph{An Adaptive Lagrangian-Based Scheme for Nonconvex Composite Optimization}, Mathematics of Operations Research, 48(4): 2337-2352, 2023.

%\bibitem{HanMan:79}{S. P. Han and O. L. Mangasarian, \emph{Exact penalty functions in nonlinear programming,} M. P., 17, No. 3, 251-269, 1979.}

\bibitem{Hes:69}{ M. Hestenes, \emph{Multiplier and gradient methods},  Journal of Optimization Theory and Applications, 4: 303–320, 1969.}

\bibitem{HonHaj:17}{ M. Hong,  D. Hajinezhad and  M. Zhao, \emph{ Prox-PDA: The proximal primal-dual algorithm for fast distributed nonconvex optimization and learning over networks}, Proceedings of International Conference on Machine Learning,  70: 1529–1538, 2017.}

\bibitem{JiaLin:19}{ B. Jiang, T. Lin, S. Ma and S. Zhang, \emph{Structured nonconvex and nonsmooth optimization: algorithms and iteration complexity analysis}, Computational Optimization and Applications, 72(1): 115–157, 2019.}

\bibitem{KruMeh:22} A.Y. Kruger and P. Mehlitz, \emph{Optimality conditions, approximate stationarity, and applications — a story
beyond Lipschitzness}, ESAIM: Control, Optimisation and Calculus of Variations, 28, 42, 2022.

\bibitem{LiChe:21} Z. Li, P. Chen, S. Liu, S. Lu and Y. Xu, \emph{Rate-improved inexact augmented Lagrangian method for constrained nonconvex optimization},  Proceedings of  International Conference on Artificial Intelligence and Statistics, 130:  2170-2178, 2021.


\bibitem{Lu:22} S. Lu, \emph{A single-loop gradient descent and perturbed ascent algorithm for nonconvex functional constrained optimization},  Proceedings of International Conference on Machine Learning, 2022. 

\bibitem{LukSab:19}{ D.R. Luke, S. Sabach and M. Teboulle, \emph{Optimization on spheres: models and proximal algorithms with computational performance comparisons}, SIAM Journal on Mathematics of Data Science, 1(3): 408–445, 2019.}

%\bibitem{MarOku:24} N. Marumo, T. Okuno and A. Takeda, \emph{Accelerated-gradient-based generalized Levenberg–Marquardt method with oracle complexity bound and local quadratic convergence}, Mathematical Programming, 2024. 

\bibitem{MesBau:21}{ F. Messerer, K. Baumgärtner and  M. Diehl,  \emph{Survey of sequential convex programming and generalized Gauss-Newton methods}, ESAIM: Proceedings and Surveys, 71:64-88, 2021.}


% \bibitem{Mis:23} K. Mishchenko,   \emph{Regularized Newton Method with Global O(1/k2) Convergence}, SIAM Journal on Optimization, 33(3): 1440-1462, 2023.

\bibitem{Mor:06}  B.S. Mordukhovich, \emph{Variational Analysis and Generalized Differentiation: Basic Theory},  Springer, Berlin, 330,  2006.


%\bibitem{MurSau:82}{ B.A. Murtagh and M.A. Saunders, \emph{A projected Lagrangian algorithm and its implementation for sparse nonlinear constraints}, Mathematical Programming, 16, 84–117, 1982.}

%\bibitem{NecKva:15}{ I. Necoara and S. Kvamme, \emph{DuQuad: A toolbox for solving convex quadratic programs using dual (augmented) first order algorithms}, 2015 54th IEEE Conference on Decision and Control (CDC), pp. 2043-2048, 2015.}

\bibitem{NedNec:14}{ V. Nedelcu, I. Necoara  and  Q. Tran-Dinh, \emph{Computational complexity of inexact gradient augmented Lagrangian methods: application to constrained MPC}, SIAM Journal on Control and Optimization, 52(5): 3109-3134, 2014.}

\bibitem{Nes:18}
Y. Nesterov,  \emph{Lectures on Convex Optimization}, Springer, Berlin, vol. 137, 2018.

\bibitem{Par:23}{
M. Parmar, \emph{Clustering Datasets}, GitHub, 2023,  \url{https://github.com/milaan9/Clustering-Datasets}.}

\bibitem{PatNec:17}{ A. Patrascu, I. Necoara and  Q. Tran-Dinh,
\emph{Adaptive inexact fast augmented Lagrangian methods for
constrained convex optimization},  Optimization Letters, 11(3): 
609-626, 2017.}

\bibitem{PenWei:07}{ 
J. Peng and Y. Wei,  \emph{Approximating K‐means‐type Clustering via Semidefinite Programming}, SIAM Journal on Optimization, 18(1): 186–205, 2007.}

%\bibitem{Pie:69}{ T. Pietzykowski, \emph{An exact potential method for constrained maxima,} SIAM J. Numer. Anal., No. 6, 269-304, 1969.}

\bibitem{Pow:69}{ M.J.D. Powell, \emph{A method for nonlinear optimization in minimization problems},  in R. Fletcher, Ed., Optimization, Academic Press, 283–298, 1969.}

\bibitem{RocWet:98}{ R. Rockafellar and R. Wets, \emph{Variational Analysis}, Springer, Berlin, 1998.}

%\bibitem{Roy:19}{ J.O. Royset, \emph{Variational Analysis in Modern Statistics}, Special Issue Mathematical Programming, Series B, 174, 2019.}

%bibitem{RoyOne:19}{ C. W. Royer, M. O’Neill and S. J.Wright,  \emph{A Newton-CG algorithm with complexity guarantees for smooth unconstrained optimization}, Mathematical Programming, 2019.}


\bibitem{SahEft:19}
M.F. Sahin, A. Eftekhari, A. Alacaoglu, F.L. Gomez and V. Cevher, \emph{ An Inexact Augmented Lagrangian Framework for Nonconvex Optimization with Nonlinear Constraints}, Proceedings of Neural Information Processing Systems, 13943–13955, 2019.

%\bibitem{SalVil:12} S. Salzo, S. Villa,  \emph{Convergence analysis of a proximal Gauss-Newton method}, Computational Optimization and Applications, 53: 557–589, 2012.

%\bibitem{SheTeb:14}{ R. Shefi and M. Teboulle, \emph{Rate of convergence analysis of decomposition methods based on the proximal method of multipliers for convex minimization}, SIAM Journal on Optimization, 24(1): 269–297, 2014.}

%\bibitem{TraDie:10}{ Q. Tran-Dinh and M. Diehl,  \emph{Local convergence of sequential convex programming for nonconvex optimization}. In M. Diehl, F. Glineur, E. Jarlebring, and W. Michiels, editors, \emph{Recent advances in optimization and its application in engineering}, pages 93–103. Springer-Verl., 2010.}

%\bibitem{WacBie:06}{ A. Wächter and L. T. Biegler, \emph{On the Implementation of a Primal-Dual Interior Point Filter Line Search Algorithm for Large-Scale Nonlinear Programming}, Mathematical Programming 106(1),  25-57, 2006.} 


\bibitem{XieWri:21}{ Y. Xie and S.J. Wright,  \emph{Complexity of proximal augmented Lagrangian for nonconvex optimization with nonlinear equality constraints},  Journal of Scientific Computing, 86, 2021.}

\bibitem{YanSun:15}{ L. Yang,  D. Sun and K.C. Toh, \emph{ SDPNAL+: a majorized semismooth Newton-CG augmented Lagrangian method for semidefinite programming with nonnegative constraints}, Mathematical Programming Computation, 7: 331–366, 2015.} 

%\bibitem{XieWri:21}{ Y. Xie, S. J. Wright,  \emph{Complexity of Proximal Augmented Lagrangian for Nonconvex Optimization with Nonlinear Equality Constraints},  Journal of Scientific Computing, 86, 38, 2021.}

%\bibitem{Yas:22}{ M. Yashtini, \emph{Convergence and rate analysis of a proximal linearized ADMM for nonconvex nonsmooth optimization}, Journal of Global Optimization, 84(4): 913-939, 2022.}

\bibitem{Zha:10}{ C. Zhang, \emph{Nearly unbiased variable selection under minimax concave penalty},  Annals of Statistics, 38: 894–942, 2010.}

\bibitem{ZhaLuo:20}{ J. Zhang and  Z.Q. Luo, \emph{A proximal alternating direction method of multiplier for linearly constrained nonconvex minimization}, SIAM Journal on Optimization, 30(3): 2272–2302, 2020.}

\end{thebibliography}
\end{document}